\documentclass[11pt]{article}

\usepackage{amsmath}
\usepackage{amssymb}
\usepackage{latexsym}
\usepackage{amsxtra}
\usepackage{amscd}
\usepackage{ntheorem}
\usepackage{xypic}
\usepackage{verbatim}
\usepackage{fullpage}
\usepackage[nottoc]{tocbibind}
\theoremstyle{plain}
\numberwithin{equation}{section}
{\theorembodyfont{\slshape}

        \newtheorem{thm}[equation]{Theorem}
        \newtheorem{cor}[equation]{Corollary}
        \newtheorem{lem}[equation]{Lemma}
        \newtheorem{prop}[equation]{Proposition}

}

{\theorembodyfont{\rmfamily}

        \newtheorem{defn}[equation]{Definition}
        \newtheorem{prob}[equation]{Problem}
        \newtheorem{rem}[equation]{Remark}
        
        \newtheorem{exa}[equation]{Example}
        \newtheorem{ass}[equation]{Assumption}
        \newtheorem{notation}[equation]{Notation}
        \newtheorem{claim}[equation]{Claim}

        \newtheorem{conj}[equation]{Conjecture}
}

\renewcommand{\em}{\sl}

\newcommand{\proof}{{\bf Proof:\ }}
\newcommand{\Endproof}{\hspace*{\fill} $\Box$ \vspace{1ex} \noindent }

\makeatletter
\renewcommand{\subsection}{\@startsection{subsection}{2}%
        {\z@}{-3.25ex plus -1ex minus-.2ex}{-1em}{\bf}}
\makeatother

\newcommand{\NN}{\mathbb{N}}
\newcommand{\ZZ}{\mathbb{Z}}
\newcommand{\QQ}{\mathbb{Q}}
\newcommand{\FF}{\mathbb{F}}

\newcommand{\PP}{\mathbb{P}}
\newcommand{\proj}{\mathbb{P}}
\newcommand{\KK}{\mathbb{K}}

\renewcommand{\AA}{\mathbb{A}}
\newcommand{\aff}{\mathbb{A}}
\newcommand{\BB}{\mathbb{B}}

\newcommand{\nats}{\mathbb{N}}
\newcommand{\ints}{\mathbb{Z}}
\newcommand{\rats}{\mathbb{Q}}
\newcommand{\reals}{\mathbb{R}}

\newcommand{\C}{\mathcal{C}}
\newcommand{\G}{\mathcal{G}}

\newcommand{\OO}{\mathcal{O}}

\newcommand{\K}{\mathfrak{K}}

\newcommand{\Aut}{{\rm Aut}}

\newcommand{\Gal}{{\rm Gal}}

\newcommand{\Tr}{{\rm Tr}}

\newcommand{\m}{\mathfrak{m}}

\newcommand{\Spec}{{\rm Spec\,}}
\newcommand{\MaxSpec}{{\rm MaxSpec\,}}

\newcommand{\Frac}{{\rm Frac}}

\newcommand{\crit}{{\rm crit}}
\newcommand{\hub}{{\rm hub}}
\newcommand{\crithub}{{\rm ch}}

\newcommand{\ord}{{\rm ord}}

\newcommand{\inj}{\hookrightarrow}
\newcommand{\To}{\;\longrightarrow\;}

\newcommand{\lpfeil}[1]{\stackrel{#1}{\To}}

\newcommand{\abs}[1]{\lvert#1\rvert}

\newcommand{\ol}[1]{\overline{#1}}
\newcommand{\mc}[1]{\mathcal{#1}}
\newcommand{\blockmatrix}[4]{
    \left( \begin{array}{c|c} #1 & #2 \\ \hline #3 &
    #4 \end{array} \right)}
\newcommand{\newmatrix}[4]{
    \left( \begin{array}{cc} #1 & #2 \\  #3 &
    #4 \end{array} \right)}

\newcommand{\Xb}{\bar{X}}
\newcommand{\Yb}{\bar{Y}}

\newcommand{\Wb}{\bar{W}}
\newcommand{\Db}{\bar{D}}

\newcommand{\Kb}{\bar{K}}

\newcommand{\xb}{\bar{x}}
\newcommand{\yb}{\bar{y}}

\newcommand{\alphab}{\bar{\alpha}}
\newcommand{\chib}{\bar{\chi}}

\newcommand{\psib}{\bar{\psi}}

\newcommand{\KKh}{\hat{\KK}}

\makeatletter
\newcommand{\subjclass}[2][2010]{%
  \let\@oldtitle\@title%
  \gdef\@title{\@oldtitle\footnotetext{#1 \emph{Mathematics subject classification.} #2}}%
}
\newcommand{\keywords}[1]{%
  \let\@@oldtitle\@title%
  \gdef\@title{\@@oldtitle\footnotetext{\emph{Key words and phrases.} #1.}}%
}
\makeatother
\begin{document}

\title{A generalization of the Oort conjecture}
\author{Andrew Obus\thanks{The author was supported by an NSF
Mathematical Science Postdoctoral Research Fellowship, as well as NSF
FRG Grant DMS-1265290.}}
\maketitle

\begin{abstract}
The Oort conjecture (now a theorem of Obus-Wewers and Pop) states that
if $k$ is an algebraically closed field of characteristic $p$, then
any cyclic branched cover of smooth projective $k$-curves lifts to
characteristic zero.  This is equivalent to the local Oort conjecture,
which states that all cyclic extensions of $k[[t]]$ lift to characteristic zero.
We generalize the local Oort conjecture to the case of Galois
extensions with cyclic $p$-Sylow subgroups, reduce the conjecture to a
pure characteristic $p$ statement, and prove it in several cases.  In
particular, we show that $D_9$ is a so-called \emph{local Oort group}.
%
\end{abstract}

\tableofcontents
\section{Introduction}\label{Sintro}

This paper concerns the \emph{local lifting problem} about lifting
Galois extensions of power series rings from characteristic $p$ to
characteristic zero.  In particular, in Conjecture \ref{Cmain} we state a generalization of the
\emph{Oort conjecture} on lifting of cyclic extensions, now a theorem
of Obus-Wewers (\cite{OW:ce}) and Pop (\cite{Po:oc}).  Our main result reduces the
generalized conjecture to an easy-to-understand, pure characteristic $p$ assertion about
existence of certain meromorphic differential forms on $\proj^1$.
We prove this assertion in several cases,
exhibiting the first positive cases of the local lifting problem for a
nonabelian group with cyclic $p$-Sylow subgroup of order greater than
$p$.  In particular, we show that $D_9$ is a so-called \emph{local Oort
  group}, and we completely solve the ``inverse Galois problem'' for
the local lifting problem for groups with cyclic $p$-Sylow subgroups.  See \S\ref{Snew} for specifics.

\subsection{The local lifting problem}\label{Sllp}
For our purposes, a finite extension $B/A$ of rings is called
$\Gamma$-\emph{Galois} (or a $\Gamma$-\emph{extension}) 
if $A$ and $B$ are integrally closed integral domains and
$\Frac(B)/\Frac(A)$ is $\Gamma$-Galois.  

\begin{prob}[The local lifting problem]\label{llp}
  Let $k$ be an algebraically closed field of characteristic $p$ and $\Gamma$ a finite group.  
  Let $k[[z]]/k[[s]]$ be a $\Gamma$-Galois extension.  Does this
  extension lift to characteristic zero?  That is,
 does there exist a DVR $R$ of characteristic zero with residue field $k$ and a
  $\Gamma$-Galois extension $R[[Z]]/R[[S]]$ that reduces to
  $k[[z]]/k[[s]]$?  In other words, does the $\Gamma$-action on $R[[Z]]$ reduce to that on $k[[z]]$,
  if we assume that $Z$ reduces to $z$?
\end{prob}
We will refer to a $\Gamma$-Galois extension $k[[z]]/k[[s]]$ as a
\emph{local $\Gamma$-extension}.

\begin{rem}\label{ringfieldrem}
Suppose $k$ is an algebraically closed field, and $B/A$ is any $\Gamma$-Galois extension of $k$-algebras with
the Galois group acting by $k$-automorphisms.  Then, if \emph{either}
$B$ or $A$ is isomorphic to a power series ring in one variable over $k$, the
other is as well.  That is, $B/A$ is a local $\Gamma$-extension.
\end{rem} 

\begin{rem}
Basic ramification theory shows that any group $\Gamma$ that occurs
as the Galois group of a local extension is of the form
$P \rtimes \ints/m$, with $P$ a $p$-group and $p \nmid m$.
\end{rem}

The main motivation for the local lifting problem is the following
\emph{global lifting problem}, about deformation of curves with an action of
a finite group (or equivalently, deformation of Galois branched covers
of curves).

\begin{prob}[The global lifting problem]\label{glp}
Let $X/k$ be a smooth, connected, projective curve over an
algebraically closed field of characteristic $p$.  Suppose a finite
group $\Gamma$ acts on $X$.  Does $(X, \Gamma)$ lift to characteristic
zero?  That is, does there exist a DVR $R$ of characteristic
zero with residue field $k$ and a relative projective curve $X_R/R$
with $\Gamma$-action such that $X_R$, along with its $\Gamma$-action,
reduces to $X$?
\end{prob}

It is a major result of Grothendieck (\cite[XIII, Corollaire
2.12]{SGA1}) that the global lifting problem can be solved whenever $\Gamma$
acts with tame (prime-to-$p$) inertia groups, and $R$ can be taken to
be the Witt ring $W(k)$.  In particular, it holds when $\Gamma$ is trivial.  The wild case is much more subtle, and cannot always
be solved.  For instance, the group $\ints/p \times \ints/p$ acts
faithfully on $\proj^1_k$ whenever $k$ is algebraically closed of
characteristic $p$, but there can be no lifting of this action to a
genus zero curve when $p$ is odd.  However, the \emph{local-global
  principle} states that the global lifting problem holds for $(X,
\Gamma)$ (and a complete DVR
$R$) if and only if the
local lifting problem holds (over $R$) for each point of $X$ with nontrivial
stabilizer in $\Gamma$.  Specifically, if $x$ is such a point, then its
complete local ring is isomorphic to $k[[z]]$. The stabilizer $I_x
\subseteq \Gamma$ acts on $k[[z]]$ by $k$-automorphisms, and we check the local lifting
problem for the local $I_x$-extension $k[[z]]/k[[z]]^{I_x}$.  Thus,
the global lifting problem is reduced to the local lifting problem.

A proof of the local-global principle
for abelian $\Gamma$ is already implicit in \cite{SOS}.
Proofs for arbitrary $\Gamma$ have been given by
Bertin and M\'{e}zard (\cite{BM:df}), Green and Matignon (\cite{GM:lg}), and Garuti (\cite{Ga:pr}).

The author's paper \cite{Ob:ll} is a detailed exposition of many aspects of the local lifting problem.

\subsection{Local Oort groups and the KGB obstruction}

The \emph{Oort conjecture} (as mentioned above, now a theorem), states
that the local lifting problem holds for all \emph{cyclic}
extensions.  In \cite{CGH:ll}, Chinburg, Guralnick, and Harbater
ask which finite groups $\Gamma$ of the form $P \rtimes \ints/m$,
with $P$ a $p$-group and $p \nmid m$, have this same property.
That is, given a prime $p$, for which groups $\Gamma$ is it true that
all local $\Gamma$-actions (over
all algebraically closed fields of characteristic $p$) lift to
characteristic zero?
Such a group is called a \emph{local Oort group} (for $p$).  The paper
\cite{CGH:ll} also investigates the
notion of a \emph{weak local Oort group} (for $p$), which is a group $\Gamma$ for which
there exists \emph{at least one} local $\Gamma$-extension that lifts
to characteristic zero.  Characterizing the weak local Oort groups has
been referred to as the ``inverse Galois problem'' for the local
lifting problem in \cite{Ma:pg}, where it was proved that elementary
abelian $p$-groups are weak local Oort.

In \cite{CGH:ll}, Chinburg, Guralnick, and Harbater introduce the so-called \emph{KGB
  obstruction} to local lifting (this is related to the earlier
\emph{Bertin obstruction} from \cite{Be:ol}).  Roughly, given a local $\Gamma$-extension,
the KGB obstruction vanishes if there is a $\Gamma$-extension
of certain characteristic zero power series rings for which the
different behaves in the same way as for the original local
$\Gamma$-extension.  A lift to characteristic zero gives such an
extension, and thus causes the KGB obstruction to vanish.  Using this obstruction,
Chinburg, Guralnick, and Harbater were able to greatly restrict the possible local Oort groups.

\begin{thm}[\cite{CGH:ll}, Theorem 1.2]\label{Tlocaloort}
If a group $\Gamma$ is a local Oort group for $p$, then $\Gamma$ is
either cyclic, dihedral of order $2p^n$, the alternating group $A_4$
($p=2$), or a generalized quaternion group ($p=2$).
\end{thm}
In fact, Brewis and Wewers showed (\cite{BW:ac}) that the generalized
quaternion groups are not local Oort, so the list of possible local
Oort groups consists only of the cyclic groups, $D_{p^n}$, and $A_4$.  The cyclic case is the Oort conjecture, and the $A_4$ case
has been claimed by Bouw (see \cite{BW:ll}) and written up by the
author (\cite{Ob:A4}).  Thus only the $D_{p^n}$ have unknown ``local Oort status,'' and
showing that the local lifting problem holds for these dihedral groups
has been referred to as the ``strong Oort conjecture'' (\cite{CGH:og}).   We propose a
somewhat different generalization (Conjecture \ref{Cmain}) below.

\subsubsection{Cyclic $p$-Sylow groups}

Recall that a $\Gamma$-extension $L_n/k[[s]]$ 
gives rise to a higher ramification filtration $\Gamma^i_{i \geq 0}$ for the 
upper numbering on the group $\Gamma$ (\cite{Se:lf}, IV).  If $\Gamma
= \ints/p^n$, then the breaks in this filtration (i.e., the values $i$ for which $\Gamma^i \supsetneq \Gamma^j$ for all $j > i$) will
be denoted by $(u_1, u_2, \ldots, u_n)$.  One knows that $u_i \in \NN$ and
$$u_i\geq pu_{i-1},$$
for $i=2,\ldots,n$ (see, e.g., \cite{Ga:ls}). 

The higher ramification filtration gives us all the data we need to
check the KGB obstruction in the case where $\Gamma$ has a cyclic
$p$-Sylow subgroup.

\begin{prop}[\cite{Ob:ll}, Proposition 5.9]\label{PKGB}
Let $\Gamma$ be a semi-direct product of the form $\ints/p^n \rtimes \ints/m$, with $p \nmid m$.  Suppose $\Gamma$ is not cyclic (thus not abelian).  
Let $L_n/k[[s]]$ be a $\Gamma$-extension whose $\ints/p^n$-subextension has upper ramification breaks $(u_1, \ldots, u_n)$. 
Then the KGB obstruction vanishes for $L_n/k[[s]]$ if and only if $u_1 \equiv -1 \pmod{m}$.
\end{prop}

\begin{rem}\label{RKGBcons}
By \cite[Theorem 1.1]{OP:wt}, knowing that $u_1 \equiv -1 \pmod{m}$
and $\Gamma$ is non-abelian implies that $\Gamma$ is center-free (in
particular, $m | (p-1)$) and $u_i \equiv -1 \pmod{m}$ for all
$i$.
\end{rem}

\begin{rem}
One can also phrase the KGB obstruction in terms of the higher
ramification breaks for the lower numbering.  In this case, the
criterion for vanishing is the same --- that the first break for the
lower numbering is
congruent to $-1 \pmod{m}$.
\end{rem}

Our generalization of the Oort conjecture is the following:
\begin{conj}\label{Cmain}
For local $\Gamma$-extensions where $\Gamma$ has a cyclic $p$-Sylow
subgroup (that is, $\Gamma$ is of the form $\ints/p^n \rtimes \ints/m$), the KGB obstruction is the only obstruction to lifting.
\end{conj}

\begin{rem}
Note that, if $\Gamma = D_{p^n}$ with
$p$ an \emph{odd} prime, then all $u_i$ as above are odd
(see, e.g., \cite[Theorem 1.1]{OP:wt}).  Thus Conjecture \ref{Cmain}
(combined with Proposition \ref{PKGB})
implies that $\Gamma$ is a local Oort group for $p$.  So for $p$ odd,
Conjecture \ref{Cmain} is somewhat stronger than the ``strong Oort
conjecture'' mentioned above.  However, Conjecture \ref{Cmain} says
nothing about $D_{2^n}$.  We currently have no opinion as to whether $D_{2^n}$
is a local Oort group.  The only results toward this end are that
$D_2 \cong \ints/2 \times \ints/2$ is a local Oort group
(\cite{Pagotthesis}, or \cite{Pa:ev} for a special case)
and $D_4$ is a weak local Oort group (\cite{Br:D4}).
\end{rem}

\subsection{Some history leading to Conjecture \ref{Cmain}}\label{Sknown}

The first major result on the local lifting problem was the 1989 paper
\cite{SOS} of Sekiguchi-Oort-Suwa, which showed that $\ints/pm$ is a local Oort group when $p
\nmid m$.  That $\ints/p^2m$ is a local Oort group was proven in 1998
by Green-Matignon (\cite{GM:lg}).  The full Oort conjecture was proven
in 2014 by Obus-Wewers and Pop (\cite{OW:ce} and \cite{Po:oc}).

The local lifting problem for $\Gamma$ is much more difficult when $\Gamma$ is
non-abelian, even if we assume its $p$-Sylow subgroup is cyclic.  
Indeed, it was not until a 2006 paper that $D_p$ (for odd $p$) was proved to be
local Oort by Bouw-Wewers (\cite{BW:ll}), and this proof is
significantly more intricate than the $\ints/p$ case.  In fact, it was
proven in the two papers \cite{BW:ll} and \cite{BWZ:dd} that
Conjecture \ref{Cmain} holds when $p$ \emph{exactly divides the order
  of $\Gamma$}.  However, other than this, up until this paper, there was essentially nothing known when $\Gamma$ is non-abelian.  
Namely, if the $p$-Sylow subgroup of $\Gamma$ is cyclic of order greater
than $p$ and $\Gamma$ is non-abelian, 
then there was no local $\Gamma$-extension with vanishing KGB obstruction that was known
\emph{either to lift or not to lift} to characteristic zero.  In particular, it was not known if such $\Gamma$ were weak local Oort groups.  We show that 
they in fact are (Corollary \ref{Cminimaljumps}).  Furthermore, our
main result (Theorem \ref{Tmachine}) brings the full solution to the
local lifting problem for such $\Gamma$ within reach.

\subsection{The (isolated) differential data criterion}\label{Sdd}
While we are not yet able to present a full proof of Conjecture \ref{Cmain}, we are able to prove it conditionally on certain meromorphic
differential forms on $\proj^1_k$ existing with special properties.  We describe this condition briefly now (for more details, see \S\ref{Scrit}).

Consider quadruples $(p, m, \tilde{u}, N_1)$ of natural numbers where:
\begin{itemize}
\item $p$ is a prime number.
\item $m > 1$ divides $p-1$.
\item $\tilde{u} \equiv -1 \pmod{m}$.
\item $N_1$ is divisible by $m$.
\end{itemize}
Write $\tilde{u} = up^{\nu}$ with $u$ prime to $p$. 
Let $k$ be an algebraically closed field of characteristic $p$. 
We say that $(p, m, \tilde{u}, N_1)$ \emph{satisfies the differential data criterion} (with respect to $k$) 
if there exists a polynomial $f(t) \in k[t^m]$ of degree exactly $N_1$
in $t$, such that the meromorphic differential form
$$\omega := \frac{dt}{f(t)t^{\tilde{u}+1}} \in \Omega^1_{k(t)/k}$$ satisfies
$$\C(\omega) = \omega + u t^{-\tilde{u}-1}dt.$$  Here $\C$ is the Cartier operator on differential forms.  Note that $\omega$ has a zero of order
$N_1 + \tilde{u} -1$ at $t = \infty$.

If $(p, m, \tilde{u}, N_1)$ satisfies the 
differential data criterion with notation as above, the basic properties of the Cartier operator imply that 
$$\omega = dg/g - u \sum_{i=0}^{\nu} t^{-up^i-1} dt,$$ for some $g \in k(t)$, well-defined up to multiplication by $p$th powers.  We say that
$(p, m, \tilde{u}, N_1)$ satisfying the differential data criterion 
\emph{satisfies the isolated differential data criterion} if there are
$f$ and $\omega$ as above such that no infinitesimal deformation $\tilde{g}$ of $g$ gives rise to a differential
form $\tilde{\omega} := d \tilde{g}/\tilde{g} - u \sum_{i=0}^{\nu} t^{-up^i-1} dt$ having a zero of order at least $N_1 + \tilde{u} - 1$ at $t = \infty$
(as will be seen in \S\ref{Scrit}, this is equivalent to invertibility of a ``Vandermonde-like" matrix constructed from the roots of $f$).
This is readily seen to be independent of the choice of $g$, once $f$
is chosen.

\subsection{Main results}\label{Snew}

Throughout this section, $m \in \nats$ is not divisible by $p$.
First we adapt an argument of Pop (\cite{Po:oc}) to reduce Conjecture
\ref{Cmain} to the case where the successive upper jumps do not grow too quickly.

\begin{prop}\label{Ppopreduction}
Let $L_n'/k[[s]]$ be a non-abelian $\Gamma = \ints/p^n \rtimes \ints/m$-extension whose $\ints/p^n$-subextension has upper ramification breaks
$(u_1', \ldots, u_n')$.  
For $1 \leq i \leq n$, 
define $u_i$ inductively to be the unique integer such that $u_i \equiv u_i' \pmod{mp}$ and $pu_{i-1} \leq u_i < pu_{i-1} + mp$ (by convention, set
$u_0 = 0$).
If, for every algebraically closed field $\kappa$ of characteristic $p$, every $\Gamma$-extension $L_n/\kappa[[s]]$ whose 
$\ints/p^n$-subextension has upper ramification breaks $(u_1, \ldots, u_n)$ lifts to characteristic zero, then so does $L_n'/k[[s]]$.
\end{prop}

Thus we need only consider $\Gamma$-extensions whose upper ramification breaks satisfy $u_i < pu_{i-1} + mp$. We say that these 
extensions have \emph{no essential ramification}.  

\begin{exa}\label{Example}
For instance, if we have a $\ints/5^4 \rtimes \ints/2$-extension with
$(u_1', \ldots, u_4') = (11, 79, 433, 2165)$, then we would have
$(u_1, \ldots, u_4) = (1, 9, 53, 265)$.
\end{exa} 

\begin{rem}
Note the similarity between this definition and
\cite[Remark/Definition 3.1(2)]{Po:oc}.  In fact, if we consider the upper jumps for the entire 
$G$-extension, as opposed to just the $\ZZ/p^n$-part, then our
assumption is exactly that of ``no essential ramification'' from
\cite{Po:oc}.  Indeed, Proposition \ref{Ppopreduction} in the abelian
case is equivalent to the main result of \cite{Po:oc}.
\end{rem}

Our main result is the following:
\begin{thm}\label{Tmachine}
Let $L_n/k[[s]]$ be a non-abelian $\ints/p^n \rtimes \ints/m$-extension whose $\ints/p^n$-subextension has upper ramification breaks
$(u_1, \ldots, u_n)$.  Suppose that $L_n/k[[s]]$ has vanishing KGB obstruction and no essential ramification. 
Suppose further that for all $1 < i \leq n$, the quadruple $(p, m, u_{i-1}, N_{i,1})$ satisfies the isolated differential data criterion, where 
$N_{i,1} = (p-1)u_{i-1}$ if $u_i = pu_{i-1}$ and $N_{i,1} = (p-1)u_{i-1} - m$ otherwise.  Then the extension $L_n/k[[s]]$ lifts to characteristic zero.
\end{thm}   

\begin{rem}
Our lifts correspond to certain covers of the non-archimedian open
disk.  We discuss the geometry of the branch locus of these covers 
in \S\ref{Sbranchgeometry}.
\end{rem}

\begin{rem}
Proposition \ref{Ppopreduction} and Theorem \ref{Tmachine} reduce Conjecture \ref{Cmain} for the group $\ints/p^n \rtimes \ints/m$ 
(nonabelian) to realizing the isolated differential data criterion for quadruples $(p, m, \tilde{u}, (p-1)\tilde{u})$ and 
$(p, m, \tilde{u}, (p-1)\tilde{u} - m)$ such that $\tilde{u} \equiv -1 \pmod{m}$, that $p^{n-1} \nmid \tilde{u}$, and that 
$\tilde{u} < m(p^{n-1} + p^{n-2} + \cdots + p)$.   
Thus, once the group is fixed, one need only realize the isolated
differential data criterion for finitely many quadruples. Our proof of
Corollary \ref{CD9} below proceeds by this method.

If one believes, for a particular group $\Gamma = \ints/p^n \rtimes \ints/m$, that there is a particular finite field $\FF_q$ such
that the isolated differential data criterion in the above cases
can always be realized using a polynomial $f(t) \in \FF_q[t]$, then
proving Conjecture \ref{Cmain} for $\Gamma$ is reduced to a finite search.
\end{rem}

\begin{exa}
In order to show that all extensions as in Example \ref{Example} lift
to characteristic $0$, we would have to realize the isolated
differential data criterion for $(5, 2, 1, 2)$, $(5, 2, 9, 34)$, and $(5,
2, 53, 212)$.
\end{exa}

By realizing various instances of the isolated differential data
criterion, we are able to prove the following corollaries, which are
special cases of Conjecture \ref{Cmain}.

\begin{cor}[Theorem \ref{TD9}]\label{CD9}
The dihedral group $D_9$ is a local Oort group for $p = 3$.
\end{cor}

\begin{cor}[Theorem \ref{Tfirstjump1}] \label{Cfirstjump1}
If $p$ is an odd prime, and $L/k[[s]]$ is a $D_{p^2}$-extension whose $\ints/p^2$-subextension has first upper ramification break
$u_1 \equiv 1 \pmod{p}$, then $L/k[[s]]$ lifts to characteristic zero.
\end{cor}

\begin{cor}[Theorem \ref{Tminimaljumps}]\label{Cminimaljumps}
If $L/k[[s]]$ is a $\ints/p^n \rtimes \ints/m$-extension whose $\ints/p^n$-subextension has upper ramification breaks
congruent to $(m-1, p(m-1), \ldots, p^{n-1}(m-1))$ $\pmod{mp}$, then $L/k[[s]]$ lifts to characteristic zero. In particular, $\ints/p^n \rtimes \ints/m$ is
a weak local Oort group whenever the conjugation action of $\ints/m$ on $\ints/p^n$ is faithful.
\end{cor}

\begin{rem}
For each non-abelian $\ints/p^n \rtimes \ints/m$, Corollary
\ref{Cminimaljumps} includes the case with the smallest possible ramification breaks
causing the KGB obstruction to vanish (these breaks are in fact $(m-1, p(m-1),
\ldots, p^{n-1}(m-1))$).
\end{rem}

\begin{rem}
By Proposition \ref{PKGB} and Remark \ref{RKGBcons}, the action of
$\ints/m$ on $\ints/p^n$ must be faithful for $\Gamma$ to be a weak
local Oort group (unless $\Gamma$ is cyclic).  Corollary
\ref{Cminimaljumps} says that this condition suffices as well, and
thus solves the ``inverse Galois problem'' for the local lifting problem for groups with cyclic
$p$-Sylow subgroups.
\end{rem}

\begin{rem}\label{Ruseless}
The proof of Theorem \ref{Tmachine} follows the same basic outline as the
analogous assertion for cyclic groups in
\cite{OW:ce}.  However, we never invoke the Oort conjecture itself in
the proof.  To
emphasize this point, note that any lift of a local non-abelian $\Gamma := \ints/p^n \rtimes \ints/m$-extension
necessarily yields an ``equivariant''  lift of its unique local
$\ints/p^n$-subextension (see \S\ref{Sequivariant}).  However,
\emph{none} of the cyclic lifts constructed in \cite{OW:ce} are
equivariant, so they cannot possibly occur inside a lift of a local $\Gamma$-extension. Thus the lifts from
\cite{OW:ce} are ``useless'' for constructing non-abelian lifts as in
Theorem \ref{Tmachine}.  
\end{rem}

\subsection{Outline of the paper} 
In \S\ref{Switt}, we recall the explicit parameterization of local
$\ints/p^n \rtimes \ints/m$-extensions, and the relationship between
the parameterization and the
higher ramification filtration.  In \S\ref{Spop}, we prove Proposition
\ref{Ppopreduction}, which allows us to consider only extensions with
no essential ramification.  Then, \S\ref{Sind}--\S\ref{Sproof} are devoted
to the proof of Theorem \ref{Tmachine}.  In \S\ref{Sind}, we set up
the induction on $n$ that we will use (which is essentially the same
framework used in \cite{OW:ce}), and in \S\ref{Sbase}, we prove
the base case $n=1$.  In \S\ref{Sgeom}, we recall the
language of characters that was used in \cite{OW:ce}, and adapt it to
our new situation of non-abelian groups.  The main part of the proof
is in \S\ref{Sproof}, and we give a further, more detailed outline in
\S\ref{Sstrategy}.  We remark that, although the basic idea of the
proof is the same as in \cite{OW:ce}, the execution is quite different
and more complicated.  To enhance the flow of the paper and clarify the
main argument, we postpone the proofs of two particularly technical
results to \S\ref{Sproofs}.

In \S\ref{Sexamples}, we give some examples of when the isolated
differential data criterion (\S\ref{Sdd}) is realized, and derive
consequences for the local lifting problem. 

\subsection{Conventions}
The letter $K$ will always be a field of characteristic zero that is complete with respect to a
discrete valuation $v:K^\times\to\QQ$. We assume that the residue field $k$ of
$K$ is algebraically closed of characteristic some fixed \emph{odd} prime $p$. We also assume that the
valuation $v$ is normalized such that $v(p)=1$.  We let $|\cdot|$ be
an absolute value on $K$ corresponding to $v$ (it does not matter how
it is normalized).  The ring of
integers of $K$ will be denoted $R$.  The maximal ideal of $R$ will be
denoted $\m$. The notation $R\{T\}$ refers to the ring of power series 
$\sum_{i=0}^{\infty}c_iT^i$ such that $\lim_{i \to \infty} |c_i| =
0$.  We write $\m\{T\}$ to refer to the subset of $R\{T\}$ for which all
$c_i$ lie in $\m$.

We fix an algebraic closure $\Kb$ of $K$, and whenever necessary, we will replace $K$ by a suitable finite
extension within $\Kb$, without changing the above notation.  
Furthermore, we fix 
once and for all a compatible system of elements $p^r \in \Kb$ for $r \in \rats$, such that $p^{r_1}p^{r_2} = p^{r_1 + r_2}$.
The letter $m$ will always refer to a prime-to-$p$ integer.  The symbol $\zeta_n$ denotes a 
primitive $n$th root of unity.  A curve is always (geometrically)
connected.  

These are the same conventions used in \cite{OW:ce}.

\section*{Acknowledgements}
I thank Ted Chinburg, Johan de Jong, Bob Guralnick, David Harbater, and
Florian Pop for useful conversations.  I especially thank Stefan
Wewers and Irene Bouw, not only for useful conversations, but also for
providing hospitality in Ulm when some of this work was done.
Some of the computations were done in SAGE, and I thank Julian R\"{u}th
for assistance.  Lastly, I thank the referees for helpful expository improvements.

\section{$\ints/p^n \rtimes \ints/m$-extensions in characteristic $p$}\label{Switt}

In this section, we recall the cyclic theory of local extensions in characteristic $p$, and then show how to adapt it to the metacyclic
case considered in this paper.  Let $\Gamma = \ints/p^n \rtimes \ints/m$.

If $L/k[[s]]$ is a $\Gamma$-extension, then, after a possible change of variables, we may assume that
the subextension corresponding to the normal subgroup $\ints/p^n \subseteq \Gamma$ can be written as $k[[t]]/k[[s]]$, with $t^m = s$.
Let $M = \Frac(L)$.  Since $\Gal(M/k((t))) \cong \ints/p^n$, Artin-Schreier-Witt theory states that $M/k((t))$ is given by
an Artin-Schreier-Witt equation 
\[
       \wp(y_1,\ldots,y_n)=(f_1,\ldots,f_n),
\]
where $(f_1, \ldots, f_n)$ lies in the ring  $W_n(k((t)))$ of truncated Witt vectors, 
$F$ is the Frobenius morphism on $W_n(k((t)))$, and $\wp(y):=F(y)-y$ is the Artin-Schreier-Witt 
isogeny.  Then $L$ is the integral closure of $k[[t]]$ in $M$.
Adding a truncated Witt vector of the form $\wp((g_1, \ldots, g_n))$ to $(f_1,\ldots,f_n)$ does not change the extension, 
and adjusting by such Witt vectors, we may assume that the 
$f_i$ are polynomials in $t^{-1},$ all of whose terms have prime-to-$p$ degree (in this case, we say the Witt vector is in \emph{standard form}).  
If
\begin{equation}\label{rambreaks}
     u_i :=\max\{\, p^{i-j}\deg_{t^{-1}}(f_j) \mid j=1,\ldots,i\,\},
\end{equation}
then the $u_i$ are exactly the breaks in the higher ramification filtration of $M/k((t))$ (\cite{Ga:ls}, Theorem 1.1).
From this, one sees that $p \nmid u_1$, that $u_i \geq pu_{i-1}$ for $2 \leq i \leq n$, and that if $p | u_i$, then $u_i = pu_{i-1}$.

\begin{prop}\label{Ppries}
The extension $L/k[[s]]$ is $\Gamma$-Galois if and only if the degrees (in $t^{-1}$) of all terms appearing in 
the polynomials $f_i$ are in the same congruence class$\pmod{m}$. 
\end{prop}  

\proof This follows from \cite[Proposition 4.3]{OP:wt}. \Endproof

Thus we can, and will think of $\Gamma$-Galois extensions $L/k[[s]]$ as corresponding to Witt vectors
$(f_1, \ldots, f_n) \in W_n(k((t)))$ such that the $f_i$ are polynomials in $t^{-1}$ with all degrees of all terms of all $f_i$ congruent to each 
other$\pmod{m}$.  By (\ref{rambreaks}), this implies that all $u_i$ belong to this congruence class. 

Recall Proposition \ref{PKGB}, which states that, for non-abelian $\Gamma$, the KGB obstruction vanishes for $L/k[[s]]$
if and only if $u_1 \equiv -1 \pmod{m}$.  By Remark \ref{RKGBcons},
this is true for all $i$, and this implies that $\Gamma$ is center-free.  For the rest of the paper, we only consider local 
$\Gamma$-extensions of this form.

\section{Reduction to the ``no essential ramification'' case}\label{Spop}

In this section, we prove Proposition \ref{Ppopreduction}.  Recall that $L_n'/k[[s]]$ is a 
non-abelian $\Gamma = \ints/p^n \rtimes \ints/m$-extension whose $\ints/p^n$-subextension has upper ramification breaks
$(u_1', \ldots, u_n')$, and $u_i$ is defined
inductively to be the unique integer such that $u_i \equiv u_i' \pmod{mp}$,
with $u_1 < mp$ and $pu_{i-1} \leq u_i < pu_{i-1} + mp$ for $i > 1$.
We may, and do, assume that $L_n'/k[[s]]$ has vanishing KGB obstruction, i.e., that all $u_i$ and $u_i'$ are $-1 \pmod{m}$.  
Write $L_n' = k[[z]]$, and write $M = k[[t]] \subseteq L_n'$, where $t^m = s$, so that $M$ is the subextension of $L_n'/k[[s]]$ corresponding to 
$\ints/p^n \subseteq \Gamma$.
Our proof follows \cite{Po:oc}.  The key is to make a deformation in characteristic $p$ so that the generic fiber has no non-abelian essential ramification, 
in some sense (cf.\ \cite[Key Lemma 3.2]{Po:oc}).

\begin{prop}[Generalized characteristic $p$ Oort conjecture]\label{Pcharpoort}
Let $\mc{A} = k[[\varpi, s]] \supseteq k[[s]]$, and let $\mc{K} = \Frac(\mc{A})$.  
There exists a $\Gamma$-extension $\mc{L}/\mc{K}$, with $\mc{L} \supseteq L_n'$, having the following properties:

\begin{enumerate}
\item The $\ints/m$-subextension $\mc{M}/\mc{K}$ 
corresponding to the subgroup $\ZZ/p^n \subseteq \Gamma$ is given by $\mc{M} = \mc{K}[t] \subseteq \mc{L}$.
\item If $\mc{B}$ is the integral closure of $\mc{A}$ in $\mc{L}$, we have $\mc{B} \cong k[[\varpi, z]]$. 
In particular, $(\mc{B}/(\varpi))/(\mc{A}/(\varpi))$ is $\Gamma$-isomorphic to the original extension $L_n'/k[[s]]$.
\item Let $\mc{C} = \mc{A}[t] \subseteq \mc{M}$.  Let $\mc{R} = \mc{A}[\varpi^{-1}]$, let $\mc{S} = \mc{B}[\varpi^{-1}]$, 
and let $\mc{T} = \mc{C}[\varpi^{-1}]$.  
Then $\mc{S}/\mc{T}$ is a $\ZZ/p^n$-extension of Dedekind rings, branched at $m+1$ maximal ideals.  Above the ideal
$(t)$, the inertia group is $\ZZ/p^n$, and the upper jumps are $(u_1, \ldots, u_n)$.  The other $m$ branched ideals are of the form
$(\zeta_m^{\alpha} t - \mu)$, where $\mu$ can be chosen arbitrarily in $\varpi^{p^{\delta_0}} k[[\varpi^{p^{\delta_0}}]]$ for some high enough $\delta_0$,
and $\alpha$ ranges from $1$ to $m$.  
\item The only branched ideal of $\mc{S}/\mc{R}$ with noncyclic inertia group is $(s)$.  
\end{enumerate}
\end{prop}

\proof
As in \cite{Po:oc}, we will prove Proposition \ref{Pcharpoort} by deforming a standard form classifying Witt vector $(f_1, \ldots, f_n)$ 
of the extension $L_n'/k[[s]]$.  
We must take care to do everything equivariantly.  By Proposition \ref{Ppries}, each $f_i$ can be written as $t^{1-m}g_i(t^{-m}),$ where $g_i$ is a
polynomial of degree $\leq (u_i' + 1 - m)/m$ over $k$. Equality holds if $p \nmid u_i'$.
Choose a factoring
$$f_i = t^{1-m}p_i(t^{-m})q_i(t^{-m}),$$
where $$\deg p_i \leq \frac{u_i + 1 - m}{m}$$ and $$\deg q_i \leq \frac{u_i' - u_i}{m}.$$  If $u_i = u_i'$, then take $q_i = 1$.
Note that if $p \nmid u_i$, then we must have equality in both inequalities above.  Factoring, we can write
$$q_i(t^{-m}) = c\prod_{\alpha = 1}^m j_i(\zeta_m^{\alpha}t^{-1}),$$ where $c \in k$ and the $j_i$ are monic polynomials of degree $\deg(q_i)$.  
Lastly, factor $j_i$ completely to write
$$j_i(t^{-1}) = \prod_{\nu = 1}^{\deg q_i} (t^{-1} - r_{\nu, i}).$$   

Now, let $\mu \in \varpi k[[\varpi]] \backslash \{0\}$.  We lift the Witt vector $(f_1, \ldots, f_n) \in W_n(k((t)))$ to a Witt vector 
$(F_1, \ldots, F_n) \in W_n(\mc{M})$.  We choose
$$F_i = ct^{1-m}p_i(t^{-m}) \prod_{\alpha = 1}^m \prod_{\nu = 1}^{\deg q_i} ((\zeta_m^{\alpha}t - \mu)^{-1} - r_{\nu, i}).$$

Let us make some observations:  
\begin{itemize}
\item $F_i$ (viewed as an element of $\mc{A}_{(\varpi)}$) reduces to $f_i$ modulo $\varpi$.
\item All terms in $F_i$ are of degree $-1 \pmod{m}$ in $t^{-1}$.
\item $F_i$ has a pole of order $\leq u_i$ at $t=0$, and, 
for each $\alpha \in \{1, \ldots, m\}$, a pole of order $\leq (u_i' - u_i)/m$ at $t = \zeta_m^{-\alpha}\mu$.
\end{itemize}

Let $\mc{L}/\mc{M}$ be the $\ints/p^n$-extension classified by $(F_1,
\ldots, F_n)$.  By the second observation above and the discussion in
\S\ref{Switt}, this extends to a $\Gamma$-extension $\mc{L}/\mc{K}$, which
will be the extension we seek.  In order to prove this, we must show that the degree $\delta_{\mc{S}/\mc{T}}$ 
of the different of $\mc{S}/\mc{T}$ is bounded above by the
degree $\delta_{L_n'/k[[t]]}$ of the different of $L_n'/k[[t]]$.  Then (i), (ii), and (iii)
follow exactly as in the proof of \cite[Key Lemma 3.2]{Po:oc} (in fact, the argument is marginally easier, as our Witt vectors have no constant terms, so
there is no need for Pop's notion of ``quasi standard form").  And
(iv) follows immediately from (iii), since $(s)$ is the only branched
ideal of $\mc{T}/\mc{R}$.

Using Hilbert's different formula (\cite[p.\ 311]{ZS:ca} or \cite[IV, Proposition 4]{Se:lf}) and the definition of the upper numbering, we obtain 
$$\delta_{L_n'/k[[t]]} = \sum_{i=1}^n (u_i' + 1) (p^i - p^{i-1}).$$  For $\delta_{\mc{S}/\mc{T}}$, we add up the contributions from the
different branched ideals separately.  For the ideal $(t)$, we consider the extension of complete discrete valuation fields given by tensoring
$\mc{S}/\mc{T}$ with $k((\varpi))((t))$ over $\mc{T}$.  Let $(P_1, \ldots, P_n)$ be the standard form (relative to $(t)$)
of the Witt vector $(F_1, \ldots, F_n)$ classifying this extension.
Then the degree of the pole of $P_i$ at $t=0$ is bounded by $u_i$, and the upper jumps
are bounded by $(u_1, \ldots, u_n)$.  Thus the contribution
$\delta_{(t)}$ from the ideal $(t)$ to $\delta_{\mc{S}/\mc{T}}$ satisfies
$$\delta_{(t)} \leq \sum_{i=1}^n (u_i+1) (p^i - p^{i-1}).$$
For each ideal $(\zeta_m^{\alpha}t - \mu)$, we consider the extension of complete discrete valuation fields given by tensoring
$\mc{S}/\mc{T}$ with $k((\varpi))((\zeta_m^{\alpha}t - \mu))$.  Let $(P_{1, \alpha}, \ldots, P_{n, \alpha})$ be the standard form of
$(F_1, \ldots, F_n)$ relative to $(\zeta_m^{\alpha}t - \mu)$.  Then
the degree of the pole of $P_{i, \alpha}$ is bounded above by $(u_i' - u_i)/m$.  In fact
the inequality is strict, because $u_i' - u_i$ is divisible by $p$.
So the contribution $\delta_{\alpha}$ from the ideal $(t -
\zeta_m^{\alpha}\mu)$ to $\delta_{\mc{S}/\mc{T}}$ satisfies
$$\delta_{\alpha} \leq \sum_{i=1}^n (\frac{u_i' - u_i}{m} - 1 + 1)(p^i - p^{i-1}).$$  We conclude: 
\begin{eqnarray*}
\delta_{\mc{S}/\mc{T}} &=& \delta_{(t)} + \sum_{\alpha = 1}^m \delta_{\alpha} \\
&\leq& \sum_{i=1}^n (u_i+1) (p^i - p^{i-1}) + m \sum_{i=1}^n (\frac{u_i' - u_i}{m})(p^i - p^{i-1}) \\
&=& \sum_{i=1}^n (u_i' + 1)(p^i - p^{i-1}) \\
&=& \delta_{L_n'/k[[t]]}
\end{eqnarray*}
\Endproof

We omit the proof of the following proposition, which follows from Proposition \ref{Pcharpoort} exactly as \cite[Theorem 3.6]{Po:oc} follows from
\cite[Key Lemma 3.2]{Po:oc}.

\begin{prop}\label{Pcharpglobal}
Let $Y \to W$ be a branched $\Gamma$-cover of projective smooth $k$-curves.  Suppose that the local inertia at each ramification point with non-abelian
inertia group has vanishing KGB obstruction.  Set $\mc{W} = W \times_k k[[\varpi]]$.  Then there is a $\Gamma$-cover of 
projective smooth $k[[\varpi]]$-curves $\mc{Y} \to \mc{W}$ with special fiber the $\Gamma$-cover $Y \to W$ such that the ramification points
on the generic fiber $\mc{Y}_{\eta} \to \mc{W}_{\eta}$ with non-cyclic inertia have no essential ramification.
\end{prop}

{\bf Proof of Proposition \ref{Ppopreduction}:} 
Let $Y \to W = \proj^1$ be the Harbater-Katz-Gabber cover
associated to $L_n'/k[[s]]$ (this is called an \emph{HKG-cover} in \cite{Po:oc}).   This is a $\Gamma$-cover
that is \'{e}tale outside $s=0, \infty$, tamely ramified above $s = \infty$, and totally ramified above $s=0$ 
such that the formal completion of $Y \to W$ at $s=0$ yields the extension $L_n'/k[[s]]$.
Let $\mc{Y} \to \mc{W}$ be the $\Gamma$-cover guaranteed by Proposition \ref{Pcharpglobal}, and let $\mc{Y}_{\eta} \to \mc{W}_{\eta}$ be
its generic fiber. 
Recall that we assume that \emph{every} local $\Gamma$-extension $L_n/k[[s]]$ with no essential ramification lifts to characteristic zero.  
Furthermore, by the (standard) Oort conjecture, every cyclic extension
of $k[[s]]$ lifts to characteristic zero.
So if we base change $\mc{Y}_{\eta} \to \mc{W}_{\eta}$ to the algebraic closure of $k((\varpi))$, the local-global principle tells us that this cover 
lifts to characteristic zero.  Then, \cite[Proposition 4.3]{Po:oc} tells us that there is a \emph{rank two} characteristic zero valuation ring $\mc{O}$ with
residue field $k$ such that the $\Gamma$-cover $Y \to W$ has a lift over $\mc{O}$.  Note that this process works starting with any $\Gamma$-extension with
upper jumps $(u_1', \ldots, u_n')$, and that such extensions can be parameterized by some affine space $\aff^N$ 
(with one coordinate corresponding to each possible coefficient in an entry of a classifying Witt vector in standard form).  

To conclude, we remark that \cite[Proposition 4.7]{Po:oc} and its setup carry through exactly in our situation, with our $\aff^N$ playing the role of 
$\aff^{|\mathbf{\iota}|}$ in \cite{Po:oc}.  Indeed, we have that the analog of $\Sigma_{\mathbf{\iota}}$ in that proposition contains all closed points,
by the paragraph above.  Thus we can in fact
lift $Y \to W$ over a \emph{discrete} characteristic zero valuation ring.  Applying the easy direction of the local-global principle, we obtain a lift of $L_n'/k[[s]]$.  This 
concludes the proof of Proposition \ref{Ppopreduction}.

\section{The induction process}\label{Sind}
Let $L_n/k[[s]]$ be a $\Gamma = \ZZ/p^n \rtimes \ints/m$-extension, with $k[[t]]/k[[s]]$ the intermediate $\ints/m$-extension, and assume without
loss of generality that $t^m = s$.  
As in \cite{OW:ce}, the local-global principle thus shows that solvability of the local lifting problem from $L_n/k[[s]]$ is equivalent
to the following claim, which will be more convenient to work with:

\begin{claim} \label{claim1}
  Given a $\Gamma$-Galois extension $L_n/k[[s]]$, 
  then after possibly changing the uniformizer $s$ of $k[[s]]$,
  there exists a $\Gamma$-Galois cover $Y_n\to W:=\PP^1_K$ 
  (where $K$ is the fraction field of some characteristic zero DVR $R$
  with residue field $k$) with the following properties:
  \begin{enumerate}
  \item The cover $Y_n\to W$ has good reduction with respect to the standard
    model $\PP_R^1$ of $W$ and reduces to a $\Gamma$-Galois cover
    $\Yb_n\to \Wb=\PP^1_k$ (with $s$ as coordinate on $\Wb$) which is totally ramified above
    $s=0$, tamely ramified above $s = \infty$, and \'etale everywhere
    else.  In other words, $\Yb_n \to \Wb$ is the Harbater-Katz-Gabber
    cover for $L_n/k[[s]]$.
  \item
     The completion of $\bar{Y}_n\to\bar{W}$ at $s=0$ yields $L_n/k[[s]]$.
  \end{enumerate}
\end{claim}  

We write $\Yb_n \to \Xb$ (resp.\ $Y_n \to X$) for the unique $\ints/p^n$-subcover of $\Yb_n \to \Wb$ (resp.\ $Y_n \to W$).  
Then the quotient covers $\Xb \to \Wb$ and $X \to W$ are both tamely ramified $\ints/m$-covers of $\proj^1$'s, and we choose coordinates
$T$ on $X$ and $S$ on $W$ such that $T$ (resp.\ $S$) reduces to $t$ (resp.\ $s$) on $\Xb$ (resp.\ $\Wb$), and such that
$X \to W$ identifies $S$ with $T^m$. 


If $R$ is a characteristic zero DVR with residue field $k$ and fraction field $K$, set  
$D(r) = \{ T \in \Kb \mid \; v(T) > r\}$, where $v$ is the unique
valuation on $\Kb$ (with value group $\rats$) prolonging the valuation on $K$. We think of this disk as lying in $X$.

We prove Theorem \ref{Tmachine} (in the context of Claim \ref{claim1}) by induction using the following base case (Lemma \ref{Lbase})
and induction step (Theorem \ref{Tsetup}).

\begin{lem}\label{Lbase}
Let $L_1/k[[s]]$ be a $\ZZ/p \rtimes \ints/m$-extension whose $\ints/p$-subextension has upper ramification break $u_1$.
Suppose that $L_1/k[[s]]$ has vanishing KGB obstruction.   
Then there exists a $\ZZ/p \rtimes \ints/m$-cover $Y \to W$ satisfying Claim \rm{\ref{claim1}} for $L_1/k[[s]]$, such that
$Y \to X = \proj^1$ is \'{e}tale outside the open disk $D(r_1)$, where $r_1 = 1/u_1(p-1)$.
\end{lem}

\begin{thm}\label{Tsetup}
Suppose $n > 1$, and let $L_n/k[[s]]$ be a $\ZZ/p^n \rtimes \ints/m$-extension with vanishing KGB obstruction
whose $\ints/p^n$-subextension $L_n/k[[t]]$ has upper 
ramification breaks $(u_1, \ldots, u_n)$.  Let $L_{n-1}/k[[s]]$ be the unique $\ints/p^{n-1} \rtimes \ints/m$-subextension.  
Suppose there exists a $\ZZ/p^{n-1} \rtimes \ints/m$-cover $$Y_{n-1} \lpfeil{\ints/p^{n-1}} X \lpfeil{\ints/m} W$$  
satisfying Claim {\rm \ref{claim1}} for $L_{n-1}/k[[s]]$, such that $Y_{n-1} \to X$ is \'{e}tale outside the open disk $D(r_{n-1})$, 
where $r_{n-1} = 1/u_{n-1}(p-1)$.  Assume that $(p, m, u_{n-1}, N_1)$ satisfies the isolated differential data criterion, where
$N_1 = (p-1)u_{n-1}$ if $u_n = pu_{n-1}$, and $N_1 = (p-1)u_{n-1} - m$ otherwise.  Lastly, assume $u_n < pu_{n-1} + mp$.
Then there is a $\ZZ/p^n \rtimes \ints/m$-cover $Y_n \to W$ satisfying
Claim {\rm \ref{claim1}} for $L_n/k[[s]]$, such that $Y_n \to X$ is \'etale outside $D(r_n)$, where
$r_n = 1/u_n(p-1)$. 
\end{thm}

Theorem \ref{Tmachine} now follows immediately from Lemma \ref{Lbase}
and Theorem \ref{Tsetup} by induction.  After we
prove Lemma \ref{Lbase} in the next section, we devote most of the rest of the paper to proving Theorem \ref{Tsetup}.

\section{The base case} \label{Sbase}

In this section, we prove Lemma \ref{Lbase}.  Maintain the notation of
\S\ref{Sind}, and assume that we are in the situation of Lemma
\ref{Lbase}. Let $\Gamma = \Gal(L_1/k[[s]])$.  
By \cite[Theorem 2.1]{BWZ:dd}, the local lifting problem
holds for $L_1/k[[s]]$, so there is a $\Gamma$-cover $Y \to W$ satisfying Claim \ref{claim1} as desired 
(the vanishing of the KGB obstruction is exactly the 
condition in the theorem in \cite{BWZ:dd}).  So we need only check that the branch points of the $\ints/p$-subcover $Y \to X = \proj^1$ 
lie in $D(r_1) = D(1/u_1(p-1))$.  We start with a lemma.

\begin{lem}\label{Lautval}
In order to prove Lemma \ref{Lbase} for $L_1/k[[s]]$, it suffices to prove it for any $\Gamma$-extension $L'/k[[s]]$ with the same ramification break.
\end{lem}

\proof By \cite[Lemma 2.1.2]{Pr:cw3}, there is a $k$-automorphism
$\phi$ of $k[[s]]$ giving rise to an isomorphism from $L'$ to $L_1$
making the diagram below commute:
\[
\xymatrix{
L' \ar[r]^{\sim} & L_1  \\
k[[s]] \ar[r]^{\phi} \ar[u] & k[[s]] \ar[u]
}.
\]
Write $\phi(s) = a_1s + a_2 s^2 + \cdots$, where all $a_i \in k$ and $a_1 \in k^{\times}$.  Now, say 
$f: Y_1 \to X \to W$ satisfies Lemma \ref{Lbase} for $L'/k[[s]]$ with all branch points of $Y_1 \to X$ lying in $D(r_1)$.  Consider the cover
$f \times_W \Spec R[[S]]$.  Let $\Phi \in \Aut(R[[S]])$ be any $R$-automorphism lifting $\phi$. 
Identifying points of $\MaxSpec R[[S]]$ with (Galois orbits of) points of $\Kb$ of absolute value $< 1$, we have that
$\Phi^*$ preserves absolute values, because $|A_1S + A_2S^2 + \cdots| = |S|$ whenever all $A_i \in R$ with $A_1 \in R^{\times}$ and $|S| < 1$.  
Thus, the branch points of $\Phi^*(f \times_W \Spec R[[S]])$ have the same absolute values as those of
$f \times_W \Spec R[[S]]$, and $\Phi^*(f \times_W \Spec R[[S]])$ is a local lifting for $L_1/k[[s]]$.  
Clearly, if $\Phi$ is extended to $\Aut(R[[T]])$, where $T^m = S$, then $\Phi^*$ preserves absolute values as well. 
Applying the local-global principle gives Lemma \ref{Lbase}.
\Endproof

We are reduced to showing that, given $u_1 \equiv -1 \pmod{m}$ and $p \nmid u_1$, Lemma \ref{Lbase} holds for some
$L'/k[[s]]$ whose $\ints/p$-subextension $L'/k[[t]]$ has ramification break $u_1$.  We will 
freely use the terminology of \emph{Hurwitz trees} for the rest of this section (see \cite[\S3]{BW:ll}, especially Definition 3.2), as they are the key to
the proof of \cite[Theorem 2.1]{BWZ:dd}.

In particular, for any possible $u_1$ (called $h$ in \cite{BWZ:dd} and
\cite{BW:ll}), 
a Hurwitz tree is constructed in \cite{BW:ll} that gives rise to a lift of some $L'/k[[s]]$
whose $\ints/p$-subextension has ramification break $u_1$.  The valuations of the branch points of the lift (in terms of the coordinate $T$) 
can be read off from this Hurwitz tree.  This is done in the local context in \cite{BW:ll}, but the local-global principle allows us to conclude the global
result of Lemma \ref{Lbase}.  We split the proof up into the two cases
$u_1 < p$ and $u_1 > p$.  

If $u_1 < p$, then the Hurwitz tree is irreducible (\cite{GM:op}).
Thus the underlying combinatorial tree consists of two vertices: a
root vertex $v_0$ and a vertex $v_1$.  The points in the set $B$ of 
\cite[Definition 3.2]{BW:ll} all lie on $v_1$.  Since the points in $B$ represent the specializations of branch 
points of $Y \to X$, the valuation of each of these branch points is equal to $p$ times the
thickness $\epsilon$ 
of the edge connecting to $v_0$ and $v_1$ (the factor of $p$ comes from \cite[Proposition 2.3.2]{Ra:sp}).  
Since the conductor of the Hurwitz tree is $u_1$, we see that $|B| = u_1 + 1$.  Since the 
differential form $\omega_1$ on $v_1$ has simple poles at the points of $B$ and no other zeroes or poles aside from a zero at the point $z$ 
corresponding to the unique edge $e$, this zero has order $u_1 - 1$.  Then the definition of Hurwitz tree implies that
$$1 = (p-1)u_1\epsilon,$$ or $$\epsilon = 1/u_1(p-1).$$  Since $r_1 < p\epsilon = p/u_1(p-1)$, this case is proved.

If $u_1 > p$, then \cite[Theorem 4.3]{BW:ll} gives a construction of
the appropriate Hurwitz tree when $m=2$, splitting the construction into two cases. 
In both cases, the underlying combinatorial tree has a root vertex $v_0$, a vertex
$v_1$, and several other vertices.  Furthermore, in both cases, the different $\delta_{v_1}$ can be any rational number in $(0, 1)$.
Again, the valuation of each of the branch points is $p\epsilon$, where $\epsilon$ is the thickness of the edge connecting $v_0$ to $v_1$.  The definition
of Hurwitz tree implies that $$\delta_{v_1} = (p-1)u_1\epsilon,$$ and
taking $\delta_{v_1} > 1/p$ ensures that $p\epsilon > r_1$.  As is
mentioned in \cite[Proof of Theorem 2.1]{BWZ:dd}, this can be generalized easily to the case $m > 2$.  One
has the same freedom for $\delta_{v_1}$ and $\epsilon$ in this case. This completes the 
proof of Lemma \ref{Lbase}.  

\begin{rem}
The global context of \S\ref{Sind} was simply an encumbrance in this section, but it will be helpful later on.
\end{rem}
 
\section{Characters and Swan conductors} \label{Sgeom}

In this section we recall the tools of \emph{characters} and \emph{Swan conductors} from \cite[\S5]{OW:ce}.
Characters will serve as a substitute for Galois covers, as they are more convenient to manipulate algebraically.
We will also relate equivariance of characters to metacyclic extensions.

\subsection{Geometric setup}\label{Sgeomsetup}

Let $X = \proj^1_K$. 
We write $\KK = K(T)$ for the function field of $X$. 
Fix a smooth $R$-model $X_R$ of $X$, corresponding to the coordinate $T$.   We let $\Xb:=X_R\times_Rk$
denote the special fiber of $X_R$, and we let $X^{\rm an}$ denote the rigid analytic space associated to $X$.
We write $0$ for the $K$-point $T = 0$ and $\ol{0}$ for its specialization to $\Xb$.    

Let
\[
        D := ]\ol{0}[_{X_R} \subset X^{\rm an},
\]
be the open unit disk around $0$, that is, the set of points of $X^{\rm an}$ specializing to $\ol{0}\in\Xb$
(\cite{BL:sr2}). Then $\hat{\OO}_{X_R,\ol{0}}=R[[T]]$, and via $T$, we make an identification
\[
     D \cong \{\, x\in(\AA^1_K)^{\rm an} \mid v(x)>0 \,\}
\]

For $r\in\QQ_{\geq 0}$ we define
\[
     D[r] := \{\, x \in D \mid v(x)\geq r \,\}
\]
and, as in \S\ref{Sind}
\[
     D(r) := \{\, x \in D \mid v(x) > r \,\}
\]

We have $D(0)=D$. For $r>0$ the subset $D[r]\subset D$ is an affinoid
subdomain. Let $v_r:\KK^\times\to\QQ$ denote the ``Gauss
valuation" with respect to $D[r]$. This is a discrete valuation on $\KK$ which
extends the valuation $v$ on $K$ and has the property $v_r(T)=r$. It
corresponds to the supremum norm on the open subset $D[r]\subset X^{\rm an}$.

Let $\kappa_r$ denote the residue field of $\KK$ with respect to the valuation
$v_r$. For $r=0$, we have that $\kappa_0$ is naturally identified with
the function field of $\Xb$.  After replacing $K$ by a finite extension (which
depends on $r$) we may assume that $p^r\in K$. Then $D[r]$ is isomorphic to a
closed unit disk over $K$ with parameter $T_r:=p^{-r}T$. Moreover, the residue
field $\kappa_r$ is the function field of the canonical reduction $\Db[r]$ of
the affinoid $D[r]$. In fact, $\Db[r]$ is isomorphic to the affine line over
$k$ with function field $\kappa_r=k(t)$, where $t$ is the image of $T_r$ in
$\kappa_r$.  We make this identification of $t$ with the reduction of
$T_r$ throughout, whenever it is clear which $r$ we are dealing with.

For a closed point $\xb \in \Db[r]$, we let $\ord_{\xb}:
\kappa_r^\times\to\ZZ$ denote the normalized discrete valuation corresponding
to the specialization of $\xb$ on $\Db[r]$. We let $\ord_\infty$ denote the
unique normalized discrete valuation on $\kappa_r$ corresponding to the `point
at infinity'.

\begin{notation}\label{Dval}
For $F\in \KK^\times$ and $r\in\QQ_{\geq 0}$, we let $[F]_r$ denote the image
of $p^{-v_r(F)}F$ in the residue field $\kappa_r$. 
\end{notation}

\subsection{Characters}\label{Scharacters}

We fix $n\geq 1$ and assume that $K$ contains a primitive $p^n$th root of
unity $\zeta_{p^n}$ (this is true after a finite extension of $K$).  For an
arbitrary field $L$, we set
\[
      H^1_{p^n}(L):=H^1(L,\ZZ/p^n\ZZ).
\]
In the case of $\KK$, we have         
\[
      H^1_{p^n}(\KK):=H^1(\KK,\ZZ/p^n\ZZ)
         \cong \KK^\times/(\KK^\times)^{p^n}
\]
(the latter isomorphism depends on the choice of $\zeta_{p^n}$). Elements of
$H^1_{p^n}(\KK)$ are called {\em characters} on $X$.  Given an element $F\in
\KK^\times$, we let $\K_n(F)\in H^1_{p^n}(\KK)$ denote the character
corresponding to the class of $F$ in $\KK^\times/(\KK^\times)^{p^n}$.

For $i=1,\ldots,n$ the
homomorphism
\[
     \ZZ/p^i\ZZ \to \ZZ/p^n\ZZ, \qquad a \mapsto p^{n-i}a,
\]
induces an injective homomorphism $H^1_{p^i}(\KK)\inj H^1_{p^n}(\KK)$. Its
image consists of all characters killed by $p^i$. We consider $H^1_{p^i}(\KK)$
as a subgroup of $H^1_{p^n}(\KK)$ via this embedding.

A character $\chi\in H^1_{p^n}(\KK)$ gives rise to a branched
Galois cover $Y\to X$. If $\chi=\K_n(F)$ for some $F\in\KK^\times$, then $Y$
is a connected component of the smooth projective curve given generically by
the Kummer equation $y^{p^n}=F$. If $\chi$ has order $p^i$ as element of
$H^1_{p^n}(\KK)$, then the Galois group of $Y\to X$ is the unique subgroup of
$\ZZ/p^n\ZZ$ of order $p^i$.

A point $x\in X$ is called a {\em branch point} for the character $\chi\in
H^1_{p^n}(\KK)$ if it is a branch point for the cover $Y\to X$. The {\em
  branching index} of $x$ is the order of the inertia group for some point
$y\in Y$ above $x$. The set of all branch points is called the {\em branch
  locus} of $\chi$ and is denoted by $\BB(\chi)$.

\begin{defn}
  A character $\chi\in H^1_{p^n}(\KK)$ is called {\em admissible} if its branch
  locus $\BB(\chi)$ is contained in the open disk $D$.
\end{defn}

\subsubsection{Reduction of characters}

Let $\chi\in H^1_{p^n}(\KK)$ be an admissible character of order $p^n$, and let
$Y\to X$ be the corresponding cyclic Galois cover.  Let $Y_R$ be the normalization of
$X_R$ in $Y$. Then $Y_R$ is a normal $R$-model of $Y$ and we have
$X_R=Y_R/(\ints/p^n)$. 

After enlarging our ground field $K$, we may assume that the character $\chi$
is {\em weakly unramified} with respect to the valuation $v_0$, see
\cite{Ep:ew}. By definition, this means that for all extensions $w$ of $v_0$
to the function field of $Y$ the ramification index $e(w/v_0)$ is equal to
$1$. It then follows that the special fiber $\Yb:=Y_R\otimes_R k$ is reduced
(see e.g.\ \cite{AW:lp}, \S 2.2).

\begin{defn}
  We say that the character $\chi$ has \emph{\'etale reduction} if the map
  $\Yb\to\Xb$ is generically \'etale.  It has \emph{good reduction} if, in
  addition, $\Yb$ is smooth.
\end{defn}

In terms of Galois cohomology the definition can be rephrased as
follows.  Let $\KKh_0$ be the completion of $\KK$ at $v_0$.  The
character $\chi$ has \'etale reduction if and only if the image
$\chi|_{\KKh_0}$ of $\chi$ in $H^1_{p^n}(\KKh_0)$ under the restriction morphism induced by the inclusion $\Gal_{\KKh_0} \to \Gal_{\KK}$ is {\em
  unramified}.  The word ``unramified'' means that $\chi|_{\KKh_0}$ lies in the image of the
cospecialization morphism
\[
      H^1_{p^n}(\kappa_0)\to H^1_{p^n}(\KKh_0)
\]
(which is simply the restriction morphism induced by the projection
$\Gal_{\KKh_0}\to\Gal_{\kappa_0}$).  Since the cospecialization morphism is
injective, there exists a unique character $\chib\in H^1_{p^n}(\kappa_0)$
whose image in $H^1_{p^n}(\KKh_0)$ is $\chi|_{\KKh_0}$.  By construction, the
Galois cover of $\Xb$ corresponding to $\chib$ is isomorphic to an irreducible
component of the normalization of $\Yb$.

\begin{defn}
  If $\chi$ has \'{e}tale reduction, we call $\chib$ the {\em reduction} of $\chi$, and $\chi$ a {\em lift} of
  $\chib$. 
\end{defn}

\begin{rem}\label{Rchiram}
  Assume that $\chi$ is an admissible character with good reduction. The condition that $\chi$ is
  admissible implies that the cover $\Yb\to\Xb$ corresponding to the reduction
  $\chib$ is \'etale over $\Xb-\{\ol{0}\}$ (the proof uses Purity of Branch
  Locus, see e.g.\ \cite[Theorem 5.2.13]{Sz:gg}). Thus we may speak of the
  ramification breaks of $\chi$, by which we mean the ramification
  breaks above the point $\ol{0}$.
\end{rem}

\subsubsection{Equivariant characters}\label{Sequivariant}
In the context of \S\ref{Sgeomsetup}, consider a $\ints/m$-action on
$\KK$ fixing $K$, given by $\tau(T) = \zeta_m T$ for $\tau$ a generator of $\ints/m$.
This gives rise to a $\ints/m$-action on $X_R$, and we set $W_R$ (resp.\ $W$, $\Wb$) equal to $X_R/\langle \tau \rangle$ 
(resp.\ $X/\langle \tau \rangle$, $\Xb/\langle \tau \rangle$).  The action of $\ints/m$ on $\KK$ naturally gives rise to a $\ints/m$-action on 
$H^1_{p^n}(\KK) \cong \KK^{\times}/(\KK^{\times})^p$.  Let $\psi:
\ints/m \to \Aut(\ints/p^n)$ be a homomorphism.  Any automorphism of
$\ints/p^n$ is given by multiplication by an element of
$(\ints/p^n)^{\times}$, and we use this to identify $\Aut(\ints/p^n)$
with $(\ints/p^n)^{\times}$.

\begin{defn}\label{Dequiv}
A character $\chi \in H^1_{p^n}(\KK)$ is called \emph{$\psi$-equivariant} if $\tau(\chi) = \chi^{\psi(\tau^{-1})}$.
\end{defn}

\begin{rem}
Since $\chi$ is an element of a $p^n$-torsion group, the expression
$\chi^{\psi(\tau^{-1})}$ is well-defined.
\end{rem}

\begin{prop}\label{Pequiv}
Let $\Gamma = \ints/p^n \rtimes \ints/m$ via the conjugation action $\psi: \ints/m \to (\ints/p^n)^{\times} = \Aut(\ints/p^n)$.  
The $\ints/p^n$-branched cover $Y \to X$ given rise to by $\chi$ extends to a $\Gamma$-branched cover
$Y \to W$ if and only if $\chi$ is a $\psi$-equivariant character.
\end{prop}

\proof 
Letting $S = T^m$, proving the proposition is the same as showing that $K(Y)/K(S)$ is a $\Gamma$-extension if and only if $\chi$ is $\psi$-equivariant.  
Say $\Gamma$ is generated by $\tau$ and $\sigma$ of orders $m$ and $p^n$, respectively, with $\tau \sigma = \sigma^{\psi(\tau)} \tau$.  
Since $K(Y)/\KK$ is a Kummer extension, there exists a Kummer
generator $f \in K(Y)$ such that $\sigma(f) = \zeta_{p^n}f$.  Now,
$\chi$ being $\psi$-equivariant 
is equivalent to $\tau(f^{p^n}) = (f^{p^n})^{(\psi(\tau^{-1}))}g^{p^n}$, for some $g \in \KK^{\times}$ (here we abuse notation and think of 
$\psi(\tau^{-1})$ as any representative of $\psi(\tau^{-1})$ in $\ints$).  This is in turn equivalent to the possibility of extending
the action of $\tau$ from $\KK$ to $K(Y)$ via $\tau(f) = f^{\psi(\tau^{-1})}g$.  Suppose this is possible.  One calculates   
$$\tau (\sigma (f)) =  \zeta_{p^n} f^{\psi(\tau^{-1})}g = \sigma^{\psi(\tau)} (\tau (f)).$$  Since $|\Aut(K(Y)/K(S))| \leq mp^n$, we see that
the automorphism group is in fact generated by $\tau$ and $\sigma$
subject to $\tau \sigma = \sigma^{\psi(\tau)} \tau$.  
Thus, $K(Y)/K(S)$ is a $\Gamma$-extension.

On the other hand, if $K(Y)/K(S)$ is a $\Gamma$-extension, then using
$\tau \sigma^{\psi(\tau^{-1})} = \sigma \tau$, we have $\tau( \sigma^{\psi(\tau^{-1})} (f)) = \zeta_{p^n}^{\psi(\tau^{-1})}\tau(f) = \sigma(\tau(f))$, so
$$\frac{\sigma(\tau(f))}{\tau(f)} = \zeta_{p^n}^{\psi(\tau^{-1})}.$$ Kummer theory tells us that $\tau(f) = f^{\psi(\tau^{-1})}$ times an element of 
$\KK^{\times}$, which is exactly what we need to prove.
\Endproof

Note that a $\ZZ/p^n$-cover of $\PP^1_k$, unramified outside $\ol{0}$, is uniquely determined by its germ above the branch point 
(see, e.g., \cite{Ka:lg}).  Thus, with the above notation, and in light of \S\ref{Switt}, Claim \ref{claim1} may be reformulated as
follows.
 
\begin{claim}[cf.\ \cite{OW:ce}, Conjecture 5.7] \label{Creform} 
Let $\chib\in H^1_{p^n}(\kappa_0)$ (note $\kappa_0  \cong k(t)$) be a character of order $p^n$, unramified outside of
  $\ol{0}$, such that the corresponding Witt vector $(f_1, \ldots,
  f_n)$ is given by polynomials in $t^{-1}$ with all degrees of all terms congruent$\pmod{m}$.  
  Then (after replacing $K$ by a finite extension, if necessary)
  there exists a homomorphism 
  $\psi: \ints/m \to (\ints/p^n)^{\times}$ and an admissible, $\psi$-equivariant character $\chi\in H^1_{p^n}(\KK)$ with good
  reduction lifting $\chib$.
\end{claim}

By abuse of language, we will say that $\chib$ \emph{has vanishing
KGB obstruction} if the completion at $\ol{0}$ of the composite cover $\Yb \to
\Xb \to \Wb$ has vanishing KGB obstruction, where $\Yb \to \Xb$ is the
cover corresponding to $\chib$ and $\Xb \to \Wb$ is the quotient
morphism from the beginning of \S\ref{Sequivariant}.

\begin{rem}
In the case where $\chib$ has vanishing KGB obstruction, the corresponding homomorphism $\psi$ will be injective.
\end{rem}

The following lemma follows from an easy calculation, and will be useful in \S\ref{Splan}.
\begin{lem}\label{Lequiv}
Under the identification $H^1_{p^n}(\KK) \cong \KK^{\times}/(\KK^{\times})^{p^n}$, a character is $\psi$-equivariant if and only if
it can be identified with
$$\prod_{i = 1}^{m} (\tau^i(g))^{\psi(\tau^i)}$$ for some $g \in \KK^{\times}$.
\end{lem}

\begin{lem}\label{Lwhichroot}
If $\chib$ has vanishing KGB obstruction, $\tau(T) = \zeta_mT$, and $n
= 1$, 
then $\psi(\tau^{\ell}) = \zeta_m^{-\ell}$ as elements of $(\ints/p)^{\times} = \FF_p^{\times}$.
\end{lem}

\proof
This follows from \cite[Lemma 4.1(iii)]{OP:wt}, using the fact that $u_1 \equiv -1 \pmod{m}$.  In particular, our $\zeta_m$,
$u_1$, and $\psi_1(\tau)$ are the same as $\zeta^{-1}$, $j$, and $\alpha$, respectively, in \cite{OP:wt}.
\Endproof

It will at times be useful to measure how far an element of $\KK$ (in a special form) is from giving rise to a $\psi$-equivariant character of order $p$.
To this end, we make the following definition:

\begin{defn}\label{Ddiscrepancy}
Let $r \in \rats_{\geq 0}$.  Recall that $T_r = p^{-r}T$.
An element $F \in \KK \cap (1 + T_r^{-1}\m\{T_r^{-1}\})$ has \emph{$r$-discrepancy valuation $\geq \sigma$} if there exists 
$F' \in \KK \cap (1 + T_r^{-1}\m\{T_r^{-1}\})$ such that $\K_1(F')$ is 
$\psi$-equivariant and $v_r(F - F') \geq \sigma$.  If $\K_1(F)$ is
$\psi$-equivariant, we may say that the $r$-discrepancy valuation is $\infty$.
\end{defn}

\begin{defn}\label{Dpartialvaluation}
Let $r \in \rats_{\geq 0}$.  Suppose 
$F = \sum_{i= 0}^{\infty} \alpha_i T_r^{-i} \in R\{T_r^{-1}\}
\otimes_R K$.   Then we extend the valuation $v_r$ from $\KK$ to $R\{T_r^{-1}\}
\otimes_R K$ (and any subring) by setting $v_r(F) = \min_i(v(\alpha_i))$.
Furthermore, we write $v_r'(F) = \min_{i \in S}(v(\alpha_i))$, where
$S$ is the set of indices either divisible by $p$ or congruent to $-1 \pmod{m}$.
\end{defn}

\begin{lem}\label{Ldiscrepancy}
Fix $r \in \rats_{\geq 0}$.  Let $F \in \KK \cap (1 + T_r^{-1}\m\{T_r^{-1}\})$ with $\K_1(F)$ $\psi$-equivariant.
Write $[F-1]_r = \sum_{i=1}^{\infty} c_i t^{-i}$, where $t$ is the
reduction of $T_r$ in $\kappa_r$.
If $v_r(F - 1) < p/(p-1)$, then $c_i = 0$ unless $p|i$ or $i \equiv -1 \pmod{m}$.
If $v_r(F - 1) = p/(p-1)$ and $c_i = 0$ whenever $p|i$, then $c_i = 0$ unless $i \equiv -1 \pmod{m}$.
\end{lem}

\proof
Suppose $v_r(F -1) = \gamma \leq p/(p-1)$, and if equality holds, that $c_i = 0$ for $p|i$.
Write $$F = 1 + \sum_{i=1}^{\infty} a_i T_r^{-i} = 1 + A + B,$$ where $A$ consists exactly of the terms 
$a_i T_r^{-i}$ such that $v(a_i) = \gamma$, and $v_r(B) > \gamma$.  
Recall that $\tau$ is a generator of $\ints/m$ such that $\tau(T) = \zeta_mT$.  
By equivariance and Lemma \ref{Lwhichroot},
$$\tau(F) \equiv (1 + A + B)^{\zeta_m} \equiv 1 + \zeta_m A + B' \pmod{(\KK^{\times})^p},$$ where $v_r(B') > \gamma$.  
On the other hand, $\tau(F) = 1 + \tau(A) + \tau(B)$.

Assume, for a contradiction, that $A$ has
some term $a_i T_r^{-i}$ such that $i$ is neither congruent to $-1 \pmod{m}$ nor to $0 \pmod{p}$. 
In particular, $\tau(A) \neq \zeta_m A$.
We must show that $$Q := (1 + \zeta_m A + B')/(1 + \tau(A) + \tau(B))$$ is not a $p$th power in $\KK$ (in fact, we will show that it is not even a
$p$th power in $1 + T_r^{-1}\m[[T_r^{-1}]]$).    
The power series expansion of $Q$ in $T_r^{-1}$ is of the form $1 + \sum_{i=1}^{\infty} d_iT_r^{-i}$, with $v(d_i) \geq \gamma$ for all $i$. 
Since $\tau(A) \neq \zeta_m A$, there exists $i \in \nats$ such that $p \nmid i$ and $v(d_i) = \gamma$.  If 
$\gamma = v_r(Q-1) < p/(p-1)$, then $Q$ can only be a $p$th power if the $d_i$ such that $v(d_i) = \gamma$ all have $p|i$, giving a contradiction.  
If $\gamma = p/(p-1)$ and $A$ has no terms of degree divisible by $p$, then $Q$ can only be a $p$th power if there is \emph{some} $i$ with $p|i$
such that $v(d_i) = \gamma$, again a contradiction.
\Endproof

The discrepancy valuation of a power series sheds light on the valuation of its coefficients. 

\begin{cor}\label{Cdiscrepancy}
Let $r \in \rats_{\geq 0}$.  Suppose $F \in \KK \cap (1 + T_r^{-1}\m\{T_r^{-1}\})$ has discrepancy valuation $\geq \sigma$.  Then
$v_r(F-1) \geq \min(\sigma,\, p/(p-1),\, v_r'(F-1))$.
\end{cor}  

\proof
Pick $F' \in \KK \cap (1 + T_r^{-1}\m\{T_r^{-1}\})$ such that $\K_1(F')$ is $\psi$-equivariant and $v_r(F - F') \geq \sigma$.
It suffices to prove that $v_r(F'-1) \geq \min(p/(p-1), v_r'(F' - 1))$.  But this follows from Lemma \ref{Ldiscrepancy}.
\Endproof

\subsection{Swan conductors}\label{swan}
We recall some properties of the depth and differential Swan
conductors of characters.  For proofs, see \cite[\S5]{OW:ce}.
Let $\chi \in H^1_{p^n}(\KK)$ be a character.  As in \cite[\S5.3]{OW:ce}, we define the
\emph{depth Swan conductor} $\delta_{\chi}(r)$, which is a continuous, piecewise linear function
$$\delta_{\chi}: \reals_{\geq 0} \to \reals_{\geq 0}.$$
The kinks in $\delta_{\chi}(r)$ (i.e., non-differentiable points) occur only at rational values of $r$. 
As part of the definition, $\delta_\chi(r)=0$ if and
only if $\chi$ is unramified with respect to $v_r$. If this is the case then
the reduction $\chib_r\in H^1_{p^n}(\kappa_r)$ is well defined.

Let us now assume that $\delta_\chi(r)>0$, and that $r \in \rats_{\geq 0}$. 
Then, again as in \cite[\S5.3]{OW:ce}, one defines the {\em differential Swan conductor} of $\chi$ with respect to $v_r$,
\[
    \omega_\chi(r) \in\Omega_{\kappa_r}^1,
\] 
which we think of as a meromorphic differential on $\proj^1_k$
(perhaps more accurately, on $\Db[r]$).
The slopes of $\delta_{\chi}$ are determined by the orders of zeroes and poles of $\omega_{\chi}$: 
\begin{prop}[\cite{OW:ce}, Corollary 5.11]\label{Pdeltalin}
If $r \geq 0$ and $\delta_{\chi}(r) > 0$, then the left and right derivatives of $\delta_{\chi}$ at $r$ are given by
$\ord_{\infty}(\omega_{\chi}(r)) + 1$ and $-\ord_0(\omega_{\chi}(r)) - 1$, respectively.
\end{prop}

We now recall how the function $\delta_\chi$ determines whether $\chi$
has good reduction.
We fix an admissible character $\chi\in
H^1_{p^n}(\KK)$ of order $p^n$ and let $Y\to X$ denote the corresponding
Galois cover.

\begin{prop} \label{Pvccor1}
Let $\chi\in H^1_{p^n}(\KK)$ be an admissible character of order $p^n$. Then
the character $\chi$ has good reduction with upper ramification breaks $(u_1, \ldots, u_n)$
if and only if $\delta_\chi(0)=0$, the right slope of $\delta_{\chi}$ at $0$ is $u_n$, and
  \[	
    \abs{\{ x \in\BB(\chi) \, | \, \text{ramification index of } x \text{ is
      exactly } p^{n-i+1}\}} = u_i - u_{i-1},
  \]
  where we set $u_0 = -1$.
\end{prop}

\proof By definition, $\chi$ has \'{e}tale reduction if and only if $\delta_{\chi}(0) = 0$.  By \cite[Corollary 5.13(i) and Proposition 5.10(i)]{OW:ce}, 
$\chi$ has good reduction if and only if the right slope of $\delta_{\chi}$ at $0$ is equal to $|\BB(\chi)| - 1$, in which case \cite[Remark 5.8(i)]{OW:ce}
shows that this right slope is $u_n$ (note that \cite[Proposition 5.10(i)]{OW:ce}, is not stated as applying to $r = 0$, but from its proof referencing
\cite{We:sc}, it is 
clear that the right slope statement does apply).  Now the proposition follows from \cite[Corollary 5.13(ii)]{OW:ce}.  
\Endproof

\begin{prop}[\cite{OW:ce}, Corollary 5.15] \label{Pvccor2}
  Let $\chi\in H^1_{p^n}(\KK)$ be an admissible character of order $p^n$, let
  $r\in\QQ_{>0}$, and let $\xb$ be a point on the canonical reduction of $D[r]$. Suppose $\delta_{\chi}(r) > 0$.  Then
  \[
  \ord_{\xb}(\omega_{\chi}(r)) \geq -\abs{\BB(\chi)\cap U(r,\xb)},   
  \]
where $U(r, \xb)$ is the residue class of $\xb$ on the affinoid
$D[r]$.  Equality holds if $\chi$ has good reduction.
\end{prop}

The depth and differential Swan conductors behave in the following way under addition of characters:
\begin{prop}[\cite{OW:ce}, Proposition 5.9]\label{Paddchar}
  Let $\chi_1, \chi_2 \in H^1_{p^n}(\KK)$, and let $\chi_3 = \chi_1\chi_2$.
  For $i \in \{1, 2, 3\}$ and $r \in \QQ_{\geq 0}$, set $\delta_i =
  \delta_{\chi_i}(r)$. If $\delta_i>0$ then we set
  $\omega_i:=\omega_{\chi_i}(r)$. If $\delta_i=0$ then $\chib_i\in
  H^1_{p^n}(\kappa_r)$ denotes the reduction of $\chi_i$ with respect to
  $v_r$.  
  \begin{enumerate}
  \item If $\delta_1 \ne \delta_2$ then $\delta_3 = \max(\delta_1, \delta_2)$.
    If $\delta_1 > \delta_2$ then $\omega_3 = \omega_1$.
  \item Assume $\delta_1 = \delta_2>0$. Then
    \[
       \omega_1 + \omega_2 \ne 0 \quad\Rightarrow\quad
          \delta_1=\delta_2=\delta_3,\;\;\omega_3=\omega_1+\omega_2
    \]
    and
    \[
        \omega_1 + \omega_2 = 0 \quad\Rightarrow\quad
           \delta_3<\delta_1.
    \]
  \item Assume $\delta_1 = \delta_2=0$. Then $\delta_3=0$ and
    $\chib_3=\chib_1\chib_2$. 
 \end{enumerate}
\end{prop}

Lastly, we relate differential Swan conductors with equivariance.
\begin{lem}\label{Lequivswan}
Let $\tau$ and $\psi$ be as in \S\ref{Sequivariant}.  
If $\chi \in H^1_{p^n}(\KK)$ is $\psi$-equivariant and $r \in
\QQ_{\geq 0}$ such that $\delta_{\chi}(r) > 0$, then $\omega_{\tau(\chi)}(r) = \psi(\tau^{-1})\omega_{\chi}(r),$
where $\psi(\tau^{-1})$, by abuse of notation, is identified with its
image under the ``reduction mod $p$'' map $(\ints/p^n)^{\times} \to \FF_p^{\times} \subseteq k^{\times}$.
\end{lem}

\proof Since $\tau(\chi) = \chi^{\psi(\tau^{-1})}$, this follows from Proposition \ref{Paddchar}(ii).
\Endproof

\subsection{Characters of order $p$} \label{orderp}

We will now describe in the special case $n=1$ how to determine the function
$\delta_\chi$ explicitly in terms of a suitable element $F\in\KK^\times$
corresponding to the character $\chi\in
H^1_p(\KK)\cong\KK^\times/(\KK^\times)^p$. 

\begin{prop}[cf.\ \cite{OW:ce}, Proposition 5.17] \label{Pdelta} 
  Let $F\in \KK^\times\backslash(\KK^\times)^p$, let
  $\chi:=\K_1(F)\in H^1_p(\KK)$, and let $r\in\QQ_{\geq 0}$. Suppose that
  $v_r(F)=0$, and that $g := [F]_r \not\in\kappa_r^p$.  
Suppose, moreover, that $\chi$ is weakly unramified with respect to
$v_r$ (which is always the case if $K$ is chosen large enough). 
  \begin{enumerate}
  \item
    We have 
    \[
        \delta_\chi(r) = \frac{p}{p-1} - v_r(F),
    \]
  \item
   If $\delta_{\chi}(r) > 0$, then
    \[
       \omega_\chi(r) = \begin{cases}
          \;\;dg/g & \text{if $\delta_\chi(r)=p/(p-1)$,} \\
          \;\;dg   & \text{if $0<\delta_\chi(r)<p/(p-1)$.}
                        \end{cases}
    \]
    If, instead, $\delta_{\chi}(r) = 0$, then $\chib$ corresponds to the
    Artin-Schreier extension given by the equation $y^p - y = g$.
  \end{enumerate}
\end{prop}

\section{Proof of Theorem \ref{Tsetup}}\label{Sproof}
\subsection{Plan of the proof} \label{Splan}

We continue with the notation of \S\ref{Sgeom}.  
Recall that $D$ is the unit disk in $(\AA^1_K)^{\rm an}$ centered at
$0$, and $D(r)$ and $D[r]$ are, respectively, the open and closed disks of radius
$|p|^r$ centered at $0$.  
We are given a character $\chib_n \in
H^1_{p^n}(\kappa_0)$ of order exactly $p^n$, unramified outside $\ol{0}$, with upper ramification breaks $(u_1, u_2, \ldots, u_n)$, corresponding to a non-abelian
$\Gamma := \ints/p^n \rtimes \ints/m$-extension as in Claim
\ref{Creform}.  We assume that $\chib_n$ has vanishing KGB obstruction
(see after Claim \ref{Creform}).  We further assume that
$n \geq 2$.  For $1 \leq i \leq n$, set $r_i = 1/u_i(p-1)$.  Recall that $p \nmid u_1$, that $u_1 \equiv \cdots \equiv u_n \equiv -1 \pmod{m}$, and that
\[
      pu_{i-1} \leq u_i < pu_{i-1} + mp,
\]
for $i=1,\ldots,n$, where we set $u_0 = 0$. It is automatic that if the first inequality above is strict then
$p \nmid u_i$.  For $i=1,\ldots,n$ we set $\chib_i:=\chib_n^{p^{n-i}}\in H^1_{p^i}(\kappa_0)$. 
By assumption, for each $1 \leq i < n$, there is a compatible sequence of injective homomorphisms $\psi_i: \ints/m \to (\ints/p^i)^{\times}$ 
(i.e., $\psi_j$ reduces to $\psi_i$ for $j \geq i$), and a $\psi_i$-equivariant character $\chi_{i}$ lifting 
$\chib_{i}$.  We assume that $\BB(\chi_{n-1})$ lies in the disk $D(r_{n-1})$.  
Assume that $$(p, m, u_{n-1}, N_1)$$ satisfies the isolated differential data criterion (\S\ref{Sdd}), where $N_1 = (p-1)u_n$ if $u_n = pu_{n-1}$, or
$N_1 = (p-1)u_n - m$ if $u_n > pu_{n-1}$.  In order to prove Theorem
\ref{Tsetup} (using Claim \ref{Creform} in place of Claim \ref{claim1}), we must 
show that, for the unique $\psi_n: \ints/m \to (\ints/p^n)^{\times}$
compatible with the $\psi_i$, there exists an admissible
$\psi_n$-equivariant character $\chi_n \in H^1_{p^n}(\KK)$ with (good) reduction 
$\chib_n$.  Furthermore, we must show $\BB(\chi_n) \subseteq D(r_n)$. 
We will construct $\chi_n$ such that $\chi_n^p = \chi_{n-1}$.

We may assume that $\chi_{n-1}$ corresponds to an extension of $\KK$ 
given by a system of Kummer equations
\[
   y_i^p=y_{i-1}G_i, \quad i=1,\ldots,n-1
\]
with $y_0:=1$ and $G_i \in \KK$.  Any $\chi \in H^1_{p^n}(\KK)$ such that $\chi^p = \chi_{n-1}$ is given by an additional equation
\begin{equation}\label{gn}
  y_n^p=y_{n-1}G.
\end{equation}
Since we must have $\BB(\chi) \subseteq D$,
we will search for $G$ of the form 
\begin{equation}\label{Egnform0}
G=\prod_{j=1}^N(1-z_jT^{-1})^{a_j},
\end{equation}
where $a_j\in\ZZ$, $(a_j,p)=1$, and $z_j$ in the maximal ideal $\m$ of $R$ are pairwise
distinct (the ability to restrict our search to rational functions without worrying about missing anything 
is one benefit of working in the global context). 
We will say that the rational function $G$ \emph{gives rise} to the character $\chi$.  

\begin{rem}\label{Rpthpower}
Note that multiplication of $G$ by an element of $(\KK^{\times})^p$
does not change $\chi$, so when it is convenient, we will think of $G$
as an element of $\KK^{\times}/(\KK^{\times})^p$.
\end{rem}

\begin{lem}\label{Lhilb90}
It is possible to replace $G_{n-1}$ by its product with an element of $(\KK^{\times})^p$ 
so that $G = 1$ gives rise to a $\psi_n$-equivariant character $\chi$.  
\end{lem}

\proof 
Identifying $H^1_{p^n}(\KK)$ with $\KK^{\times}/(\KK^{\times})^{p^n}$, we have that choosing $G=1$ corresponds to
a character $\chi$ given by 
$A := G_1G_2^p\cdots G_{n-1}^{p^{n-2}}$.  Let $\tau$ be a generator of $\ints/m$. 
Since $\chi_{n-1}$ is assumed to be $\psi_{n-1}$-equivariant, Lemma \ref{Lequiv} tells us that
$$A \equiv \prod_{j = 1}^{m} (\tau^j(g))^{\psi_{n-1}(\tau^j)} \pmod{(\KK^{\times})^{p^{n-1}}}$$ for some $g \in \KK^{\times}$.
Since $\psi_n$ is compatible with $\psi_{n-1}$, we have that
$$A \equiv \prod_{j = 1}^{m} (\tau^j(g))^{\psi_n(\tau^j)} B^{p^{n-1}} \pmod{(\KK^{\times})^{p^n}},$$ for some $B \in \KK^{\times}$. 
Replacing $G_{n-1}$ by its product with $\prod_{j=1}^{m-1} (\tau^j(B))^{p\psi_n(\tau^j)}$ replaces $A$ by 
$$\prod_{j = 1}^{m} (\tau^j(gB^{p^{n-1}}))^{\psi_n(\tau^j)},$$ which is $\psi_n$-equivariant by Lemma \ref{Lequiv}.
\Endproof

Note that performing the replacement of Lemma \ref{Lhilb90} does not change the character $\chi_{n-1}$.
Thus, for the rest of the paper, we assume that $G_{n-1}$ is chosen in accordance with Lemma \ref{Lhilb90}.
 
\begin{prop}\label{Pwhenequiv}
Suppose a generator $\tau \in \ints/m$ sends $T$ to $\zeta_m T$.
If $G$ in the form of \eqref{Egnform0} gives rise to $\chi$, then in order for $\chi$ to be $\psi_n$-equivariant, 
it is necessary and sufficient that $G$ be of the form 
\begin{equation}\label{Egnform}
G = \prod_{j=1}^{N/m} \prod_{\ell=1}^m (1 - \zeta_m^{-\ell}z_jT^{-1})^{\psi_1(\tau^{\ell})a_j}
\end{equation}
after a possible reindexing.  Here we are viewing $\psi_1(\tau^{\ell})$ as an element of $\FF_p^{\times}$, which makes $G$ 
a well-defined element of $\KK^{\times}/(\KK^{\times})^p$.  
\end{prop}

\proof 
Identifying $H^1_{p^n}(\KK)$ with $\KK^{\times}/(\KK^{\times})^{p^n}$
via $\K_n$,
we have that $\chi = \K_n(G_1G_2^p\cdots
G_{n-1}^{p^{n-2}}G^{p^{n-1}})$.  
Since $\K_n(G_1G_2^p \cdots G_{n-1}^{p^{n-2}})$ 
is assumed to be $\psi_n$-equivariant, we have that $\chi$ being $\psi_n$-equivariant is equivalent to 
$\K_n(G^{p^{n-1}})$ being $\psi_n$-equivariant, which is equivalent to $\K_1(G)$ being
$\psi_1$-equivariant.  By Lemma \ref{Lequiv}, this is equivalent to $G$ having the desired form. 
\Endproof

\begin{rem}
We say that $G \in \KK^{\times}$ is ``of the form \eqref{Egnform}'' if its
residue class as an element of $\KK^{\times}/(\KK^{\times})^p$ is.
\end{rem}

Let us assume that none of the $z_i$ is a
branch point of $\chi_{n-1}$.  If this is the case, then Proposition \ref{Pvccor1} shows that a necessary condition
for good reduction of $\chi$ is that $N=u_n-u_{n-1}$.  We assume this.
Note that $N = |\BB(\chi) \backslash \BB(\chi_{n-1})|$.

We will try to find a choice $G_n$ for $G$ of the form
\eqref{Egnform} giving rise to a character
$\chi_n$ whose (good) reduction is $\chib_n$.  In \S\ref{Scrit}, we
give some contraints that $G_n$ will have to satisfy.  In
\S\ref{Sstrategy}, we give our strategy in more detail.

\subsection{The critical radius} \label{Scrit}

We continue with the setup of \S\ref{Splan}.  In particular, recall that $\chi_i$ is a lift of
$\chib_i$ for $1 \leq i < n$, and $\chi$ is the character arising from $G$.  
The number $r_{n-1} = 1/u_{n-1}(p-1)$ will be of the utmost importance, and we will refer to it as the \emph{critical radius}, or $r_{\crit}$.
From \cite[Eq. (15)]{OW:ce}, we know that
\begin{equation}\label{Edeltacrit}
\delta_{\chi}(r_{\crit}) = \frac{p}{p-1},
\end{equation} 
regardless of our choice of $G$ (this is, essentially, why the critical radius is ``critical." It is the minimal $r$ such that
$G$ does not affect $\delta_{\chi}(r)$).  For this section, we let $u$ be the minimal 
upper ramification break $u_i$ such that $u_{n-1}$ is a power of $p$
times $u$.  Thus $u$ is prime to $p$ (see \S\ref{Switt}, just before
Proposition \ref{Ppries}).  Set $\nu = n - 1 - i$, so that
$u_{n-1} = up^{\nu}$.  
%
%
%

We start with a proposition that is not necessary for the proof of Theorem
\ref{Tsetup}, but it guides some of our choices about how we
construct $G$.  
Namely, we already know that $G$ is of the form \eqref{Egnform}, and that $N = u_n - u_{n-1}$.  Our first decision is how many of
the branch points of $\chi$ to place at the critical radius.  That is,
how many of the $z_j$ in $\eqref{Egnform}$ should have valuation $r_{\crit}$?  

\begin{prop}[\cite{OW:ce}, Proposition 6.4] \label{Pdistrib}
  If $\chi$ has good reduction then the following hold.
  \begin{enumerate}
  \item
    For all $j$ we have $v(z_j)\leq r_{\crit}$.
  \item
    For $i,j$ with $v(z_i)=v(z_j)= r_{\crit}$ we have $\bar{x}_i\neq\bar{x}_j$
    (where $\bar{x}_j$ denotes the reduction of $x_j := z_jp^{-r_{\crit}}$).
  \item
    Write $N = N_1 + N_2$, where $N_1/m$ is the number of $z_j$ in \eqref{Egnform} with $v(z_j) = r_{\crit}$.  If $u_n=pu_{n-1}$ then $N_1=u_{n-1}
    (p-1)$ and $N_2=0$. Otherwise, $N_1<u_{n-1}(p-1)$ and $N_2>0$.
  \end{enumerate}
\end{prop}

Let $t = [T_{r_{\crit}}]_{r_{\crit}}$.  Since $\sum_{\ell = 1}^m \psi_1(\tau^{\ell}) = 0$ in $\FF_p$, it follows from \eqref{Egnform} that, up to reordering the $z_j$ and 
up to a constant factor that we may eliminate by rescaling $t$, we have
\begin{equation}\label{Egred}
     [G]_{r_{\crit}} = g=t^{a_0}\prod_{j=1}^{N_1/m} \prod_{\ell=1}^m (1 - \zeta_m^{-\ell}\bar{x}_jt^{-1})^{\psi_1(\tau^{\ell})a_j},
\end{equation}
where $p | a_0$.

\begin{cor}\label{Cgform}
  In the notation of Proposition {\rm \ref{Pdistrib}}, if $\chi$ has
  good reduction and $g = [G]_{r_{\crit}}$, then
  \begin{equation} \label{Ecrit}
      \omega_{\chi}(r_{\crit}) = \frac{dg}{g}-u\sum_{s=0}^{\nu}t^{-up^s-1}dt = 
          \frac{c\,dt}{t^{u_{n-1}+1}\prod_{j=1}^{N_1/m}(t^m-\bar{x}_j^m)},
  \end{equation}
  where $c$ is a nonzero constant.
  In particular, $\ord_{\infty} (\omega_{\chi}({r_{\crit}})) = N_1 + u_{n-1} - 1.$
\end{cor}

\proof The first equality follows from \cite[Lemma 6.3]{OW:ce}.  This
middle expression shows that
$\omega_{\chi}(r_{\crit})$ has no multiple poles outside of $t = 0$, where there is a pole of order $up^{\nu} + 1 = u_{n-1} + 1$.  
Furthermore, Proposition \ref{Pvccor2} shows that
$\omega_{\chi}(r_{\crit})$ has a simple pole at each of the $N_1$ points $\zeta_m^{\ell}\bar{x}_j$, 
no zeroes outside of $t = \infty$, and no other poles outside of 
$t = 0$.  It follows that $\omega_{\chi}(r_{\crit})$ has the form in the third expression, from which
$\ord_{\infty}(\omega_{\chi}(r_{\crit}))$ can be read off.
\Endproof

\begin{rem}
Notice that $t^{a_0}$ from \eqref{Egred} disappears in the logarithmic
derivative.
\end{rem}

Recall (\S\ref{Sdd}) that $(p, m, u_{n-1}, N_1)$ \emph{satisfies the differential data criterion} with respect to $k$ 
if there exists a polynomial $f(t) \in k[t^m]$ of degree exactly $N_1$
in $t$, such that the meromorphic differential form
$$\omega := \frac{dt}{f(t)t^{u_{n-1}+1}} \in \Omega^1_{k(t)/k}$$ satisfies
$\C(\omega) = \omega + u t^{-u_{n-1}-1}dt,$ where $u$ is the
prime-to-$p$ part of $u_{n-1}$.  Note that this implies $f(0) \neq 0$,
otherwise the order of the pole of $\omega$ at $t = 0$ will be too large compared
to that of $\C(\omega)$ and $t^{-u_{n-1}-1}dt$.  We will suppress $k$ when it is understood.  

\begin{prop}\label{Pcrit}
Suppose $p$, $m$, $u_{n-1}$ are as in this section, and $N_1$ is as in Proposition \ref{Pdistrib}.  The following are equivalent:
\begin{enumerate}
\item There exists $G$ of the form \eqref{Egnform} such that $g := [G]_{r_{\crit}}$ satisfies 
(\ref{Ecrit}). 
\item The quadruple $(p, m, u_{n-1}, N_1)$ satisfies the differential data criterion.
\item There exists a solution to the following system of equations:
\begin{equation}\label{Eass71}
\sum_{j=1}^{N_1/m} a_j\xb_j^{q} = 
\begin{cases} 
\;\;u/m, & q=u \\ \;\;0, & \text{otherwise,}
\end{cases} 
\end{equation}
where the $\xb_j \in k$, the $a_j \in \FF_p^{\times}$, and $q$ ranges over those numbers from $1$ to $N_1 + u_{n-1} - 1$ that are
congruent to $-1 \pmod{m}$ and not divisible by $p$.
\end{enumerate}
\end{prop}

\proof We first prove (i) implies (ii).  
Suppose $g$ is a solution to (\ref{Ecrit}). Taking $\omega = \omega_{\chi}(r_{\crit})$ and $f(t) = c^{-1}\prod_{j=1}^{N_1/m}(t^m-\bar{x}_j^m)$, 
and noting that the Cartier operator preserves logarithmic differential forms, it is clear that $\C(\omega) = \omega + ut^{-u_{n-1}-1}dt$.

Now we prove (ii) implies (i).  Suppose $(p, m, u_{n-1}, N_1)$ 
satisfies the differential data criterion via a differential form $\omega = dt/f(t)t^{u_{n-1} +1}$.
By the properties of the Cartier operator, $\omega$ is equal to a logarithmic differential form minus $u\sum_{s=0}^{\nu} t^{-up^s-1}dt$.
Since $t^{-up^s-1}$ has trivial residues, the residues of $\omega$ are the same as those of a logarithmic differential form.  In particular,
they lie in $\FF_p$ (\cite[Lemma 1.5]{BWZ:dd}).

Factor $f$ as $c^{-1}\prod_{j=1}^{N_1/m} (t^m - \bar{x}_j^m)$. Let $a_j$ be the residue of $\omega$ at $\xb_j$.  
An easy calculation shows that the residue of $\omega$ at $\zeta_m^{-\ell} \xb_j$ is $\zeta_m^{-\ell} a_j$, which is $\psi_1(\tau^{\ell})a_j$, by
Lemma \ref{Lwhichroot}.  Since $m | (p-1)$ (see Remark \ref{RKGBcons}), all of these 
residues lie in $\FF_p$.  Now, take 
\begin{equation}\label{Edefg}
g =  \prod_{j=1}^{N_1/m} \prod_{\ell=1}^m (1 -
\zeta_m^{-\ell}\bar{x}_jt^{-1})^{\psi_1(\tau^{\ell})a_j},
\end{equation}
where by abuse of notation we take an arbitrary lift of each
$\psi_1(\tau^{\ell})a_j$ and consider it as an element of $\ints$. 
Then $dg/g$ has the same residues at the simple poles $\zeta_m^{\ell}\xb_j$ as $\omega$.  So
$$\beta := dg/g - u\sum_{s=0}^{\nu} t^{-up^s-1}dt - \omega$$ is a logarithmic differential form with no poles outside of $0$.  Since a 
nonzero logarithmic differential form has only simple poles, and at least two of them, we conclude that $\beta = 0$.  So $g$ is a solution
to (\ref{Ecrit}).  Let $x_1, \ldots, x_{N_1/m}$ be lifts of $\xb_1,
\ldots, \xb_{N_1/m}$ to $R$.  Then we take $G$ to be anything in the
form \eqref{Egnform} such that $z_i = p^{r_{\crit}}x_i$ for $1 \leq i
\leq N_1/m$ and $v(z_i) < r_{\crit}$ for $i > N_1/m$.

Lastly, we prove that (i) and (iii) are equivalent (cf.\ \cite[p.\ 266]{OW:ce}).  We identify the choices of the $a_j$ and $\xb_j$ in (i) and (iii).
If $G$ is of the form \eqref{Egnform} with $g = [G]_{r_{\crit}}$, then differentiating logarithmically, we obtain
$$\frac{dg}{g} = \sum_{j=1}^{N_1/m} \sum_{\ell = 1}^m \frac{\psi_1(\tau^{\ell})a_j\zeta_m^{-\ell}\xb_jt^{-2}dt}{1-\zeta_m^{-\ell}\xb_jt^{-1}}.$$
Since $\psi_1(\tau^{\ell}) = \zeta_m^{-\ell}$ by Lemma \ref{Lwhichroot}, we obtain
\begin{equation}\label{Elaurent}
    \frac{dg}{g} = \sum_{q=1}^\infty \left(\sum_{j=1}^{N_1/m} \sum_{\ell =1}^m \zeta_m^{(-q-1)\ell}a_j\bar{x}_j^q\right) t^{-q-1}dt.
\end{equation}
Thus all terms in the expansion \eqref{Elaurent} disappear unless $q \equiv -1 \pmod{m}$.  
In particular, $\omega_{\chi}(r_{\crit}) = dg/g - u\sum_{s=0}^{\nu} t^{-up^s-1}dt$ has a zero of order at least $N_1 +u_{n-1} - 1$ at $\infty$ if and only if,
for all $q \equiv -1 \pmod{m}$ between $1$ and $N_1 + u_{n-1} - 1$ inclusive, we have
\begin{equation}\label{Ecoefftest}
\sum_{j=1}^{N_1/m} a_j\xb_j^{q} = 
\begin{cases} 
\;\;u/m, & q=u, up, \ldots, up^{\nu} \\ \;\;0, & \text{otherwise.}
\end{cases} 
\end{equation}
Now, if an equation in \eqref{Ecoefftest} holds for $q$, then it also holds for $pq$, as replacing $q$ with $pq$ simply raises both the left hand and 
right hand sides of the equation to the $p$th power.  So $\omega_{\chi}(r_{\crit})$ has a zero of order at least $N_1 + u_{n-1} - 1$ at $\infty$ if and only
if (iii) holds.  But $\omega_{\chi}(r_{\crit})$ cannot have a zero of order greater than $N_1 + u_{n-1} - 1$ at $\infty$, as it has at worst a pole of order
$u_{n-1} + 1$ at $0$ and $N_1$ simple poles at the $\zeta_m^{\ell}\xb_j$.  So (iii) is equivalent to $\omega_{\chi}(r_{\crit})$ having a zero of order
exactly $N_1 + u_{n-1} - 1$ at $\infty$, poles in the aforementioned places, and no other zeroes.  That is, (iii) is equivalent to (i).    
\Endproof 

\begin{rem}\label{Rnotjustexistence}
From its proof, it is clear that Proposition \ref{Pcrit} is not just an existence result.  In particular, any $f$ realizing the differential data criterion
gives rise to a $g$ satisfying \eqref{Ecrit}, which in turn gives rise
to a solution to the system \eqref{Eass71} (the $\xb_j$ in \eqref{Eass71} are
representatives from the $\mu_m$-equivalence classes of roots of $f(t)$).
The definition of realizing the differential data criterion
(Proposition \ref{Pcrit}(ii)) is easier
to state than the criterion in Proposition \ref{Pcrit}(iii), 
and is also usually easier to work with computationally, 
but it is the criterion of Proposition \ref{Pcrit}(iii) that we will
mostly use in our proofs.
\end{rem}

\begin{rem}\label{Rsquare}
One checks that \eqref{Eass71} is a system of $N_1/m$ equations in $N_1/m$ variables if and only if if $u_{n-1}(p-1) - mp < N_1 \leq u_{n-1}(p-1)$.
\end{rem}

\begin{rem}
The choice of the $a_j \in \FF_p^{\times}$ in \eqref{Egnform} is known as the ``type''
(cf.\ \cite{BW:ll}, \cite{BWZ:dd}, \cite{OW:ce}).  One of the
advantages of phrasing the differential data criterion
in terms of the Cartier operator, rather than in terms of the equations
\eqref{Eass71}, is that this phrasing is ``type independent.''  That is, one does not have
to determine the $a_j$ separately---they fall out automatically as the
residues of $\omega$, which is determined solely in terms of the roots
of $f$ (which correspond to the $\xb_j$).  In the papers mentioned above, one of the
difficulties is guessing the correct type in an analogous situation. 

Furthermore, since the problem is naturally symmetric in the $\xb_j$,
it makes sense to ``symmetrize'' things by thinking in
terms of $f$ instead.  The coefficients of $f$ will in general lie in smaller
fields than the $\xb_j$.
\end{rem}


In \S\ref{Scrit2hub}, it will become important not only to be able to satisfy the equivalent 
criteria of Proposition \ref{Pcrit}, but to do so in an ``isolated" fashion, that is,
to choose the $a_j$ and $\xb_j$ as in Proposition \ref{Pcrit}(iii) such that no infinitesimal deformation of the $\xb_j$ yields a solution to \eqref{Eass71}.
For fixed $a_j$, the Jacobian matrix of \eqref{Eass71} at a solution $(\xb_j)_j$ is the $N_1/m \times N_1/m$ matrix
\begin{equation}\label{Einvertiblematrix}
      \left(qa_j\xb_j^{q-1}\right)_{q,j} 
\end{equation}
over $k$, where $j$ ranges from $1$ to $N_1/m$ and $q$ ranges over those numbers from $1$ to $N_1 + u_{n-1} - 1$ that are
congruent to $-1 \pmod{m}$ and not divisible by $p$.  Thus,
in light of Remark \ref{Rnotjustexistence}, we make a definition (cf.\ \cite[Assumption 7.2]{OW:ce}).

\begin{defn}\label{Disolated}
Suppose $p$, $m$, $u_{n-1}$ are as in this section, and $N_1$ is as in Proposition \ref{Pdistrib}.  The quadruple $(p, m, u_{n-1}, N_1)$ satisfies
the \emph{isolated differential data criterion} if there is a
polynomial $f \in k[t^m]$ realizing the differential data criterion
for $(p, m, u_{n-1}, N_1)$ (equivalently, a $g \in k(t)$ satisfying \eqref{Ecrit})
that gives rise to a solution to the system of equations \eqref{Eass71} for which the matrix \eqref{Einvertiblematrix} is invertible
over $k$ (or is empty).
\end{defn}

\begin{rem}\label{Rconstants}
Dividing by (nonzero) constants, one sees that the isolated differential data criterion
holds if the matrix
$\left(\xb_j^{q-1}\right)_{q, j}$ (for the same $q$ and $j$ as in
\eqref{Einvertiblematrix}) is invertible.
\end{rem}

\begin{rem}
The differential data criterion is analogous to \cite[Assumption
7.1]{OW:ce} in the cyclic case, and the isolated differential data criterion is analogous
to \cite[Assumption 7.2]{OW:ce}.
\end{rem} 

\begin{defn} \label{Dgcrit} 
  If $g \in k(t)$ is a solution to \eqref{Ecrit} realizing the isolated
  differential data criterion, then we define $\G_{\crit, g} \subseteq
  \KK$ to be the set of all $G$ of the
  form \eqref{Egnform} (but with $N_1$ replacing $N$) with $[G]_{r_{\crit}}=g.$  
\end{defn}

To sum up, we have shown that a $G$ of the form \eqref{Egnform} can
only give rise to a character with good reduction if it lies in
$\G_{\crit, g}$, for some $g$ solving \eqref{Ecrit}.

\subsection{Plan of the proof, part II}\label{Sstrategy}
Maintain the notation of \S\ref{Splan} and \S\ref{Scrit}.  Recall that we are searching for $G_n$ of the form \eqref{Egnform} giving rise to a character 
$\chi_n$ with good reduction $\chib_n$.  The proposition below follows immediately from Proposition \ref{Pvccor1} and the
discussion at the beginning of \S\ref{Splan}.

\begin{prop}\label{Psufficient}
If $G$ is of the form \eqref{Egnform} such that all $z_j$ satisfy $v(z_j) > r_n = 1/u_n(p-1)$, such that $\delta_{\chi}(0) = 0$, and such that the right slope
of $\delta_{\chi}$ at $0$ is $u_n$, then $G$ gives rise to a
$\psi_n$-equivariant character $\chi$ with good reduction $\chib$ having 
ramification breaks $(u_1, \ldots, u_n)$, and such that $\BB(\chi) \subseteq D(r_n)$. 
\end{prop}

The argument outlined in the remainder of this section is the most
important difference between this paper and \cite{OW:ce}.

Recall that $pu_{n-1} \leq u_n < pu_{n-1} + mp$ (in fact, since all $u_i$ are congruent to $-1 \pmod{m}$, we
have $u_n \leq pu_{n-1} + m(p-1)$).  As was mentioned before Proposition \ref{Pdistrib}, we must decide how many of the $z_j$ to choose such that 
$v(z_j) = r_{\crit} = r_{n-1} = 1/(p-1)u_{n-1}$.  Recall that there
are $N/m$ $z_j$ in total.
Let $N_1$ and
$N_2$ be two multiples of $m$ such that $N_1 + N_2 = N = u_n -
u_{n-1}$.  If $u_n = pu_{n-1}$, we choose $N_1 = u_n - u_{n-1}$ and
$N_2 = 0$.
Otherwise, we take some $N_1$ such that $(p-1)u_{n-1} - mp < N_1 < (p-1) u_{n-1}$.  This gives $0 \leq N_2 \leq 2mp -2m$, with the first equality 
holding if and only if $u_n = pu_{n-1}$.  Note that $N_1 + u_{n-1}$ is divisible by $p$ if and only if $u_n = pu_{n-1}$. 
We will construct $G$ and a rational number $r_{\hub}$ such that $N_1/m$ of the $z_j$ satisfy $v(z_j) = r_{\crit}$ and the other $N_2/m$ of the
$z_j$ satisfy $v(z_j) = r_{\hub}$.  If $u_n = pu_{n-1}$, we declare
$r_{\hub} = 0$.  Otherwise, $0 < r_{\hub} < r_{\crit}$ is defined by the following proposition.  

\begin{prop}\label{Phublocation}
In the notation above, suppose $v(z_j) = r_{\crit}$ for $N_1/m$ of the $z_j$ and $v(z_j) = r_{\hub}$ for $N_2/m$ of the $z_j$.  Suppose further
that $N_2 > 0$ and $\chi$ has good reduction.  Then $$r_{\hub} = \frac{1}{N_2} - \frac{N_1}{(p-1)u_{n-1}N_2}.$$  Furthermore, $r_n < r_{\hub} < 
r_{\crit}$.
\end{prop}

\proof
Under the assumptions in the proposition, $\BB(\chi)$ has $N_2$ points
with valuation $r_{\hub}$ and another $N_1$ points with valuation
$r_{\crit}$.  Using Proposition \ref{Pvccor1} along with the fact
that $\chi_{n-1}$ has good reduction, $\BB(\chi)$ has exactly $u_{n-1} +
1$ other points, all of which have valuation greater than $r_{\crit}$.  
Since $\chi$ has good reduction, Propositions \ref{Pdeltalin} and \ref
{Pvccor2} imply that the right slope of
$\delta_{\chi}$ is $N_1 + N_2 + u_{n-1} = u_n$ for $0 \leq r < r_{\hub}$ and $N_1 + u_{n-1}$ for $r_{\hub} \leq r < r_{\crit}$.  
Furthermore, $\delta_{\chi}(r_{\crit}) = p/(p-1)$ by \eqref{Edeltacrit}, and $\delta_{\chi}(0) = 0$. Thus we obtain the equation
$$(N_1 + N_2 + u_{n-1}) r_{\hub} + (N_1 + u_{n-1})(r_{\crit} - r_{\hub}) = \frac{p}{p-1}.$$  This yields $r_{\hub} = 1/N_2 - N_1/(p-1)u_{n-1}
N_2$, proving the first part of the proposition.   

Since $(p-1)u_{n-1} - N_1 < u_n - u_{n-1} - N_1 = N_2$, it follows easily that $r_{\hub} < r_{\crit}$.  On the other hand, we know
$$N_2 = u_n - u_{n-1} - N_1 \leq (p-1)u_{n-1} - N_1 + m(p-1).$$ So 
$$r_{\hub} \geq \frac{(p-1)u_{n-1} - N_1}{(p-1)u_{n-1}((p-1)u_{n-1} - N_1 + m(p-1))}.$$  Now, since $(p-1)u_{n-1} - N_1 \geq m$, we have
$$\frac{(p-1)u_{n-1} - N_1}{(p-1)u_{n-1} - N_1 + m(p-1)} \geq \frac{1}{p}.$$  Putting everything together, we obtain
$$r_{\hub} \geq \frac{1}{(p-1)u_{n-1}p} > \frac{1}{(p-1)u_n} = r_n.$$
\Endproof

From the proof of Proposition \ref{Phublocation}, it is clear that if $\chi$ has good reduction, then 
\begin{equation}\label{Edhub}
\delta_{\chi}(r_{\hub}) = u_n r_{\hub} = \frac{p}{p-1} - (N_1 +
u_{n-1})(r_{\crit} - r_{\hub}) =: \delta_{\hub},
\end{equation}
regardless of whether $u_n = pu_{n-1}$.

\emph{We will work under the running assumption that $(p, m, u_{n-1}, N_1)$ satisfies the isolated differential data 
criterion.}  Thus, we let $g$ be a solution to \eqref{Ecrit} realizing the isolated differential data criterion, and we define $\G_{\crit, g}$ as in 
Definition \ref{Dgcrit}.  Our first step, to which \S\ref{Scrit2hub}
is devoted (and which parallels \cite{OW:ce} very closely), is to find $G_{\crit} \in \G_{\crit, g}$ such that 
$G_{\crit}$ gives rise to a character $\chi_{\crit}$ with
$\delta_{\chi_{\crit}}(r_{\hub}) = \delta_{\hub}$.

If $u_n = pu_{n-1}$, then $N_1 = N$, so $G_{\crit}$ is already of the form \eqref{Egnform}. In this case, we set $G = G_{\crit}$, from which
$\chi = \chi_{\crit}$ satisfies $\delta_{\chi}(0) = 0$, and the right slope of $\delta_{\chi}$ at $0$ is $u_n$.  
Since $G$ is already in the
form \eqref{Egnform}, Proposition \ref{Psufficient} shows that $\chi$ has good reduction
$\chib$ with upper ramification breaks $(u_1, \ldots, u_n)$. 
We then show quite easily that $G$ can be replaced by some $G_n \in \G_{\crit,g}$ (and thus still of the form \eqref{Egnform}) 
such that $G_n$ gives rise to a character $\chi_n$ with good reduction $\chib_n$ and $\BB(\chi_n) \subseteq D(r_n)$.  By Proposition
\ref{Pwhenequiv} and the discussion at the beginning of \S\ref{Splan},
this proves Theorem
\ref{Tsetup} (using Claim \ref{Creform} instead of Claim \ref{claim1}) when $u_n = pu_{n-1}$.

If $u_n > pu_{n-1}$, our next step (\S\ref{Shub}) is to construct a space $\G_{\hub} \subseteq \KK$ consisting of certain functions 
whose images in $\KK^{\times}/(\KK^{\times})^p$ have the form
$$\prod_{j=N_1/m+1}^{N/m} \prod_{\ell=1}^m (1 - \zeta_m^{-\ell}z_jT^{-1})^{\psi_1(\tau^{\ell})a_j},$$ where $v(z_j) = r_{\hub}$ for all
$j$.  To do this, we will need to assume $N_2 \leq mp$.  
This assumption will always be satisfied if $N_1 = (p-1)u_{n-1} - m$.
Now, $\G_{\hub}$ will have the property that 
if $G_{\hub} \in \G_{\hub}$, then the character $\chi$ given rise to
by $G_{\crit}G_{\hub}$ satisfies
$\delta_{\chi}(r_{\hub}) = \delta_{\hub}$, and the left slope of
$\delta_{\chi}$ at $r_{\hub}$ is $u_n$.  This puts us on the right
track for having the right-slope of $\delta_{\chi}$ at $0$ be $u_n$.  Furthermore,
$G_{\crit}G_{\hub}$ will be of the form \eqref{Egnform}.

In \S\ref{Shub2boundary}, in the case $u_n > pu_{n-1}$, we will
construct a particular function $G_{\hub} \in \G_{\hub}$ and modify our original choice of
$G_{\crit} \in \G_{\crit, g}$ such that if $G = G_{\hub}G_{\crit}$ gives rise to $\chi$, then $\delta_{\chi}(0) = 0$, and $\delta_{\chi}$ is linear of
slope $u_n$ on the interval $[0, r_{\hub}]$.  Since $G$ is in the
form \eqref{Egnform}, Proposition \ref{Psufficient} shows that $\chi$ has good reduction
$\chib$ with upper ramification breaks $(u_1, \ldots, u_n)$.

In \S\ref{Sboundary}, still in the case $u_n > pu_{n-1}$, we replace $G$ with $G_n$, where $G_n$ is still a product of
an element of $\G_{\crit, g}$ and one of $\G_{\hub}$ (and thus still of the form \eqref{Egnform}), 
such that $G_n$ gives rise to a $\psi_n$-equivariant character
$\chi_n$ with good reduction $\chib_n$ and $\BB(\chi_n) \subseteq
D(r_n)$ (recall that having good reduction specifically equal to $\chib_n$ is what we seek,
whereas \S\ref{Shub2boundary} only gives us \emph{some} good reduction).  
This is analogous to what happens in the case $u_n = pu_{n-1}$, but a
little more difficult.  In particular, it is tricky to deal with the
coefficient of $t^{-u_n}$ in the last component of the standard form Witt vector
corresponding to $\chib_n$ (no such issue arises in the $u_n =
pu_{n-1}$ case because this coefficient is always zero).  The
underlying calculations concerning this coefficient are deferred to \S\ref{Scalc}.

By Proposition \ref{Pwhenequiv} and the discussion at the beginning of \S\ref{Splan},
we obtain a proof of Theorem \ref{Tsetup} (using Claim \ref{Creform}
instead of Claim \ref{claim1}) in the case $u_n > pu_{n-1}$.  In
\S\ref{Sbranchgeometry}, we summarize the geometry of the branch locus
for the lifts we construct.  

\begin{rem}\label{Rnotoort}
In \cite{OW:ce}, the construction in the case $u_n = pu_{n-1}$ is used
to obtain the proof of lifting in the case $u_n > pu_{n-1}$.  Our
technique here is different, in that it proves both cases independently.
In fact, our method here can
be used to give an alternate proof of \cite[Theorem 1.4]{OW:ce}, and thus (combining with \cite{Po:oc}) of the Oort conjecture.
\end{rem}

\subsection{Controlling $\delta_{\chi}$ between $r_{\crit}$ and $r_{\hub}$}\label{Scrit2hub}
Maintain the previous notation.  In particular, $g$ is a solution to \eqref{Ecrit} realizing the isolated differential data criterion, and 
$\G_{\crit, g}$ is defined as in Definition \ref{Dgcrit}. 

Recall that any $G \in\G_{\crit, g}$
gives rise to a character $\chi$ of order $p^n$ lifting $\chi_{n-1}$. by adjoining the equation
$y_n^p = y_{n-1}G$.  By \eqref{Ecrit} and Proposition \ref{Pdeltalin},
we know that the left derivative of $\delta_{\chi}$ at $r_{\crit}$ is $N_1 + u_{n-1}$.
Recall also from \eqref{Edeltacrit} that $\delta_{\chi}(r_{\crit}) = p/(p-1)$.  

Let $\lambda(G)$ be the minimal 
$\lambda$ in the interval $[r_{\hub}, r_{\crit}]$ such that $\delta_{\chi}(r) = p/(p-1) - (N_1 + u_{n-1})(r_{\crit} - r)$
for all $r\in [\lambda,r_{\crit}]$.  
In other words, $\lambda(G)$ is the largest element in $[r_{\hub},
  r_{\crit}]$ where $\delta_{\chi}$ has left slope less than $N_1 +
  u_{n-1}$ (or is $r_{\hub}$ if there is no such point).
Since $G \in G_{\crit, g}$, we have $\lambda(G) < r_{\crit}$.
Note that $$\delta_{\chi}(\lambda(G)) = \frac{p}{p-1} -
(r_{\crit} - \lambda)(N_1 + u_{n-1}) < \frac{p}{p-1}.$$

\begin{lem}\label{Lreducelambda}
Suppose $G \in \G_{\crit, g}$ with $\lambda := \lambda(G) > r_{\hub}$.  
Identify $\kappa_{\lambda}$ with $k(t)$.  Then $\omega_{\chi}(\lambda)$ can be written
in the form 
\begin{equation}\label{E10}
   \omega_{\chi}(\lambda)=\frac{c\,dt}{t^{N_1 + u_{n-1}+1}} + df,
\end{equation}
where $c \in k^{\times}$ and $f \in t^{1-m}k[t^{-m}]$ has degree less than $N_1 + u_{n-1}$ in $t^{-1}$.
\end{lem}

\proof
%
If $u_n = pu_{n-1}$ then the same argument as in the proof of
  \cite[Proposition 6.13]{OW:ce} shows that $\omega_{\chi}(\lambda)$
  is as in \eqref{E10}
for some $c \in k^{\times}$ and $f\in\kappa_\lambda$, with  
$f$ a polynomial in $t^{-1}$ of
degree $< N_1 + u_{n-1}$ and without constant term.  

If $u_n > pu_{n-1}$, so that $N_1 < (p-1)u_{n-1}$ (see the beginning
of \S\ref{Sstrategy}), then $\delta_{\chi}(\lambda) > p
\delta_{\chi_{n-1}}(\lambda) = \lambda pu_{n-1}$ (the equality is due
to \cite[Lemma 6.1]{OW:ce}).
  Thus $\C(\omega_{\chi}(\lambda)) = 0$ (also from \cite[Proposition
  4.3(ii)]{We:fr}) and $p \nmid N_1 + u_{n-1}$.  Since the
  differential form in \eqref{E10} is exact in this case, \eqref{E10}
  holds as well, with the same conditions on the terms of $f$. 

By Lemmas
\ref{Lwhichroot} and \ref{Lequivswan}, the
$\tau$-equivariance of $\chi$ implies that $\tau(df) = \zeta_m df$, where $\tau$ acts on $t$ and $dt$ by multiplication by $\zeta_m$.
That is, we may assume that $f$ only has terms of degree $t^{-q}$ where $q \equiv -1 \pmod{m}$. 
\Endproof

The following proposition is crucial, and will be proved in
\S\ref{Sproofs}.

\begin{prop}[cf.\ \cite{OW:ce}, Corollary 7.5]\label{Preduce} 
Let $G \in \G_{\crit, g}$, let $r \in [r_{\hub}, r_{\crit}) \cap \rats$, and let $f \in t^{1-m}k[t^{-m}]$ be a polynomial of degree less than $N_1 + u_{n-1}$ in 
$t^{-1}$, which we regard as the reduction of $T_r$ in $\kappa_r$ (\S\ref{Sgeomsetup}).  Assume $f$ has no terms of degree divisible by $p$.  
Let $\beta = (N_1 + u_{n-1})(r_{\crit} - r)$.  After a possible finite extension of $K$, 
there exist $G' \in \G_{\crit, g}$ and $F \in \KK$ with $v_r(F) = 0$ and $[F]_r = f$ 
such that 
$$\frac{G'}{G} \equiv 1 - p^{\beta}F \pmod{(\KK^{\times})^{p}}.$$
\end{prop}

We now show that $\lambda$ can be reduced.
\begin{prop}[cf.\ \cite{OW:ce}, Proposition 6.13]\label{Preducelambda} Suppose $G \in\G_{\crit, g}$ with $\lambda(G)> r_{\hub}$. Then there
  exists $G'\in\G_{\crit, g}$ with $\lambda(G')<\lambda(G)$.
\end{prop}

\proof
This follows from Lemma \ref{Lreducelambda} and Proposition
\ref{Preduce} exactly as in the proof of \cite[Proposition
6.13]{OW:ce} with $N_1 + u_{n-1}$ playing the role of $m_n$ there.
\Endproof


\begin{prop}[cf.\ \cite{OW:ce}, Proposition 6.15]\label{Pkinkmin} 
  The function $G\mapsto \lambda(G)$ takes a minimal value on $\G_{\crit, g}$.
\end{prop}

\proof 
Recall that $g =  \prod_{j=1}^{N_1/m} \prod_{\ell=1}^m (1 -
\zeta_m^{-\ell}\xb_jt^{-1})^{\psi_1(\tau^{\ell})a_j}$, with the
$\psi_1(\tau^{\ell})a_j$ viewed as lying in $\ints$ \eqref{Edefg}.  
Let $U' \subseteq \aff^{N_1/m}(k)$ be the open subset consisting of those
$(\yb_1, \ldots, \yb_{N_1/m})$ such $\yb_i^m \neq \yb_j^m$ if $i \neq
j$, and let $V \subset \aff^{N_1/m}_k$ be the subvariety such that the 
$\yb_j$ and $a_j$ give a solution to \eqref{Ecrit}.  Since $g$ realizes
the isolated differential data criterion, the point $\xb = (\xb_1, \ldots,
\xb_{N_1/m})$ is an isolated point of $V$.  In particular, $V' := V
\backslash \{\xb\} \subseteq \aff^{N_1/m}_k$ is closed and $U =
U' \backslash V' \subseteq \aff^{N_1/m}_k$ is open.
 
Let $\G'_{\crit,g} \supseteq \G_{\crit,g}$ be the set of $G \in \KK$ such that
$$G = \prod_{j=1}^{N_1/m} \prod_{\ell=1}^m (1 - \zeta_m^{-\ell}z_jT^{-1})^{\psi_1(\tau^{\ell})a_j},$$
where if $y_j = z_jp^{-r_{\crit}}$, then the reductions $\yb_j$ give a
point $(\yb_1, \ldots, \yb_{N_1/m}) \in U$.  By identifying each $G
\in \G'_{\crit, g}$ with $(y_1, \ldots, y_{N_1/m})$, we identify $\G'_{\crit,
  g}$ with the rigid analytic space 
$$U^{\text{rig}} := \{y \in (\aff^{N_1/m})^{\text{an}} \ | \ \yb \in U \},$$
where $\yb$ is the canonical reduction of $y$.  Since $U$ is open,
$U^{\text{rig}}$ is a finite union of open
affinoid subdomains of $(\aff^{N_1/m})^{\text{an}}$.  In particular,
it is quasi-compact and quasi-separated.

Extend the domain of $\lambda$ from $\G_{\crit, g}$ to $\G'_{\crit,g}$,
keeping the definition the same.  The family of $\ints/p^n$-covers of $\proj^1_K$
parameterized by $U$ via taking the Kummer extensions given rise to by
points in $U$ is a good relative Galois cover in the language of $\cite[\S5]{OW:wr}$.
By \cite[Corollary 5.3(ii)]{OW:wr} (taking $r_0 = r_{\crit}$ and
$m_{\text{Swan}} = N_1 + u_{n-1} + 1$ in the notation of that
corollary), $\lambda$ ($ = \lambda_{\text{Swan}}$) achieves its minimal value on
$\G'_{\crit, g}$ after a possible extension of $K$.  
On the other hand, our construction of $\G'_{\crit,g}$ shows that $\lambda(G) =
r_{\crit}$ if $G \in \G'_{\crit, g} \backslash \G_{\crit, g}$ and $\lambda(G) < r_{\crit}$ if $G \in \G_{\crit, g}$.  
Thus our minimal value must achieved on $\G_{\crit, g}$.  
\Endproof

\begin{cor}\label{Cgoodred}
\begin{enumerate}
\item There exists $G_{\crit} \in \G_{\crit,g}$ giving rise to a
character $\chi_{\crit}$ such that $\delta_{\chi_{\crit}}(r_{\hub}) =
\delta_{\hub}$.
\item
If $u_n = pu_{n-1}$, then $\chi_{\crit}$ has good reduction.
\end{enumerate}
\end{cor}

\proof
Proposition \ref{Preducelambda}, combined with Proposition \ref{Pkinkmin}, shows that there exists $G_{\crit} \in \G_{\crit,g}$ giving rise to a
character $\chi_{\crit}$ such that $\lambda(\chi_{\crit}) = r_{\hub}$.
In other words, $\delta_{\chi_{\crit}}(r_{\hub}) = \delta_{\hub}$.
This proves (i).  

If $u_n = pu_{n-1}$, then $r_{\hub} = 0$, and $\lambda(G_{\crit}) =
0$.  That is, $\delta_{\chi_{\crit}}$ is linear of slope $N_1 + u_{n-1} = pu_{n-1}$ on the 
interval $[0, r_{\crit}]$, with $\delta_{\chi_{\crit}}(0) = 0$.  Part
(ii) then follows from Proposition \ref{Psufficient}.
\Endproof

\begin{defn}\label{Dconstant}
We let $C \in k^{\times}$ be the coefficient of $t^{-(N_1 + u_{n-1} + 1)}dt$ in
$\omega_{\chi_{\crit}}(r_{\hub})$.  
\end{defn}

In fact, $C$ is independent of the choice of $G_{\crit} \in \G_{\crit,
  g}$.  This is the statement of Proposition \ref{PCconstant}, whose
proof will be deferred to \S\ref{Scalc}.

\subsection{Controlling $\omega_{\chi}$ at $r_{\hub}$}\label{Shub}
The material in this section is only necessary if $u_n > pu_{n-1}$. So we now assume that
$pu_{n-1} < u_n < pu_{n-1} + mp$ (this is the assumption of no essential ramification).  
Recall that this means that $p \nmid u_n$, that $(p-1)u_{n-1} - mp < N_1 < (p-1)u_{n-1}$, and that
$N_2 = u_n - u_{n-1} - N_1 \leq 2mp - 2m$.  Throughout this section,
let $s = (N_1+u_{n-1})(r_{\crit} - r_{\hub})$, and let $t =
[T]_{r_{\hub}}$.  We have constructed a rational function $G_{\crit}
\in \G_{\crit, g}$ (Corollary \ref{Cgoodred}) giving
rise to a character $\chi_{\crit}$ such that (see \eqref{Edhub}) 
\begin{equation}\label{Es}
\delta_{\chi_{\crit}}(r_{\hub}) = \delta_{\hub} = u_n r_{\hub} =
\frac{p}{p-1} - (N_1 + u_{n-1})(r_{\crit} - r_{\hub}) = \frac{p}{p-1}
- s.
\end{equation}
Let $C \in k^{\times}$ be the constant from Definition
\ref{Dconstant}.

For the rest of this section, we will make a further assumption.

\begin{ass}\label{AN1N2}
$N_2 \leq mp$.
\end{ass}
Note that Assumption \ref{AN1N2} is always satisfied when $N_1 = (p-1)u_{n-1} - m$ (and that for any other choice of $N_1$, there will be values of
$u_n$ leading to a violation of Assumption \ref{AN1N2}).

As mentioned in \S\ref{Sstrategy}, our eventual goal is (after
possibly modifying $G_{\crit}$), to construct a rational function $G_{\hub} \in \KK$ with $N_2$ zeroes and poles away from 
$T = 0$, all of which have
valuation $r_{\hub}$, so that if we let $G_n = G_{\hub}G_{\crit}$, then $G_n$ gives rise to a character $\chi_n$
with good reduction $\chib_n$.  In this case, $\delta_{\chi_n}$ would be linear of slope $u_n$ on the interval $[0, r_{\hub}]$ and linear of slope
$N_1 + u_{n-1}$ on the interval $[r_{\hub}, r_{\crit}]$.  In particular, the differential form 
$\omega_{\chi_n}(r_{\hub})$ would have to have a zero of order $u_n - 1$ at $t=\infty$ and a pole of order $N_1 + u_{n-1} + 1$ at $t=0$ (Proposition
\ref{Pdeltalin}).  By Proposition \ref{Pvccor2}, there can be no zeroes away from $\infty$.  By Lemma \ref{Lequivswan}, $\omega_{\chi_n}(r_{\hub})$
must transform equivariantly under $t \mapsto \zeta_m t$.  So we will
search for $G_{\hub}$ such that $G_n := G_{\hub}G_{\crit}$ can give
rise to a character $\chi_n$ such that
$\delta_{\chi_n}(r_{\hub})$ is still $\delta_{\hub}$ and 

\begin{equation}\label{Eomegahub}
\omega_{\chi_n}(r_{\hub}) = \frac{c\, dt}{(t^m - \alphab^m)^{N_2/m}t^{N_1 + u_{n-1} + 1}},
\end{equation}
where $c$ and $\alphab$ are in $k^{\times}$.  Proposition
\ref{Phubstart} will show how a valid $G_{\hub}$ arises.

\begin{lem}\label{Lomegaexact}
The differential form on the right hand side of \eqref{Eomegahub} has a zero of order $u_n - 1$ at $\infty$, 
a pole of order $N_1 + u_{n-1} + 1$ at $0$, no zeroes away
from $\infty$, transforms as in Lemma \ref{Lequivswan}, and is exact for all choices of $c$ and $\alphab$.
\end{lem}

\proof
Once we note that $m | (N_1 + u_{n-1} + 1)$, all assertions become trivial except the last one. 
Multiplying a differential form by a $p$th power does not change its exactness, so it suffices to show that
$$\frac{(t^m - \alphab^m)^{p - N_2/m}dt}{t^{N_1 + u_{n-1} + 1}}$$ is exact.  By Assumption \ref{AN1N2}, the numerator is a polynomial in 
$t$. Expanding everything out, $t$ occurs to degrees
$-(N_1 + u_{n-1} + 1)$ through $mp - u_n - 1,$ counting by $m$'s.
Since $-(N_1 + u_{n-1} + 1) > -pu_{n-1} -1$ and $mp - u_n - 1 <
-pu_{n-1} + mp - 1$, and since all the degrees in question are
divisible by $m$, we see that 
none of the above degrees is congruent to $-1 \pmod{p}$.  This means that the differential form is exact.
\Endproof

\begin{lem}\label{Ldiffexact}
Let $\omega$ be the differential form on the right hand side of \eqref{Eomegahub}.  Then
$\omega - \omega_{\chi_{\crit}}(r_{\hub})$ is exact.  Furthermore, we can write
\begin{equation}\label{Ehubstart}
\omega - \omega_{\chi_{\crit}}(r_{\hub}) = \left(\frac{a(t)}{(t^m - \alphab^m)^{N_2/m}} + \frac{b(t)}{t^{N_1 + u_{n-1} + 1}}\right)dt,
\end{equation}
where both fractions are proper, $a(t)$ and $b(t)$ are in $k[t^m]$,
and each of the two summands is exact.
Lastly, choosing $c =
C(-\alphab^m)^{N_2/m}$ on the right hand side of \eqref{Eomegahub} results in $b(t)$ having no
constant term.
\end{lem}

\proof
The exactness of $\omega - \omega_{\chi_{\crit}}(r_{\hub})$
follows from Lemmas \ref{Lreducelambda} and \ref{Lomegaexact}, noting that $p \nmid N_1 + u_{n-1}$.
Since $\omega_{\chi_{\crit}}(r_{\hub})/dt$ is a proper fraction in $t$
with denominator $t^{N_1 + u_{n-1} + 1}$ and $\alphab \neq 0$, 
the theory of partial fractions gives the desired decomposition into the two summands.  The polynomials $a(t)$ and $b(t)$ lie in 
$k[t^m]$ by Lemma \ref{Lequivswan}, combined with the fact that $\chi_{\crit}$ is $\tau$-equivariant.  Each of the two summands is
exact because their sum is, and a sum of two proper fractions (times $dt$) with relatively prime denominators can only be exact if each one is.
Lastly, in order for $b(t)$ not to have a constant term, we need only ensure that when $\omega - \omega_{\chi_{\crit}}(r_{\hub})$ is written as
$f(t)dt/(t^m - \alphab^m)^{N_2/m} t^{N_1 + u_{n-1} + 1}$, that $f(t)$ has no constant term.  This is accomplished by taking
$c = C(-\alphab^m)^{N_2/m}$, where $C$ is the coefficient of the
$dt/t^{N_1 + u_{n-1} +1}$ term of $\omega_{\chi_{\crit}}(r_{\hub})$.
\Endproof

\begin{prop}\label{Pcritreplace}
Let $\omega$ be the right hand side of \eqref{Eomegahub}, with $c$
chosen as in Lemma \ref{Ldiffexact} and $C$ chosen as in Definition \ref{Dconstant}. 
By modifying $G_{\crit}$ within $\G_{\crit, g}$, 
we can ensure $$\omega - \omega_{\chi_{\crit}}(r_{\hub}) = \frac{a(t) dt}{(t^m - \alphab^m)^{N_2/m}},$$ with $a(t)$ as in Lemma \ref{Ldiffexact}.
\end{prop}

\proof
By our choice of $c$, we may assume that $b(t)$ has no constant term in the notation of Lemma \ref{Ldiffexact}. Thus we can write
$$\frac{b(t)dt}{t^{N_1 + u_{n-1} + 1}} = df,$$ where
$f \in t^{1-m}[t^{-m}]$ has degree less than $N_1 + u_{n-1}$ in
$t^{-1}$.  By Proposition \ref{Preduce}, there exists $G_{\crit}' \in\G_{\crit, g}$ such that
\[
    \frac{G_{\crit}'}{G_{\crit}} \equiv 1+p^sF \pmod{(\KK^\times)^p},
\]
where $v_{\lambda}(F) = 0$, where $[F]_{\lambda} = f$.  As in the
proof of \cite[Proposition 6.13]{OW:ce}, replacing $G_{\crit}$ by $G_{\crit}'$ has the effect of adding $df$ to $\omega_{\chi_{\crit}}(r_{\hub})$, which in turn has the
effect of subtracting $df$ from the right hand side of \eqref{Ehubstart}.  This proves the proposition.
\Endproof

By Proposition \ref{Pcritreplace} we may, and do, assume that $b(t) = 0$ in \eqref{Ehubstart}.  We do a further partial fractions decomposition 
on the other term to obtain
$$\frac{a(t) dt}{(t^m-\alphab^m)^{N_2/m}} = \sum_{\ell=0}^{m-1} \frac{\beta_{\ell}(t)dt}{(\zeta_{m}^{\ell}t - \alphab)^{N_2/m}},$$ where 
$\beta_{\ell}(t)$ is a polynomial of degree less than $N_2/m$.  Using equivariance under 
$t \mapsto \zeta_m t$, it is not hard to check that $\beta_{\ell}(t) = \beta_0(\zeta_m^{\ell}t)$ for all $\ell$.  

The following definition is the key idea of \S\ref{Shub}.                                
\begin{defn}\label{Dghub}
Let $\alphab \in k^{\times}$.
Let $T_{\hub} = p^{-r_{\hub}}T$ and $s = p/(p-1) - \delta_{\hub}$.
For any lift $\alpha$ of $\alphab$ to $R$ and fixed lifts of
$\zeta_m^{\ell}$ from $\FF_p$ to $\ints$ (denoted again by
$\zeta_m^{\ell}$ by abuse of notation),
let $\G_{\hub, \alpha} \subseteq \KK^{\times}$ be
the set of all rational functions 
of the form
$$\prod_{\ell=0}^{m-1} \left(1 + p^sA(\zeta_m^{\ell}T_{\hub})\right)^{\zeta_m^{-\ell}},$$
such that 

\begin{equation}\label{Ehubform}
A(T_{\hub}) \text{ is of the form }
\frac{B(T_{\hub})}{(T_{\hub} - \alpha)^{N_2/m - 1}},
\end{equation}
where $B(T_{\hub})$ is a polynomial of degree at most $N_2/m - 2$ and $v_{r_{\hub}}(A) =
0$ with $d([A]_{r_{\hub}})/dt = \beta_0(t)/(t - \alphab)^{N_2/m}$.
\end{defn}

\begin{defn}\label{Dlargehub}
Let $\G_{\hub} = \bigcup_{\alpha \in R^{\times}}
\G_{\hub, \alpha}$.
\end{defn}

\begin{rem}\label{Rghubnonempty}
By Lemma \ref{Ldiffexact}, $a(t) dt/(t^m-\alphab^m)^{N_2/m}$ is exact.  It is straightforward then to show that
$\beta_0(t) dt/(t - \alphab)^{N_2/m}$ is also exact.  In particular,
$\G_{\hub, \alpha}$ is nonempty for all $\alpha \in R^{\times}$.  
\end{rem}

The following definition will be useful in the proof of Proposition \ref{Pcrithubkinkmin}.
\begin{defn}\label{Dghubbig}
For $\alpha \in R^{\times}$, define $\G'_{\hub, \alpha} \supseteq
\G_{\hub, \alpha}$ exactly as in Definition \ref{Dghub}, except that we
impose the condition $v_{r_{\hub}}(A) \geq 0$ instead of
$v_{r_{\hub}}(A) = 0$, and we place no condition on $d([A]_{r_{\hub}})/dt$.
\end{defn} 

We prove the major result of this section.

\begin{prop}\label{Phubstart}
If $G_{\hub} \in \G_{\hub, \alpha}$ and $G_{\crit}$ is chosen as in
Proposition \ref{Pcritreplace}, then $G_{\crit}
G_{\hub}$ gives rise to a character $\chi$ such that
$\delta_{\chi}(r_{\hub}) = \delta_{\hub}$, and $\omega :=
\omega_{\chi}(r_{\hub})$ is the right hand side of \eqref{Eomegahub}
with $c$ chosen as in Lemma \ref{Ldiffexact}.  
Consequently, the left slope of $\delta_{\chi}$ at $r_{\hub}$ is
$u_n$.  

Furthermore, $G_{\crit}G_{\hub}$ is of the form \eqref{Egnform}.
\end{prop}				

\proof
The product $$\prod_{\ell=0}^{m-1} \left(1 + p^sA(\zeta_m^{\ell}T_{\hub})\right)^{\zeta_m^{-\ell}}$$ can be written as
$$1 + p^s \sum_{\ell = 0}^{m-1} \zeta_m^{-\ell}A(\zeta_m^{\ell}T_{\hub}) + D,$$ where $v_{r_{\hub}}(D) > s$.
By the definition of $A$, the derivative of $\sum_{\ell = 0}^{m-1} [\zeta_m^{-\ell}A(\zeta_m^{\ell}T_{\hub})]_{r_{\hub}}$
is $$\sum_{\ell = 0}^{m-1} \frac{\beta_0(\zeta_m^{\ell}t)}{(\zeta_m^{\ell}t - \alphab)^{N_2/m}} = 
\sum_{\ell= 0}^{m-1} \frac{\beta_{\ell}(t)}{(\zeta_{m}^{\ell}t - \alphab)^{N_2/m}} = \frac{a(t) dt}{(t^m - \alphab^m)^{N_2/m}}.$$
By Propositions \ref{Paddchar} and \ref{Pcritreplace}, we get that $\delta_{\chi}(r_{\hub}) = \delta_{\hub}$ and $\omega_{\chi}(r_{\hub}) = \omega$. 

Since $\omega$ has a zero of order $N_2 + N_1 + u_{n-1} - 1 = u_n -1$ at $\infty$, Proposition \ref{Pdeltalin} shows that the left slope of 
$\delta_{\chi}$ at $r_{\hub}$ is $u_n$.

Showing that $G_{\hub}G_{\crit}$ is of the form \eqref{Egnform} is equivalent to showing that $1 + p^sA(T_{\hub})$ is, up to multiplication by a $p$th power,
a polynomial in $T^{-1}$ with constant term $1$ and at most $N_2/m$ distinct roots.  In order to do this, we multiply $1 + p^sA(T_{\hub})$ by 
$(T_{\hub}-\alpha)^p/T_{\hub}^p$, and we leave it to the reader to verify
that everything works (the roots will be the $N_2/m - 1$ roots of $(T_{\hub} - \alpha)^{N_2/m-1} + p^sB(T_{\hub})$, along with $\alpha$).   
\Endproof

\begin{rem}\label{Ronlysmall}
Let $G_{\hub} \in \G'_{\hub, \alpha}
\backslash \G_{\hub, \alpha}$ and $G_{\crit} \in \G_{\crit, g}$.  By
the discussion above, if $G_{\crit}G_{\hub}$ gives rise
to a character $\chi$, then it is \emph{not} the case that both
$\delta_{\chi}(r_{\hub}) = \delta_{\hub}$ and the left slope of $\delta_{\chi}$
at $r_{\hub}$ is $u_n$.
\end{rem}

\subsection{Controlling $\delta_{\chi}$ between $r_{\hub}$ and $0$}\label{Shub2boundary}
We maintain the assumption of \S\ref{Shub} that $pu_{n-1} < u_n <
pu_{n-1} + mp$, as well as all the notation so far.
Fix $\alpha \in R^{\times}$ and $g$ a solution to \eqref{Ecrit}
realizing the isolated differential data criterion.  By Proposition \ref{Phubstart}, there exists $G_{\crit} \in \G_{\crit, g}$ such that, for any $G_{\hub} \in \G_{\hub, \alpha}$, 
the character $\chi$ given rise to by $G_{\hub}G_{\crit}$ has $\delta_{\chi}(r_{\hub}) = \delta_{\hub}$ and $\delta_{\chi}$ 
has a left slope of $u_n$ at $r_{\hub}$.  The goal of this section
is to find a particular $\tilde{G}_{\crit} \in \G_{\crit, g}$ and $\tilde{G}_{\hub} \in \G_{\hub, \alpha}$ such that $\tilde{G}_{\hub}\tilde{G}_{\crit}$ 
gives rise to a character $\chi$ with $\delta_{\chi}(0) = 0$.  
Since $\delta_{\hub} = u_n r_{\hub}$ by definition, one can test this by seeing if $\delta_{\chi}$ is
linear of slope $u_n$ on the interval $[0, r_{\hub}]$.  Let $\G_{g, \alpha}$ be the (nonempty) subset of $\G_{\hub, \alpha} \G_{\crit,g}$ consisting of elements giving 
rise to characters $\chi$ with $\delta_{\chi}(r_{\hub}) = \delta_{\hub}$ and such that the left slope of $\delta_{\chi}$ at $r_{\hub}$ is $u_n$.  
Note that every $G \in \G_{g, \alpha}$ is of the form \eqref{Egnform}. 
If $G \in \G_{g, \alpha}$ gives rise
to $\chi$, then we define $\mu(G)$ to be the minimal element of
$[0, r_{\hub}]$ such that $\delta_{\chi}(\mu(G)) = u_n\mu(G)$ (that
is, $\mu(G)$ is the largest element of $[0, r_{\hub}]$ where $\delta_{\chi}$ has
left slope less than $u_n$, or $0$, if no such element exists).  This
is analogous to the definition of $\lambda(G)$ in \S\ref{Scrit2hub}.
For any $G \in G_{g, \alpha}$, we have $\mu(G) < r_{\hub}$.  The goal of this section is to prove the existence of
$G \in \G_{g, \alpha}$ such that $\mu(G) = 0$.  Then $G$ will give rise
to a character with good reduction and upper jumps $(u_1, \ldots,
u_n)$.  The argument is parallel to that of \S\ref{Scrit2hub}.

\begin{lem}\label{Lreducemu}
Suppose $G \in \G_{g, \alpha}$ with $\mu := \mu(G) > 0$.  
Identify $\kappa_{\mu}$ with $k(t)$.  Then $\omega_{\chi}(\mu)$ can be written
in the form $$\frac{c\,dt}{t^{u_n+1}} + df,$$
where $c \in k^{\times}$ and $f \in t^{1-m}k[t^{-m}]$ has degree less than $u_n$ in $t^{-1}$.
\end{lem}

\proof
After noting that
  $$\delta_{\chi}(\mu) = p/(p-1) - s - u_n(r_{\hub} - \mu) = \mu u_n >
  \mu pu_{n-1},$$ where $s$ is as in Definition \ref{Dghub}, the proof is exactly the same as the $u_n > pu_{n-1}$
  case of Lemma \ref{Lreducelambda}.
%
\Endproof

As in \S\ref{Scrit2hub}, we also postpone the proof of the following
crucial result 
to \S\ref{Sproofs}.
 
\begin{prop}\label{Pcrithubreduce} 
Suppose $N_1 = (p-1)u_{n-1} - m$ (this is consistent with Assumption
\ref{AN1N2}).
Let $G_{\crit}, G_{\hub} \in G_{\crit, g}, \G_{\hub, \alpha}$, respectively. 
Let $r \in [0, r_{\hub}) \cap \QQ$, and let $f \in t^{1-m}k[t^{-m}]$ have degree less than $u_n$ in 
$t^{-1}$, which we regard as the reduction of $T_r$ in $\kappa_r$ (\S\ref{Sgeomsetup}).  Assume $f$ has no terms of degree divisible by $p$.  
Let $\beta = p/(p-1) - u_nr$.  After a possible finite extension of $K$, 
there exist $G_{\crit}', G_{\hub}' \in \G_{\crit, g}, \G_{\hub, \alpha}$
respectively, and $F \in \KK$ with $v_r(F) = 0$ and $[F]_r = f$ 
such that 
$$\frac{G_{\crit}'G_{\hub}'}{G_{\crit}G_{\hub}} \equiv 1 - p^{\beta}F \pmod{(\KK^{\times})^{p}}.$$
\end{prop}

\begin{rem}
The proofs of Propositions \ref{Preduce} and \ref{Pcrithubreduce} are
the only places where the \emph{isolatedness} in the isolated
differential data criterion is used.
\end{rem}  

This has the following consequence:

\begin{prop}\label{Preducemu} Suppose $G \in \G_{g, \alpha}$ with $\mu(G)> 0$. Then there
  exists $G' \in\G_{g, \alpha}$ with $\mu(G')<\mu(G)$.
\end{prop}

\proof
The proof is exactly the same as the proof of Proposition
\ref{Preducelambda}, using Lemma
\ref{Lreducemu} and Proposition \ref{Pcrithubreduce} in place of Lemma
\ref{Lreducelambda} and Proposition \ref{Preduce}, and taking $\beta = p/(p-1) - u_n\mu$ instead of
$\beta = (N_1 + u_{n-1})(r_{\crit} - \lambda)$.
\Endproof

\begin{prop}\label{Pcrithubkinkmin} 
  The function $G \mapsto \mu(G)$ takes a minimal value on $\G_{g, \alpha}$.
\end{prop}

\proof 
We identify $\G'_{\hub, \alpha}$ with the rigid $(N_2/m -
1)$-dimensional closed unit polydisc corresponding to the coefficients
of $A$ in Definition \ref{Dghubbig}.  This is an affinoid space.  
Let $\G'_{\crit, g}$ be as in the proof of Proposition
\ref{Pkinkmin}, and let $\G'_{g, \alpha} = 
\G'_{\crit, g}\G'_{\hub, \alpha} \subseteq \KK$.  It is easy to see that $\G'_{g, \alpha}
\cong \G'_{\crit, g} \times \G'_{\hub, \alpha}$, and is thus
identified with a quasi-compact, quasi-separated rigid-analytic space.

Extend the domain of $\mu$ from $\G_{g, \alpha}$ to $\G'_{g, \alpha}$,
keeping the definition the same.
The family of $\ints/p^n$-covers of $\proj^1_K$
parameterized by $\G'_{g, \alpha}$ via taking the Kummer extensions given rise to by
points in $\G_{g, \alpha}$ is a good relative Galois cover in the language of $\cite[\S5]{OW:wr}$. 
By \cite[Corollary 5.3(ii)]{OW:wr}, $\mu$ achieves its minimal value on
$\G'_{g, \alpha}$, after a possible extension of $K$.  
On the other hand, suppose $G = G_{\crit}G_{\hub} \in \G'_{g, \alpha}
\backslash \G_{g, \alpha}$ with $G_{\crit} \in \G'_{\crit, g}$ and
$G_{\hub} \in \G'_{\hub, \alpha}$.  We claim that the
left-slope of $\delta_{\chi}$ at $r_{\hub}$ is less than $u_n$, which
means that $\mu(G) = r_{\hub}$.  Since $\mu(G) < r_{\hub}$ when $G \in
\G_{g, \alpha}$, this means that the minimal value of $\mu$ on
$G'_{g, \alpha}$ must be achieved on $\G_{g, \alpha}$, thus completing the proof.  

To prove the claim, first assume that $\delta_{\chi}(r_{\hub}) \neq
\delta_{\hub}$.  Then $\delta_{\chi}(r_{\hub}) > \delta_{\hub}$.
Since $\delta_{\chi}$ is concave up on $[r_{\hub}, r_{\crit}]$
(combine Propositions \ref{Pdeltalin} and \ref{Pvccor2} with the fact
that $G$ has no zeroes or poles with valuation in $(r_{\hub}, r_{\crit})$), we
have that the right-slope of 
$\delta_{\chi}$ at $r_{\hub}$ is less than $N_1 + u_{n-1}$. 
Since $G_{\hub}$ has at most $N_2$
zeroes and poles with valuation $r_{\hub}$, the left-slope of
$\delta_{\chi}$ at $r_{\hub}$ is less than $N_1 + N_2 + u_{n-1} = u_n$
at $r_{\hub}$ (again, combine Propositions \ref{Pdeltalin} and \ref{Pvccor2}).

Now, assume $\delta_{\chi}(r_{\hub}) = \delta_{\hub}$.  Then the left-slope of
$\delta_{\chi}$ at $r_{\crit}$ is $N_1 + u_{n-1}$, and $G_{\crit} \in
\G_{\crit, g}$.  Remark \ref{Ronlysmall} shows that, if
$G_{\hub} \in \G'_{\hub, \alpha} \backslash \G_{\hub, \alpha}$, then the left-slope of
$\delta_{\chi}$ at $r_{\hub}$ is less than $u_n$.  So assume $G_{\hub}
\in \G_{\hub, \alpha}$.  The definition of $G_{g, \alpha}$ shows
that if $G \in G'_{g, \alpha} \backslash G_{g, \alpha}$, then the
left-slope of $\delta_{\chi}$ at $r_{\hub}$ is less than $u_n$.  The
claim, and thus the proposition, is proved.
\Endproof

\begin{cor}\label{Ccrithubgoodred}
\begin{enumerate}
\item Suppose $u_n > pu_{n-1}$.  Then there exists $G \in \G_{g, \alpha}$ giving rise to a
character $\chi$ having good reduction.
\end{enumerate}
\end{cor}

\proof
Proposition \ref{Preducemu}, combined with Proposition
\ref{Pcrithubkinkmin}, shows that there exists $G \in \G_{g, \alpha}$ giving rise to a
character $\chi$ such that $\mu(\chi) = 0$.  That is,
$\delta_{\chi}(0) = 0$.  The corollary then follows from Proposition \ref{Psufficient}.
\Endproof

\subsection{Ensuring the correct reduction on the boundary}\label{Sboundary}
In this section, we prove Theorem \ref{Tmincase}, which will complete
the proof of Theorem \ref{Tsetup}.  Maintain all notation from the previous sections, including $C$ as the
constant from Definition \ref{Dconstant}. First, we prove two lemmas.

\begin{lem}\label{Lartinschreier}
If $G_1 \in \G_{g, \alpha_1}$ gives rise to a character $\chi$ with good
reduction $\chib$ that corresponds (after completion at $t=0$) to the Witt
vector $(f_1, \ldots, f_{\chi}) \in W_n(k((t)))$, and
$G_2 \in \KK$ is such that $G_2/G_1 \equiv 1 + p^{p/(p-1)}F \pmod{(\KK^{\times})^p},$
with $v_0(F) = 0$ and $[F]_0 = f$ in $k(t) \subseteq k((t))$, then
$G_2$ gives rise to a character $\chi'$ with \'{e}tale reduction $\chib'$
that corresponds (after completion at $t=0$) to the Witt vector $(f_1,
\ldots, f_{n-1}, f_{\chi} + f) \in W_n(k((t)))$.
\end{lem}

\proof
Replacing $G_1$ by $G_2$ has the effect of multiplying $\chi$ by
$$\psi_n := \K_n((G_2/G_1)^{p^{n-1}}) \in H^1_{p^n}(\KK).$$
This is just the image of $\psi_1 := \K_1(G_2/G_1) \in H^1_p(\KK)$.
Proposition \ref{Pdelta} shows that $\delta_{\psi_1}(0) = 0$ and the reduction $\psib_1$ corresponds 
to the Artin-Schreier extension given by $y^p - y = f$.  
Consequently, $\delta_{\psi_n}(0) = 0$ and its reduction $\psib_n$
corresponds to the extension encoded by the Witt vector $(0, \ldots, 0, f)$.  

By Proposition \ref{Paddchar}(iii), we conclude that the reduction of
$\chib' = \chib\psib_n$ corresponds to the sum of the Witt vectors
$(f_1, \ldots, f_{\chi})$ and $(0, \ldots, 0, f)$.  This is $(f_1,
\ldots, f_{n-1}, f_{\chi}  + f)$, as desired.
\Endproof

The proof of the following lemma relies on Lemma \ref{Luncoeff},
proven in \S\ref{Scalc}.
\begin{lem}\label{Luncoeff2}
Suppose $u_n > pu_{n-1}$, and let $\alpha_1 \in R^{\times}$ with
reduction $\alphab_1 \in k^{\times}$.  Suppose
$G_1 \in \G_{g, \alpha_1}$ gives rise to a character
$\chi$ with good reduction $\chib$.  Then if $(f_1, \ldots, f_{n-1},
f_{\chi})$ is the Witt vector in standard form corresponding to
$\chib$ (after completion at $t=0$), the $t^{-u_n}$ coefficient of $f_{\chi}$ is
$-Cu_n^{-1}(-\alphab_1^m)^{N_2/m}$.
\end{lem}

\proof
Let $\gamma$ be the (nonzero) $t^{-u_n}$ coefficient of $f_{\chi}$, and assume
for a contradiction that $\gamma \neq -Cu_n^{-1}(-\alphab_1^m)^{N_2/m}$.
Let $\alphab_2$ be such that
$$Cu_n^{-1}((-\alphab_1^m)^{N_2/m} - (-\alphab_2^m)^{N_2/m}) = -\gamma.$$
By our assumption, $\alphab_2 \neq 0$.  Let $\alpha_2 \in R$
be a lift of $\alphab_2$.  Using Corollary
\ref{Ccrithubgoodred}, choose $G_2 \in \G_{\alpha_2}$
giving rise to a character $\chi'$ with good reduction $\chib'$.
Since $\gamma \neq 0$, we have $\alpha_1^{N_2} \neq \alpha_2^{N_2}$,
so Lemma \ref{Luncoeff} applies.  In particular, Lemma \ref{Luncoeff}
implies that 
$$G_2/G_1 \equiv 1 + p^{p/(p-1)}F \pmod{(\KK^{\times})^p},$$ where $F
  \in \KK$ satisfies $v_0(F) = 0$ and $[F]_0$ is a polynomial in
  $t^{-1} = [T^{-1}]_{0}$ of degree $u_n$ with leading term $-\gamma t^{-u_n}$.
By Lemma \ref{Lartinschreier}, replacing $G_1$ with $G_2$ replaces $f_{\chi}$ in the
Witt vector for $\chib$ with $f_{\chi} + [F]_0$, which has degree less than $u_n$ in
$t^{-1}$.  This means that the $n$th higher ramification jump for the
upper numbering of $\chib'$ is less than $u_n$ (\S\ref{Switt}), which
contradicts Proposition \ref{Pvccor1}.  
\Endproof

If $u_n = pu_{n-1}$, let $\G = \G_{\crit,g}$.  Otherwise, let $\G =
\bigcup_{\alpha \in R^{\times}} \G_{g, \alpha}$.  Note that all
elements of $\G$ are of the form \eqref{Egnform}.  Recall that $\chib_n$ is our original character, with upper ramification breaks $(u_1, \ldots, u_n)$. 
Furthermore, we saw in \S\ref{Switt} that $\chib_n$ corresponds (upon completion at $t=0$) to a (truncated) Witt vector 
$w_n := (f_1, \ldots, f_n) \in W_n(k((t)))$, and we may assume that each $f_i \in t^{1-m}k[t^{-m}]$, and all terms of $f_i$ have prime-to-$p$ degree.
 
If $u_n = pu_{n-1}$, Corollary \ref{Cgoodred} shows that there exists
$G \in \G$ giving rise to a character $\chi$ with good reduction
$\chib$ corresponding (after completion at $t=0$) to the Witt vector
$w_{\chi}: = (f_1, \ldots, f_{n-1}, f_{\chi})$, where $f_{\chi} \in
t^{1-m}k[t^{-m}]$ has degree less than $u_n$ in $t^{-1}$, and all
terms of prime-to-$p$ degree.  If $u_n > pu_{n-1}$, 
Corollary \ref{Ccrithubgoodred} and Lemma \ref{Luncoeff2} guarantee (after
a possible finite extension of $R$) the existence
of $\alpha \in R^{\times}$ and $G \in \G_{g, \alpha}$ such that
$G$ gives rise to a character $\chi$ with good reduction $\chib$
corresponding (after completion at $t=0$) to the Witt vector $w_{\chi}
:= (f_1, \ldots, f_{n-1}, f_{\chi})$, where $f_{\chi} \in t^{1-m}k[t^{-m}]$ has
degree $u_n$ in $t^{-1}$, all terms of prime-to-$p$ degree, and the
coefficient of $t^{-u_n}$ in $f_{\chi}$ is the same as that in $f_n$.
In both cases, $f_{\chi}$ and $f_n$ differ by a polynomial of degree
less than $u_n$ in $t^{-1}$.

\begin{thm}\label{Tmincase}
There exists $G_n \in \G$ giving rise to a ($\psi_n$-equivariant) character $\chi_n$ with good reduction $\chib_n$.
\end{thm}

\proof
Let $f = f_n - f_{\chi} \in t^{1-m}k[t^{-m}]$, 
which has degree less than $u_n$ in $t^{-1}$.
By Proposition \ref{Preduce} (in the case $u_n = pu_{n-1}$) or \ref{Pcrithubreduce} (in the case
$u_n > pu_{n-1}$), there exists $G_n \in \G$ such that 
$$\frac{G_n}{G} \equiv 1 + p^{p/(p-1)}F \pmod{(\KK^{\times})^p},$$ 
for some $F \in \KK$ satisfying 
$v_0(F) = 0$ and $[F]_0 = f$.  By Lemma \ref{Lartinschreier},
replacing $G$ by $G_n$ gives rise to a character $\chi_n$ whose
reduction corresponds to the Witt vector $(f_1, \ldots, f_{n-1},
f_{\chi} + f) = (f_1, \ldots, f_{n-1}, f_n)$.  In other words, the
reduction of $\chi_n$ is $\chib_n$.  Since $G_n$, by virtue of being
in $\G$, is of the form of \eqref{Egnform}, Proposition \ref{Psufficient}
shows that $\chi_n$ has good reduction and is $\psi_n$-equivariant.
\Endproof

Since the $G_n$ guaranteed by Theorem \ref{Tmincase} lies in $\G$, all the zeroes and poles of $G_n$ have valuation
$r_{\crit}$ or (in the case $u_n > pu_{n-1}$), valuation $r_{\hub}$.  By
Proposition \ref{Phublocation}, we conclude that
$\BB(\chi_n) \subseteq D(r_n)$. In particular, $\chi_n$ is admissible.  
Thus, we have proven Theorem \ref{Tsetup}, using Claim \ref{Creform}
in place of Claim \ref{claim1}.  Since Lemma \ref{Lbase} was proven in \S\ref{Sbase}, we obtain Theorem \ref{Tmachine}.

\subsection{Calculations}\label{Scalc}
Maintain the notation of the previous sections.  The purpose of
this section is twofold: to prove Lemma \ref{Luncoeff}, which is used in the
proof of Lemma \ref{Luncoeff2}, and to prove
Proposition \ref{PCconstant}, which shows that $C$ as defined in
Definition \ref{Dconstant} depends only on $g$.  Propositions
\ref{Preduce} and \ref{Pcrithubreduce} have much more complicated
proofs, and are deferred to \S\ref{Sproofs}.

\begin{lem}\label{Lleadingcoeff}
Suppose $\alphab_1, \alphab_2, C \in k^{\times}$.  For $i = 1, 2$, write
$$\omega_i  = \frac{C(-\alphab_i^m)^{N_2/m}\, dt}{(t^m - \alphab^m)^{N_2/m}t^{N_1 + u_{n-1} + 1}}.$$
Then, when expanded out
as a power series in $t^{-1}$, one obtains
$$\omega_2 - \omega_1 = (C(-\alphab_2^m)^{N_2/m} -
C(-\alphab_1^m)^{N_2/m})t^{-(u_n+1)}dt + \text{higher order terms.}$$
\end{lem}

\proof This is a straightforward computation, using the fact that $u_n
= N_1 + N_2 + u_{n-1}$.
\Endproof

\begin{cor}\label{Cleadingcoeff}
Let $\alphab_1, \alphab_2 \in k^{\times}$ with
$\alphab_1^{N_2} \neq \alphab_2^{N_2}$.  Choose lifts $\alpha_i$ of
the $\alphab_i$ in $R$.  If $G_{\crit} \in \G_{\crit,g}$ and
$G_{\hub, i} \in \G_{\hub, \alpha_i}$ ($i = 1, 2$) are chosen as in
Proposition \ref{Pcritreplace} such that $G_{\crit} G_{\hub, i}$ gives
rise to a character $\chi_i$
such that $\omega_i := \omega_{\chi_i}(r_{\hub})$ ($i = 1, 2$), then
$\delta_{\chi_i}(r_{\hub}) = \delta_{\hub}$ and the expansion
of $\omega_{\chi_2 \chi_1^{-1}}(r_{\hub})$ as a power series in $t^{-1}$
is $$\omega_2 - \omega_1 = (C(-\alphab_2^m)^{N_2/m} -
C(-\alphab_1^m)^{N_2/m})t^{-(u_n+1)}dt + \text{higher order terms.}$$
\end{cor}

\proof
This follows from Lemma \ref{Lleadingcoeff}, using Propositions
\ref{Paddchar} and \ref{Phubstart}.
\Endproof

\begin{lem}\label{Lcomparecoeffs}
Suppose $F_1 = 1 + \sum_{i=1}^{\infty} a_iT^{-i}$ and $F_2 = 1 +
\sum_{i=1}^{\infty} b_iT^{-i}$ lie in $\KK \cap (R\{T^{-1}\} \otimes K)$.
Suppose that for some $0 < \alpha \leq p/(p-1)$ and $M > 0$ not
divisible by $p$, we have 
$v(a_i) \geq \alpha$ for all $i > 0$ (with strict inequality holding
when $p | i$), and that $v(b_i) \geq \alpha$ for all $i
\geq M$ (with strict inequality holding for $i > M$). Lastly, suppose that $F_1H^p = F_2$ for some $H \in
\KK^{\times}$.  Then $v(a_i - b_i) > \alpha$ for all $i \geq M$.
\end{lem}

\proof
We may assume $H = 1 + \sum_{i=1}^{\infty} c_iT^{-i}$ and $H^p = 1 +
\sum_{i=1}^{\infty} d_iT^{-i}$ as power series expansions.  It
suffices to show that $v(d_i) > \alpha$ for all $i \geq M$.  
If $v(c_i) > \alpha/p$ for all $i$, then we are done.  
If not, let $i_0$ be the maximal $i$ such that $v(c_i) \leq \alpha/p$
(such an $i_0$ must exist).  If $i_0 < M/p$, then $v(d_i) > \alpha$
for all $i \geq M$.  If $i_0 > M/p$, then $v(d_{pi_0}) \leq \alpha$,
and $F_1H^p = F_2$ then shows that $v(b_{pi_0}) \leq \alpha$, contradicting our
assumptions on $F_2$.   
\Endproof

We now show that the value of $C$ from Definition \ref{Dconstant} only
depends on $g$.  Recall from \eqref{Es} that $s = p/(p-1) -
\delta_{\hub} = u_nr_{\hub}$.

\begin{prop}\label{PCconstant}
Assume $u_n > pu_{n-1}$.  When $\omega_{\chi_{\crit}}(r_{\hub})$ is
expanded as a Laurent series in $t$, the coefficient of $t^{-(N_1
  + u_{n-1} + 1)}dt$  
does not depend on the choice of
$G_{\crit} \in \G_{\crit, g}$, so long as
$\delta_{\chi_{\crit}}(r_{\hub}) = \delta_{\hub}$.
\end{prop}

\proof
For $i = 1, 2$, suppose $G_{\crit, i} \in \G_{\crit, g}$ gives rise to
a character $\chi_{\crit, i}$ with $\delta_{r_{\hub}}(\chi_{\crit, i})
= \delta_{\hub}$.  Proposition \ref{Paddchar} shows that
$\delta_{\chi_{\crit, 2}\chi_{\crit, 1}^{-1}}(r_{\hub}) \geq \delta_{\hub}$.
If $\delta_{\chi_{\crit, 2}\chi_{\crit, 1}^{-1}} > \delta_{\hub}$, then by
Proposition \ref{Paddchar}, we have $\omega_{\chi_{\crit,1}}(r_{\hub})
=\omega_{\chi_{\crit,2}}(r_{\hub})$, and we are done.  
If $\delta_{\chi_{\crit, 2}\chi_{\crit, 1}^{-1}}(r_{\hub}) = \delta_{\hub}$, then
Proposition \ref{Paddchar} shows that we must prove that
$\omega_{\chi_{\crit, 2}\chi_{\crit, 1}^{-1}}(r_{\hub})$ does not have
a nontrivial term of the form $ct^{-(N_1 + u_{n-1} + 1)}dt$.  

Now, $\chi_{\crit, 2}\chi_{\crit, 1}^{-1} =
\K_n((G_{\crit,2}/G_{\crit, 1})^{p^{n-1}})$, which is identified with
$\K_1(G_{\crit, 2}/G_{\crit,1})$. 
By Proposition \ref{Pdelta}, we can write 
$$G_{\crit, 2}/G_{\crit, 1} \equiv 1 + p^{s}F
\pmod{(\KK^{\times})^p}$$ 
for some $F \in \KK$ with $v_{r_{\hub}}(F) = 0$.  
If we write $1 + p^{s}F = 1 + \sum_{i=1}^{\infty} b_iT_{\hub}^{-i}$, then by multiplying by a $p$th power, we may assume
that $v(b_i) > s$ whenever $p | i$.  According to Proposition \ref{Pdelta}, we must show that either $p | (N_1 + u_{n-1})$ or
$v(b_{N_1+ u_{n-1}}) > s$.  So we assume $p \nmid (N_1 + u_{n-1})$.

On the other hand, since both $G_{\crit, i}$ lie in $\G_{\crit, g}$,
we may assume that the $G_{\crit, i}$ are chosen in $\KK$ such that 
the quotient $G_{\crit, 2}/G_{\crit, 1}$ lies in $1 + T_{\crit}^{-1}\m\{T_{\crit}^{-1}\}$.
If we write $G_{\crit, 2}/G_{\crit, 1} = 1 + \sum_{i=1}^{\infty} a_iT_{\hub}^{-i}$, 
then 
\begin{equation}\label{Efirsteq}
v(a_i) > (N_1 + u_{n-1})(r_{\crit} - r_{\hub}) = s
\end{equation}
for all $i \geq N_1 + u_{n-1}$.

The proposition now follows from Lemma \ref{Lcomparecoeffs}, taking
$M$, $\alpha$, $F_1$, and $F_2$ to be $N_1 + u_{n-1}$, $s$, $G_{\crit,
  2}/G_{\crit, 1}$, and $1 + p^sF$, respectively.
\Endproof

Let $C \in k^{\times}$ be the coefficient of
$t^{-(N_1+u_{n-1}+1)}dt$ from Proposition \ref{PCconstant} (equivalently, Definition \ref{Dconstant}). 

\begin{lem}\label{Luncoeff}
Suppose $u_n > pu_{n-1}$.
Let $\alphab_1, \alphab_2 \in k^{\times}$ with $\alphab_1^{N_2} \neq
\alphab_2^{N_2}$.  Choose lifts $\alpha_i$ of the $\alphab_i$ to $R$.
For $i = 1, 2$, suppose $G_i \in \G_{g, \alpha_i}$ 
gives rise to a character with good reduction.  Then we can write
$$G_2/G_1 \equiv 1 + p^{p/(p-1)}F \pmod{(\KK^{\times})^p},$$ where $F
  \in \KK$ satisfies $v_0(F) = 0$ and $[F]_0$ is a polynomial in
  $t^{-1} = [T^{-1}]_{0}$ of degree $u_n$ with leading term 
$$Cu_n^{-1}((-\alphab_1^m)^{N_2/m} - (-\alphab_2^m)^{N_2/m})t^{-u_n}.$$
\end{lem}

\proof
Since $G_1$ and $G_2$ both give rise to characters with \'{e}tale
reduction, Proposition \ref{Paddchar} shows that $\K_{n}((G_2/G_1)^{p^{n-1}})$
has \'{e}tale reduction as well.  In particular, $\K_1(G_2/G_1)$ has \'{e}tale reduction.
By Proposition \ref{Pdelta}, we can write $$G_2/G_1 \equiv 1 +
p^{p/(p-1)}F \pmod{(\KK^{\times})^p}$$ for some $F \in \KK$ with $v_0(F)
\geq 0$.  Furthermore, by multiplying by a $p$th power, we may assume
that when $F$ is expanded out as a power series in $T^{-1}$, no terms with
valuation $0$ have degree divisible by $p$.

On the other hand, Corollary \ref{Cleadingcoeff} and Proposition
\ref{Pdelta} show that we can write
$$G_2/G_1 \equiv 1 + p^s\Phi \pmod{(\KK^{\times})^p},$$
where $v_{r_{\hub}}(\Phi) = 0$ and $[\Phi]_{r_{\hub}}$ has derivative
$$C((-\alphab_2^m)^{N_2/m} -
(-\alphab_1^m)^{N_2/m})t^{-(u_n+1)} + \text{higher order terms})dt$$
 when expanded out as a power series in $t^{-1} =
 [T^{-1}]_{r_{\hub}}$.  

Write $1 + p^s\Phi = 1 + \sum_{i=1}^{\infty} a_iT^{-i}$.
Then $v(a_{u_n}) =  s + u_nr_{\hub} = p/(p-1)$ and $v(a_i) > p/(p-1)$ when $i > u_n$.
Also, $$[a_{u_n}T^{-u_n}]_0 = Cu_n^{-1}((-\alphab_1^m)^{N_2/m} - (-\alphab_2^m)^{N_2/m})t^{-u_n}.$$
The lemma now follows from Lemma \ref{Lcomparecoeffs}, taking $M$,
$\alpha$, $F_1$, and $F_2$ to be $u_n$, $p/(p-1)$, $1 + p^s \Phi$, and
$1 + p^{p/(p-1)}F$, respectively.
%
%
\Endproof


\subsection{Geometry of the branch locus}\label{Sbranchgeometry}

In this section, we briefly summarize the geometry of the branch locus
of the lifts our method gives for extensions as in Theorem
\ref{Tmachine}.  We only sketch the arguments.  Recall that we start with a $\Gamma$-extension
$k[[z]]/k[[s]]$ whose $\ints/p^n$-subextension has upper ramification
breaks $(u_1, \ldots, u_n)$ and no essential ramification.  We have shown that
we can lift this to a $\Gamma$-extension $R[[Z]]/R[[S]]$.  Let
$k[[t]]/k[[s]]$ and $R[[T]]/R[[S]]$ be the respective intermediate
subextensions.  The generic
fiber of $\Spec R[[Z]] \to \Spec R[[T]] \to \Spec R[[S]]$ corresponds
to a tower of branched covers of non-archimedean disks.  Since $\Spec R[[T]]
\to \Spec R[[S]]$ is simply a $\ints/m$-cover totally branched at $S =
0$, we describe the branch locus of the $\ints/p^n$-cover $\Spec
R[[Z]] \to \Spec R[[T]]$.

For each $1 \leq j \leq n$, let $i = n + 1 - j$, and let $N_{i,1}$ be
as in Theorem \ref{Tmachine}.  For $j < n$, there
are $u_i - u_{i-1}$ branch points of index $p^j$ arranged as follows:
$N_{i,1}$ of these branch points are equidistant from each other and from
the origin, at a mutual distance of $|p|^{r_{i, \crit}}$, where $r_{i,
  \crit} = {1/u_{i-1}(p-1)}$.  If we let
$N_{i,2} = u_i - u_{i-1} - N_{i,1}$, then the other
$N_{i,2}$ branch points lie at a distance of $p^{r_{i, \hub}}$ from
the origin, where $$r_{i, \hub} = \frac{1}{N_{i,2}} - \frac{N_{i,1}}{(p-1)u_{n-1}N_{i,2}}$$
(assuming $N_{i,2} > 0$).  This all follows from the discussion before
Proposition \ref{Phublocation} and the inductive nature of the proof
of Theorem \ref{Tmachine}.  The $N_{i,2}$ branch points at a distance
of $p^{r_{i, \hub}}$ from the origin come in $m$ families,
with any two points in distinct families at a distance of exactly
$|p|^{r_{i,\hub}}$ from each other.  Within a family, the distance between any two
of them is not easy to calculate exactly, but it cannot be more than $|p|^{r_{i, \hub} + s_i/(N_{i,2}/m - 1)}$, where $s_i = p/(p-1) - u_ir_{i, \hub}$.
This follows from \eqref{Es} and Definition \ref{Dghub} via a
Newton polygon argument, along with the inductive nature of the proof
of Theorem \ref{Tmachine}.

The $u_1 + 1$ branch points of index $p^n$ are arranged as in
\S\ref{Sbase}.  Specifically, if
$u_1 < p$, then all these points are equidistant from the origin
and from each other, at a mutual distance of $|p|^{1/u_1(p-1)}$.  If
$u_1 > p$, then these points all lie at a distance $|p|^{p\epsilon}$
from the origin, where $\epsilon$ can be chosen in the interval 
$(1/pu_1(p-1), 1/u_1(p-1))$.  We refer the reader to \cite[Theorem 4.3]{BW:ll} for
the finer geometry of this situation when $m = 2$, and leave the
generalization to $m > 2$ as an exercise.

Since any two branch points of distinct indices lie at distinct
distances from the origin, the ultrametric inequality determines their
distance from each other uniquely.

\section{Examples of lifting}\label{Sexamples}

In this section, we write down several examples where the isolated
differential data criterion holds, and we derive consequences for the
local lifting problem.  

\subsection{Instances of the isolated differential data criterion}\label{Sisolated}
Because of Theorem \ref{Tsetup}, the quadruples for which the isolated
differential data criterion (Definition \ref{Disolated}) is of interest to us are those of the form
$(p, m, \tilde{u}, N_1)$, where $N_1 = (p-1)\tilde{u}$ or $N_1 =
(p-1)\tilde{u} - m$.  Recall that we always assume $m | (p-1)$ and
$\tilde{u} \equiv -1 \pmod{m}$.

We begin with some small examples:
\begin{prop}\label{PD9}
The quadruples $(3, 2, 1, 2)$, $(3, 2, 1, 0)$, $(3, 2, 5, 8)$, and
$(3, 2, 5, 10)$ satisfy the isolated differential data criterion.
\end{prop}

\proof
The first two cases are covered under Proposition \ref{Pisosmall},
which does not depend on this proposition, so
suppose we are in the third or fourth cases.

Let $f_8 = t^8 + t^6 + 1$, and let $f_{10} = 2t^{10} + t^8 + t^6
+ 1$.  We claim that $f_{N_1}$ realizes the isolated differential
data criterion for $(3, 2, 5, N_1)$ with $N_1 \in \{8, 10\}$.  Let $\omega_{N_1} = dt/f_{N_1}t^6$.
The assertion that
$$\C(\omega) = \omega - \frac{dt}{t^6}$$ can be checked using a
computer algebra system, for instance SAGE (it is easier to verify
that
$$\C(f_{N_1}^3\omega) = f_{N_1}\omega - \frac{f_{N_1}dt}{t^6},$$ as no
power series are necessary--- indeed, the
computation is small enough to be checked by hand).  Thus $f_{N_1}$ realizes
the differential data criterion for $(3, 2, 5, N_1)$.

Let $\xb_1, \ldots, \xb_{N_1/2}$ be a complete set of representatives of the
$\mu_2$-equivalence classes of the roots of $f_{N_1}$.  We note that
the $\xb_j^2$ are pairwise distinct.  Then $f_{N_1}$
realizes the isolated differential data criterion if the matrix in
Remark \ref{Rconstants} is invertible over $k$.  That is, we must show
that the matrix
$$A_{N_1} = \left((\xb_j^2)^i\right)_{i, j}$$with $j \in \{1, \ldots,
  N_1/2\}$ and $i \in \{0, 2, 3, 5\}$ ($N_1 = 8$) or $i \in \{0, 2,
  3, 5, 6\}$ ($N_1 = 10$) is invertible.  Heinemann's formula for generalized
  Vandermonde determinants 
  (\cite[Theorem IV]{He:gvd} --- take $n=4$ and $s=3$ in the formula
  if $N_1 = 8$ and $n=5$ and $s=4$ if $N_1 = 10$) shows that 
$$\det(A_{N_1}) = \begin{cases} D \det\newmatrix{e_3}{e_4}{e_0}{e_1} & N_1
  = 8 \\ D \det\newmatrix{e_4}{e_5}{e_1}{e_2} & N_1 = 10 \end{cases}$$
where
\begin{itemize}
\item $D$ is an integral power of the (standard) Vandermonde
  determinant corresponding to $\xb_1^2, \ldots, \xb_{N_1/2}^2$.
\item For all $s$, the number $e_s$ is the $s$th
  elementary symmetric polynomial in the $\xb_j^2$'s.
\end{itemize}
Since the standard Vandermonde determinants are invertible,
we need only show that $e_3e_1 - e_4e_0$ is invertible when $N_1 = 8$,
and that $e_4e_2 - e_5e_1$ is invertible when $N_1 = 10$.  If $N_1 =
8$, then
$\prod_{j=1}^4(t - \xb_j^2) = t^4 + t^3 + 1$, and thus
$e_3e_1 - e_4e_0  = -1$.  If $N_1 = 10$, then
$\prod_{j=1}^5(t - \xb_j^2) = t^5 + 2t^4 + 2t^3 + 2$, and $e_4e_2 -
e_5e_1 = -1$.  We are done.
\Endproof

\begin{prop}\label{Pisosmall}
For any (odd) prime $p$, the isolated differential data criterion
holds for $(p, 2, 1, N_1)$, when $N_1 = p-1$ or $N_1 = p-3$. 
\end{prop}

\proof
Using Proposition \ref{Pcrit}, we may realize the differential data
criterion by solving the system of equations \eqref{Eass71}. If $N_1 =
0$ there is nothing to do, so assume otherwise.  Let $\xb_j =
j$ for all $j \in \{1, 2, \ldots, N_1/2\}$.  The square matrix
$\left(\xb_j^q\right)_{q, j}$ for $j \in \{1, 2, \ldots, N_1/2\}$ and $q
\in \{1, 3, \ldots, N_1 - 1\}$ is Vandermonde (up to multiplication of
each column by a nonzero scalar).  Since $N_1 < p$, the $\xb_j$ all have
distinct squares and thus the columns of this Vandermonde matrix are
all distinct. 
So the system \eqref{Eass71}
has a unique solution for the $a_j$ with the $a_j \in k$.  Since the
$\xb_j^q$ in fact lie in $\FF_p$, so do the $a_j$.  We must show
that no $a_j$ is zero.  

For a contradiction, assume, after possibly renumbering the $\xb_j$, that $a_{N_1/2}
= 0$.  If $N_1 = 2$, this is clearly a contradiction, and we are
done.  Assume otherwise.  Since $u = 1$ in
\eqref{Eass71},
we must have
$$\sum_{j=1}^{N_1/2 - 1} a_j \xb_j^q = 0$$ for all $q \in \{3, 5,
\ldots, N_1 - 1\}$.  Since $\xb_j \neq 0$, we once again have (up
to rescaling) a Vandermonde system of linear equations for the $a_j$,
$j \in \{1, 2,
  \ldots, N_1/2 - 1\}$.  Thus all $a_j$ are zero, which is a
  contradiction.  This gives the differential data criterion.

To prove isolatedness, we remark that the matrix 
$\left(\xb_j^{q-1}\right)_{q,j}$ for $q \in \{1, 3, \ldots,
N_1-1\}$ and $j \in \{1, 2, \ldots, N_1/2\}$ in Remark \ref{Rconstants} is,
up to scaling, Vandermonde with distinct columns.  So it is invertible.
\Endproof   

\begin{lem}\label{Lmainiso}
The quadruple $(p, m, \tilde{u}, (p-1)\tilde{u})$ satisfies the
differential data criterion for all odd primes $p$, all $m
| (p-1)$, and all $\tilde{u} \equiv -1 \pmod{m}$.
\end{lem}

\proof
As in Proposition \ref{Pisosmall}, we will realize the differential data
criterion by solving the system of equations \eqref{Eass71}.  
Write $\tilde{u} = up^{\nu}$, with $p \nmid u$.
Note that the set $S$ of $u(p^{\nu + 1} - 1)$th roots of unity
whose $-u$th powers have trace zero 
(under $\Tr_{\FF_{p^{\nu+1}}/\FF_p}$) 
has cardinality $u(p^{\nu} - 1)$.  Thus, we have $$|\mu_{u(p^{\nu+1} -1)} \backslash S| = u(p^{\nu+1} - p^{\nu}).$$  
Furthermore, multiplication by $m$th roots of unity (which all lie in $\FF_p$) 
preserves $S$ and $\mu_{u(p^{\nu+1} -1)} \backslash S$.  We take the $\xb_i$ to be any complete
set of orbit representatives for the multiplicative action of $\mu_m$ on $\mu_{u(p^{\nu+1} -1)} \backslash S$.  Note that there are 
$$u(p^{\nu+1} - p^{\nu})/m = (p-1)\tilde{u}/m$$ of these orbits, so we have the correct number of $\xb_j$.  Furthermore, for each $\xb_j$, 
let the associated $a_j$ be given by the formula
$$a_j = -\Tr(\xb_j^{-u}) = -\sum_{i = 0}^{\nu} \xb_j^{-up^i},$$ where
for simplicity, we write $\Tr$ for $\Tr_{\FF_{p^{\nu+1}}/\FF_p}$.  
This is a nonzero element of $\FF_p$.

We then have, for any $q \equiv -1 \pmod{m}$:
\begin{eqnarray*}
\sum_{j=1}^{N/m} a_j\xb_j^{q} &=& \sum_{j=1}^{N/m} -\Tr(\xb_j^{-u})\xb_j^{q} \\
&=& \frac{1}{m} \sum_{x \in \mu_{u(p^{\nu+1} -1)} \backslash S} -\Tr(x^{-u})x^{q} \\
&=& \frac{1}{m} \sum_{x \in \mu_{u(p^{\nu+1} -1)}} -\Tr(x^{-u})x^{q} \\
&=& \frac{1}{m} \sum_{x \in \mu_{u(p^{\nu+1}-1)}} -\left(x^{q-u} + x^{q - up} + \cdots + x^{q - up^{\nu}}\right)\\
&=& \begin{cases} u/m & q \equiv u, up, \ldots, up^{\nu} \pmod{u(p^{\nu+1} - 1)} \\ 0 & \text{otherwise.} \end{cases}
\end{eqnarray*}
The second equality above comes from the fact that $u \equiv -1 \pmod{m}$ and 
$q \equiv -1 \pmod{m}$, so multiplying any $\xb_j$ by any $m$th root
of unity leaves $\Tr(\xb_j^{-u})\xb_j^{q}$ invariant. 
This solves system (\ref{Eass71}) when we restrict to the case $1 \leq
q \leq (p-1)\tilde{u} + \tilde{u}  - 1 = up^{\nu + 1} -1$ and $p \nmid q$.  
\Endproof

\begin{prop}\label{Pmainiso}
In the situation of Lemma \ref{Lmainiso}, if $\tilde{u} = (m-1)p^{\nu}$
for some $\nu \geq 0$, then $(p, m, \tilde{u}, (p-1)\tilde{u})$ satisfies the
isolated differential data criterion.
\end{prop}

\proof
Let the $a_j$ and $\xb_j$ be as in Lemma \ref{Lmainiso}.  Recall that each 
$\xb_j$ is in $\mu_{u(p^{\nu+1} -1)}$.

For any set $\Sigma \subseteq \ZZ$, let $\ol{\Sigma}$ be its image as a subset of $\ZZ/u(p^{\nu+1} - 1)$.
Write $u = m-1$.  Let 
$$A = \left(\xb_j^{q-1}\right)_{q,j}$$ with $j$ ranging from $1$ to $(p-1)up^{\nu}/m$
and $q$ ranging from $1$ to $up^{\nu+1}-1$ over those numbers congruent
to $-1 \pmod{m}$ and not divisible by $p$.  By Remark
\ref{Rconstants}, it suffices to show that $A$ is invertible.

We first claim that the set
$$B := \{q \ | \ q \text{ corresponds to a row of }A\}$$ and the set
$$C := \{m-1 + imp\}_{0 \leq i < \frac{N}{m}}$$ satisfy 
$\ol{B} = \ol{C}$.   
To prove the claim, note that if
$$C' = \{m-1 + imp\}_{0 \leq i < \frac{u(p^{\nu+1}-1)}{m}},$$ 
then $\ol{C'}$ is exactly the set of elements of $\ZZ/u(p^{\nu+1}-1)$
congruent to $-1 \pmod{m}$ (this abuse of language is justified since
$m \mid (p-1) \mid u(p^{\nu+1}-1)$).  Furthermore, a straightforward computation
shows that $$\ol{C'} \backslash \ol{C} = \{up, (u+m)p, \ldots, up^{\nu+1} - mp\}.$$
Now, by the Chinese Remainder Theorem,
the set $$\{up-1, (u+m)p -1, \ldots, up^{\nu+1} - mp - 1\},$$ viewed as a subset of $\ZZ$, is exactly the set of integers between
$1$ and $up^{\nu+1} - 1$ which are congruent to $-1 \pmod{m}$ and to $0 \pmod{p}$ (this is where we use $u = m-1$).
Thus $\ol{C'} \backslash \ol{C} = \ol{C'} \backslash \ol{B}$.  Since $\ol{B} \subseteq \ol{C}'$ is clear, we have $\ol{B} = \ol{C}$, proving the claim.

The claim shows that the elements of $\ol{B}$, arranged appropriately, form an arithmetic progression with common difference $mp$.  
If $A'$ is the matrix obtained by rearranging the rows of $A$ to correspond to this ordering, then the definition of $A$ shows that the $j$th column of $A'$
is a geometric progression with common ratio $\xb_j^{mp}$.  The common ratios of the columns are pairwise distinct, 
as the $x_j$ are all $u(p^{\nu+1}-1)$th roots of unity lying in pairwise distinct multiplicative $\mu_m$-orbits, and $p \nmid u(p^{\nu+1}-1)$.
We can scale each column to make a new matrix $A''$ where the first entry in each column is equal to $1$.  Then $A''$ is a Vandermonde matrix with 
pairwise distinct column ratios.  So $A''$ is invertible, which means $A$ is invertible.
\Endproof

\begin{rem}\label{Rnotiso}
It is not hard to show, in the context of Lemma \ref{Lmainiso}, that
if $\tilde{u}$ is \emph{not} a $p$th power times $(m-1)$, then the
proposed solution in Lemma \ref{Lmainiso} will \emph{never} realize the
isolated differential data criterion.  Indeed, the matrix $A$ from the
proof of Proposition \ref{Pmainiso} can be shown to have at least two identical
rows.
\end{rem}
 
\subsection{Affirmative local lifting results}\label{Sllresults}

\begin{thm}\label{TD9}
The dihedral group $D_9$ is a local Oort group for $p = 3$.
\end{thm}

\proof
By Proposition \ref{Ppopreduction}, we need only consider
$D_9$-extensions whose $\ints/9$-subextension has upper jumps $(1,
3)$, $(1, 5)$, $(1, 7)$, $(5, 15)$, $(5, 17)$, or $(5, 19)$.  By
Theorem \ref{Tmachine}, it suffices to show that the isolated
differential data criterion holds for $(3, 2, 1, 2)$, $(3, 2, 1, 0)$,
$(3, 2, 5, 10)$, and $(3, 2, 5, 8)$.  This follows from Proposition \ref{PD9}.
\Endproof

\begin{thm}\label{Tfirstjump1}
If $p$ is an odd prime, and $L/k[[s]]$ is a $D_{p^2}$-extension whose $\ints/p^2$-subextension has first upper ramification break
$u_1 \equiv 1 \pmod{p}$, then $L/k[[s]]$ lifts to characteristic zero.
\end{thm}

\proof
Since $u_1$ is odd, we have that $u_1 \equiv 1 \pmod{2p}$.
By Proposition \ref{Ppopreduction}, we need only consider
$D_{p^2}$-extensions whose $\ints/p^2$-subextension has first upper
jump $1$.
By Theorem \ref{Tmachine}, it suffices to show that the isolated
differential data criterion holds for $(p, 2, 1, p-1)$ and $(p, 2, 1,
p-3)$
This follows from Proposition \ref{Pisosmall}.
\Endproof

\begin{thm}\label{Tminimaljumps}
If $L/k[[s]]$ is a $\ints/p^n \rtimes \ints/m$-extension whose
$\ints/p^n$-subextension has upper ramification breaks
congruent to $(m-1, p(m-1), \ldots, p^{n-1}(m-1))$ $\pmod{mp}$, then $L/k[[s]]$ lifts to characteristic zero. In particular, $\ints/p^n \rtimes \ints/m$ is
a weak local Oort group whenever the conjugation action of $\ints/m$ on $\ints/p^n$ is faithful.
\end{thm}

\proof
By Proposition \ref{Ppopreduction}, we need only consider
$\ints/p^n \rtimes \ints/m$-extensions whose
$\ints/p^n$-subextension has upper ramification breaks
$(m-1, p(m-1), \ldots, p^{n-1}(m-1))$ (such extensions exist by
\cite[Theorem 1.1]{OP:wt}).  By Theorem \ref{Tmachine}, it suffices to show that the isolated
differential data criterion holds for $$(p, m, (m-1)p^{\nu-1},
(p-1)(m-1)p^{\nu-1})$$ for $0 \leq \nu < n$.  This follows from Proposition \ref{Pmainiso}.
\Endproof

\section{Proof of Propositions \ref{Preduce} and \ref{Pcrithubreduce}}\label{Sproofs}
We use the notation of \S\ref{Sgeomsetup}, \S\ref{Scharacters} and \S\ref{Sproof} throughout.
In particular, recall that 
\begin{itemize}
\item $pu_{n-1} \leq u_n < pu_{n-1} + mp$ (no essential ramification).
\item $N = N_1 + N_2 = u_n - u_{n-1}$, and both $N_1$ and $N_2$ are
  divisible by $m$ (Proposition \ref{Pdistrib}).
\item $N_1 \leq (p-1)u_{n-1}$ with strict inequality unless $u_n =
  pu_{n-1}$ (Proposition \ref{Pdistrib}). 
\item $N_2 \leq mp$ (Assumption \ref{AN1N2}).
\item $r_{\crit} = 1/(p-1)u_{n-1}$ (beginning of \S\ref{Scrit}).
\item $r_{\hub} = 1/N_2 - N_1/(p-1)u_{n-1}N_2$, or $r_{\hub} = 0$ if
  $N_2 = 0$ (Proposition \ref{Phublocation}).
\item $s = (N_1 + u_{n-1})(r_{\crit} - r_{\hub}) = p/(p-1) -
  u_nr_{\hub} = p/(p-1) - \delta_{\hub}$ (Equation \eqref{Es}).
\item If $r \in \rats_{\geq 0}$, then $T_r = p^{-r}T$.  For short,
  $T_{\crit} = p^{-r_{\crit}}T$ and $T_{\hub} = p^{-r_{\hub}}T$
  (beginning of \S\ref{Scharacters}).
\item $v_r$ and $v_r'$ are defined as in Definition
  \ref{Dpartialvaluation} (and $v_r$ is a valuation on $R\{T_r^{-1}\} \otimes_R K$).
\item $\G_{\crit, g}$ and $\G_{\hub, \alpha}$ are defined as in Definitions
  \ref{Dgcrit}, \ref{Dghub}, respectively.  Here $g$ is a solution to
  \eqref{Ecrit}, corresponding to an $f$ realizing the isolated
  differential data criterion (Remark \ref{Rnotjustexistence}).
\end{itemize}

As a matter of notation, in the context of a congruence between two
power series or polynomials in $T_r^{-1}$, the symbol $\equiv'$
(resp.\ $='$)
means that the congruence (resp.\ equality) need only hold for terms of degree congruent to $-1 \pmod{m}$ or $0 \pmod{p}$ in $T_r^{-1}$.

While neither Proposition \ref{Preduce} nor Proposition
\ref{Pcrithubreduce} follows directly from the other, their proofs are
very similar, and we will prove them simultaneously.  Essentially, the
proof of Proposition \ref{Preduce} is an easier version of the proof
of Proposition \ref{Pcrithubreduce}.

\subsection{Preliminaries}\label{Stechnical}
We start by defining the ch (think ``crit-hub'') ``valuation,'' which is not actually a
valuation, but has many similar properties.
\begin{defn}\label{Dcrithub}
For a power series
$$F = \sum_{q=0}^{N_1 + u_{n-1} - 1} c_q T_{\crit}^{-q} + \sum_{q = N_1 + u_{n-1}}^{\infty} p^s c_q T_{\hub}^{-q}
\in R\{T_{\hub}^{-1}\} \otimes_R K,$$ we write 
$v_{\crithub}(F) = \min_q v(c_q)$.
We define $v_{\crithub}'(F)$ in the same way, except we only take the minimum over $q$ that are congruent either to $-1 \pmod{m}$ or $0 \pmod{p}$. 
\end{defn}

\begin{defn}\label{Dnegligible}
\begin{enumerate}
\item
An element $f \in R\{T_{\hub}^{-1}\} \otimes_R K$ is called
\emph{hub-negligible} if, for all $r \in [0, r_{\hub}] \cap \rats$, we
have $v_r(f) > p/(p-1) - u_nr$.  If $f \in R\{T_{\hub}^{-1}\}
\otimes_R K$, then making a \emph{hub-negligible adjustment} to $f$ means
replacing it with some $f' \in R\{T_{\hub}^{-1}\} \otimes_R K$ where
$f' - f$ is hub-negligible.
\item
An element $f \in R\{T_{\crit}^{-1}\} \otimes_R K$ is called
\emph{crit-negligible} if, for all $r \in [r_{\hub}, r_{\crit}] \cap \rats$, we
have $v_r(f) > (N_1 + u_{n-1})(r_{\crit} - r)$.  If $f \in R\{T_{\hub}^{-1}\}
\otimes_R K$, then making a \emph{crit-negligible adjustment} to $f$ means
replacing it with some $f' \in R\{T_{\crit}^{-1}\} \otimes_R K$ where
$f' - f$ is crit-negligible.
\end{enumerate}
\end{defn} 

\begin{lem}\label{Lnegligibletest}
\begin{enumerate}
\item Let $f = \sum_{i=0}^{\infty} c_iT_{\hub}^{-i} \in R\{T_{\hub}^{-1}\}
\otimes_R K$.  If $v_0(c_iT_{\hub}^{-i}) > p/(p-1)$ for $i < u_n$ and
$v_{\crithub}(f) > 0$, then $f$ is hub-negligible.
\item Let $f = \sum_{i=0}^{\infty} c_iT_{\crit}^{-i} \in R\{T_{\crit}^{-1}\}
\otimes_R K$.  If $v_0(c_iT_{\crit}^{-i}) > p/(p-1)$ for $i < N_1 + u_{n-1}$ and
$v_{r_{\crit}}(f) > 0$, then $f$ is crit-negligible.
\end{enumerate}
\end{lem}

\proof
It suffices to check each monomial in $f$. In case (i), if $i < u_n$, the
definition of hub-negligibility yields
$v_r(c_iT_{\hub}^{-i}) > p/(p-1) - ir$, proving the lemma for these terms.  
When $i \geq u_n$, the fact that
$v_{\crithub}(c_iT_{\hub}^{-i}) > 0$ implies that
$v_{r_{\hub}}(c_iT_{\hub}^{-i}) = v(c_i) > s$.  Now, $s = p/(p-1) -
u_nr_{\hub}$.
Thus, for $r \leq r_{\hub}$, $$v_r(c_iT_{\hub}^{-i}) > s + i(r_{\hub} - r) \geq s +
u_n(r_{\hub} - r) = p/(p-1) - u_n r.$$  

In case (ii), if $i < N_1 + u_{n-1}$, then 
$$v_r(c_iT_{\crit}^{-i}) > p/(p-1) - ir > p/(p-1) - (N_1 + u_{n-1})r
\geq (N_1 + u_{n-1})(r_{\crit} - r),$$
because $N_1 + u_{n-1} \leq pu_{n-1}$.  This proves the lemma for
these terms.  If $i \geq N_1 + u_{n-1}$, then 
$$v_r(c_iT_{\crit}^{-i}) > i(r_{\crit} - r) \geq (N_1 +
u_{n-1})(r_{\crit} - r),$$ and we are done.
\Endproof 

\begin{rem}\label{Rnegligible}
Lemma \ref{Lnegligibletest} shows that if $f \in R\{T_{\hub}^{-1}\}
\otimes_R K$ and $v_{\crithub}(f) > 0$, then removing the terms of $f$
of degree at least $u_n$ in $T^{-1}$ is a hub-negligible adjustment.
Likewise, if $f \in R\{T_{\crit}^{-1}\} \otimes_R K$ and
$v_{r_{\crit}}(f) > 0$, then removing the terms of $f$ of degree at least
$N_1 + u_{n-1}$ in $T^{-1}$ is a crit-negligible adjustment.
\end{rem}

\begin{lem}\label{Lcrithub1}
Suppose $f_1$ and $f_2$ are in $R\{T_{\hub}^{-1}\} \otimes_R K$.
\begin{enumerate}
\item
We have $v_{\crithub} (f_1+f_2) \geq \min (v_{\crithub}(f_1),
v_{\crithub}(f_2))$, with equality if $v_{\crithub}(f_1) \neq
v_{\crithub}(f_2)$, and the same holds
for $v_{\crithub}'$.   
\item
We have $v_{\crithub}(f_1f_2) \geq v_{\crithub}(f_1) + v_{\crithub}(f_2)$.
\end{enumerate}
\end{lem}

\proof
Part (i) is obvious, and reduces part (ii) to the case of monomials.  The only non-obvious case is if $f_1 = aT_{\crit}^{-b}$ and $f_2 = cp^sT_{\hub}^{-d}$,
where $b < N_1 + u_{n-1}$ and $b+d \geq N_1 + u_{n-1}$.  Then $v_{\crithub}(f_1) = v(a)$ 
and $v_{\crithub}(f_1f_2) = v(a) + v(c) + b(r_{\crit} - r_{\hub})$.
If $d \geq N_1 + u_{n-1}$, then $v_{\crithub}(f_2) = v(c)$, which proves part (ii) since $r_{\crit} > r_{\hub}$.  If $d < N_1 + u_{n-1}$, then 
$v_{\crithub}(f_2) = v(c) + s - d(r_{\crit} - r_{\hub})$.  So $v_{\crithub}(f_1f_2) - v_{\crithub}(f_1) - v_{\crithub}(f_2) = (b+d)(r_{\crit} - r_{\hub}) - s$.  
This is nonnegative, since $b+d \geq N_1 + u_{n-1}$.  This proves part (ii).
\Endproof

\begin{cor}\label{Ccrithub}
If $f_1$ and $f_2$ are in $R\{T_{\hub}^{-1}\} \otimes_R K$ with $v_{\crithub}(f_1)$ and $v_{\crithub}(f_2) \geq 0$, then
$v_{\crithub}(f_1f_2-1) \geq \min(v_{\crithub}(f_1-1), v_{\crithub}(f_2-1))$.
\end{cor}

\proof
Since $f_1f_2 - 1 = (f_1-1)(f_2-1) + (f_1-1) + (f_2-1)$, the corollary follows from Lemma \ref{Lcrithub1}.
\Endproof

\begin{rem}\label{Rcrithub}
Of course, since $v_{r_{\crit}}$ is a valuation, Lemma \ref{Lcrithub1} and Corollary \ref{Ccrithub} are also
true when applied to $R\{T_{\crit}^{-1}\} \otimes_R K$, with
$v_{r_{\crit}}$ and $v_{r_{\crit}}'$ replacing $v_{\crithub}$ and
$v_{\crithub}'$, respectively.
\end{rem}

\begin{lem}\label{Lcrithub2}
If $f \in R\{T_{\hub}^{-1}\} \otimes_R K$, then $v_{\crithub}(f) \geq \max(v_{r_{\crit}}(f), v_{r_{\hub}}(f) - s)$, and the same is true when 
$v_{\crithub}$, $v_{r_{\crit}}$, and $v_{r_{\hub}}$ are replaced by
$v_{\crithub}'$, $v_{r_{\crit}}'$, and $v_{r_{\hub}}'$, respectively.
Furthermore,
$v_{r_{\hub}}(f) \geq v_{\crithub}(f)$ and $v_{r_{\hub}}'(f) \geq v_{\crithub}'(f)$.
\end{lem}

\proof
It suffices to prove the statements for $v_{\crithub}$ applied to monomials $f =
T_{r_{\hub}}^{-i}$.  Then $v_{r_{\hub}}(f) = 0$ and 
$v_{r_{\crit}}(f) = i(r_{\hub} - r_{\crit})$.  If $i \geq N_1 + u_{n-1}$, then $v_{\crithub}(f) = -s$, which is greater than
$i(r_{\hub} - r_{\crit})$ and nonpositive.  If $i < N_1 +
u_{n-1}$, then $v_{\crithub}(f) = i(r_{\hub} - r_{\crit})$, which is
greater than or equal to $-s$
and nonpositive.
\Endproof

\begin{lem}\label{Lnegligibletest2}
\begin{enumerate}
\item Let $f,g \in R\{T_{\hub}^{-1}\}
\otimes_R K$ such that $f$ is hub-negligible and $v_{\crithub}(g)\geq 0$.
Then $fg$ is hub-negligible.
\item Let $f,g \in R\{T_{\crit}^{-1}\}
\otimes_R K$ such that $f$ is crit-negligible and $v_{r_{\crit}}(g)\geq 0$.
Then $fg$ is crit-negligible.
\end{enumerate}
\end{lem}

\proof
For part (i), since $v_{\crithub}(g) \geq 0$, we have $v_r(g) \geq 0$ for all $r
\leq r_{\hub}$.  Thus $$v_r(fg) = v_r(f) + v_r(g) \geq v_r(f) >
p/(p-1) - u_n r.$$  Part (ii) is similar and just as easy.
\Endproof

\begin{lem}\label{Linvert}
If $f  = 1 + h \in 1 + T_{\hub}^{-1}R\{T_{\hub}^{-1}\} \otimes_R K$ with
$v_{\crithub}(h) = \beta > 0$, then $v_{\crithub}(f^{-1} - 1) =
\beta$.
\end{lem}

\proof
We have $f^{-1} - 1 = 1/(1 + h) - 1 = -h + h^2 - h^3 + \cdots$.  Now
the result follows from Lemma \ref{Lcrithub1}.
\Endproof 

\begin{lem}\label{LItoH}
Assume that $N_1 = (p-1)u_{n-1} - m$, so that $r_{\hub} > 0$.
Let $$I := 1 + \sum_{l=1}^{\lfloor(N_1 + u_{n-1} - 1) /p\rfloor} b_l T_{\crit}^{-pl} +  \sum_{l = \lceil(N_1 + u_{n-1})/p\rceil}^{\lfloor(u_n - 1)/p\rfloor} p^s b_l T_{\hub}^{-pl},$$
with all $b_l \in \m$. 
Let $$H := 1 + \sum_{l=1}^{\lfloor(N_1 + u_{n-1} - 1) /p\rfloor} b_l^{1/p} T_{\crit}^{-l} +  \sum_{l = \lceil(N_1 + u_{n-1})/p\rceil}^{\lfloor(u_n - 1)/p\rfloor} (p^s b_l)^{1/p} T_{\hub}^{-l},$$
for any choice of $p$th roots of the coefficients.  
Then each term $c_iT^{-i}$ in $I - H^p$ for $i \geq 1$ satisfies
$$v_{\crithub}(c_iT^{-i}) > \theta_i + \frac{v_{\crithub}(I)}{p},$$ where
$$\theta_i = 
\begin{cases} 
\frac{p-1}{p}\left(\frac{p}{p-1} - i
    r_{\crit}\right) & i < N_1 + u_{n-1} \\
\frac{p-1}{p}\left(\frac{p}{p-1} - ir_{\hub} - s\right) & i \geq N_1 +
u_{n-1}
\end{cases}$$

The same holds when, instead of taking $I - H^p$,  we expand out $I/H^p - 1$ as a power series in
$T^{-1}$.   
\end{lem}

\proof
The terms in $H^p - I$ are the cross-terms in $H^p$.  We consider the
two cases separately.  Note that the multinomial coefficient in any
cross-term of $H^p$ has valuation at least $1$. 

Suppose $i < N_1 + u_{n-1}$.  Then we must show that
$v_{r_{\crit}}(c_iT^{-i}) \geq 1 + v_{\crithub}(I)/p - i/pu_{n-1}$.  
Each term in $H$ can be written either as $b_l^{1/p} T_{\crit}^{-l}$
or as $b_l^{1/p}p^{s/p}p^{-l(r_{\crit} - r_{\hub})}T_{\crit}^{-l}$.
Note that $v(b_l^{1/p}) \geq v_{\crithub}(I)/p$.
If no terms of the second form factor into the given $c_iT^{-i}$, then
the result is obvious.  If at least one such term factors in, then
$$v_{r_{\crit}}(c_iT^{-i}) \geq 1 + \frac{v_{\crithub}(I)}{p} +
\frac{s}{p} - i(r_{\crit} - r_{\hub}).$$  Since $s = (N_1 +
u_{n-1})(r_{\crit} - r_{\hub})$, it suffices to show that
$$((N_1 + u_{n-1})/p - i)(r_{\crit} - r_{\hub}) > - i/pu_{n-1}.$$
Substituting in $r_{\crit} = 1/(p-1)u_{n-1}$ and $r_{\hub} =
m/N_2(p-1)u_{n-1}$, and multiplying both sides by $-pu_{n-1}$, 
we are reduced to showing that
$$i > (pi - N_1 - u_{n-1})\left(\frac{N_2 - m}{N_2(p-1)}\right).$$
Since $N_2 \leq mp$, the right hand side is at most
$i - (N_1 + u_{n-1})/p$ (if it is positive), from which the result follows.

Now suppose $i \geq N_1 + u_{n-1}$.  The we must show that
$$v_{r_{\hub}}(c_iT^{-i}) > 1 + v_{\crithub}(I)/p -
\frac{p-1}{p}(ir_{\hub} + s) + s.$$  At least one term of the
form $b_l^{1/p}p^{s/p}T_{\hub}^{-l}$ factors into $c_iT^{-i}$ and all
terms factoring in have nonnegative valuation at $r_{\hub}$.  So
$v_{r_{\hub}}(c_iT^{-i}) \geq 1 + v_{\crithub}(I)/p + s/p$.   The desired
inequality follows immediately.

To prove the statement for $I/H^p - 1$, note that 
$I/H^p - 1 = (I-H^p)H^{-p}$.  By Lemma \ref{Linvert}, we have
$v_{\crithub}(H^{-p}) = 0$.  Write $(I/H^p - 1) = \sum_{i=1}^{\infty} d_i T^{-i}$.
By Lemma \ref{Lcrithub1}, 
$$v_{\crithub}(d_iT^{-i}) \geq \min_{j \leq i} v_{\crithub}(c_jT^{-j})
> \min_{j \leq i} \theta_j + \frac{v_{\crithub}(I)}{p} = \theta_i + \frac{v_{\crithub}(I)}{p}.$$

%
\Endproof

\begin{rem}\label{RItoH}
If we take $\theta_i$ to be any number less than $1$, then Lemma \ref{LItoH} also holds for $I \in 1 + T_{\crit}^{-p}\m[T_{\crit}^{-p}]$ and
$v_{r_{\crit}}$ replacing $v_{\crithub}$.  In particular, if we assume $I$
has degree less than $N_1 + u_{n-1}$ in $T^{-1}$, then we may define 
$\theta_i$ (for $i < N_1 + u_{n-1}$) as in Lemma \ref{LItoH}.
\end{rem}

\subsection{The underlying Hensel's lemma calculation}
For each of Propositions \ref{Preduce} and \ref{Pcrithubreduce}, we
get most of the way to a proof via an application of Hensel's lemma.
For Proposition \ref{Preduce}, the necessary result is as follows.

\begin{lem}[cf.\ \cite{OW:ce}, Lemma 7.4(i)]\label{Limprovelem} 
Let $G\in\G_{\crit, g}$, and let $J \in 1 + T_{\crit}^{-1}\m\{T_{\crit}^{-1}\}$. 
There exists a unique $G' \in \G_{\crit, g}$ and a unique 
polynomial $I \in 1 +  T_{\crit}^{-p}\m[T_{\crit}^{-p}]$ of degree $< N_1 + u_{n-1}$ in $T_{\crit}^{-1}$
such that $$\frac{G'}{G}I \equiv' J \pmod{T_{\crit}^{-(N_1 + u_{n-1})}}.$$ 

  If $J \equiv' 1 \pmod{p^{\beta}, T_{\crit}^{-(N_1 + u_{n-1})}}$ for 
$\beta \in \QQ_{> 0}$, then $v_{r_{\crit}}(G'/G-1) \geq \beta$ and $v_{r_{\crit}}(I-1) \geq \beta$.
\end{lem}     

\proof
By assumption we have
\[
 G = \prod_{j=1}^{N_1/m} \prod_{\ell=1}^m (1 - \zeta_m^{-\ell}x_jT_{\crit}^{-1})^{\psi_1(\tau^{\ell})a_j},
\]
where $x_j\in R$ is a lift of $\xb_j$, where the $\xb_j$ are a solution to \eqref{Eass71} corresponding to $g$. We set
\[
   G' = \prod_{j=1}^{N_1/m} \prod_{\ell=1}^m (1 - \zeta_m^{-\ell}x_j'T_{\crit}^{-1})^{\psi_1(\tau^{\ell})a_j},
 \quad x_j':=x_j+\epsilon_j,
\]
and where the $\epsilon_j$ are for the moment considered as indeterminates. We also set
\[
      I := 1 + \sum_{l=1}^{\lfloor(N_1 + u_{n-1} - 1) /p\rfloor} b_l
      T_{\crit}^{-pl}
\]
for another system of indeterminates $b_l$. Write
\[
      \frac{G'}{G}I = 1 +\sum_{q=1}^\infty c_q T_{\crit}^{-q},
\]
where $c_q$ is a formal power series in $(\epsilon_j,b_l)$.  One computes,
using $\psi_1(\tau^{\ell}) = \zeta_m^{-\ell}$ (Lemma
\ref{Lwhichroot}), that
\begin{equation} \label{Ederivs} 
\begin{split}
  \frac{\partial c_q}{\partial \epsilon_j}|_{\epsilon_j=b_l=0} =
        & \;\sum_{\ell=1}^m \zeta_m^{-(q+1)\ell} a_jx_j^{q-1}  \\
  \frac{\partial c_q}{\partial b_l}|_{\epsilon_j=b_l=0} = &\;
    \begin{cases} \;\;1, & q=pl \\ \;\;0, & q\neq pl. \end{cases}
\end{split}
\end{equation} 
In particular, when $q \equiv -1 \pmod{m}$, we have
$$\frac{\partial c_q}{\partial \epsilon_j}|_{\epsilon_j=b_l=0} = ma_jx_j^{q-1},$$
and otherwise $\partial c_q/\partial \epsilon_j = 0$.
The congruence 
\begin{equation} \label{Emaincong}
  \frac{G'}{G}I \equiv' J \pmod{T_{\crit}^{-(N_1 + u_{n-1})}}
\end{equation}
corresponds to a system of equations in the indeterminates $(\epsilon_j,b_l)$, one equation for each 
$c_{q}T_{\crit}^{-q}$ for $q \equiv -1 \pmod{m}$ or $q \equiv 0 \pmod{p}$, with $1 \leq q < N_1 + u_{n-1}$.
The Jacobian matrix $M_{\crit}$ of this system of
equations is invertible over $R$ if and only if its reduction is invertible over $k$.  From
\eqref{Ederivs} it is easy to see that this is true iff the matrix from \eqref{Einvertiblematrix} is
invertible (One obtains the matrix in \eqref{Einvertiblematrix} from the Jacobian matrix by eliminating 
all of the columns corresponding to the $b_l$, which are standard
basis vectors, along with the rows corresponding to equations for which
$p|q$).  The matrix from \eqref{Einvertiblematrix} is invertible
because we are assuming that $g$ realizes 
the isolated differential data criterion for $(p, m, u_{n-1}, N_1)$.
By Hensel's lemma, we conclude that
\eqref{Emaincong} has a (unique) solution with $\epsilon_j,b_l\in\m$,
proving the first statement of the lemma.
In fact, by the effective Hensel's Lemma, the second statement holds
as well.
\Endproof

The analogous result toward Proposition \ref{Pcrithubreduce} is the following:

\begin{lem}\label{Lcrithubimprovelem} 
Suppose $u_n > pu_{n-1}$.  Let $G_{\crit} \in \G_{\crit,g}$, and $G_{\hub}\in \G_{\hub, \alpha}$, and let 
$$J = 1+ \sum_{q=1}^{N_1 + u_{n-1} - 1} c_q T_{\crit}^{-q} + \sum_{q = N_1 + u_{n-1}}^{\infty} p^s c_q T_{\hub}^{-q}$$ with
all $c_q \in \m$ and $\lim_{q \to \infty} c_q = 0$.  
There exists a unique $G_{\crit}' \in \G_{\crit, g}$, a unique $G_{\hub}' \in \G_{\hub, \alpha}$ and a unique 
polynomial 
$$I := 1 + \sum_{l=1}^{\lfloor (N_1 + u_{n-1} - 1) /p\rfloor} b_l
T_{\crit}^{-pl} +  \sum_{l = \lceil (N_1 + u_{n-1})/p\rceil}^{\lfloor (u_n - 1)/p\rfloor} p^s b_l T_{\hub}^{-pl},$$
such that 
\begin{equation}\label{Emaincong2}
\frac{G_{\crit}'G_{\hub}'}{G_{\crit}G_{\hub}}I \equiv' J \pmod{T_{\hub}^{-u_n}}. 
\end{equation}

If for some $\beta \in \QQ_{> 0}$, we have $v(c_q) \geq \beta$ for all $q < u_n$,
then $v_{r_{\crit}}(G_{\crit}'/G_{\crit} - 1) \geq {\beta}$, $v_{\crithub}(I-1) \geq {\beta}$, and $v_{r_{\hub}}(G_{\hub}'/G_{\hub} - 1) \geq s + \beta$.
\end{lem}     

\proof
As in Lemma \ref{Limprovelem}, we have
$$G_{\crit} = \prod_{j=1}^{N_1/m} \prod_{\ell=1}^m (1 - \zeta_m^{-\ell}x_jT_{\crit}^{-1})^{\zeta_m^{-\ell}a_j}.$$
Furthermore, by Definition \ref{Dghub}, we have
$$G_{\hub}  = \prod_{\ell=0}^{m-1} \left(1 + p^s\frac{\sum_{j=1}^{N_2/m-1} y_j \zeta_m^{-j\ell}T_{\hub}^{-j}}{(1 -
  \alpha\zeta_m^{-\ell}T_{\hub}^{-1})^{N_2/m -1}}\right)^{\zeta_m^{-\ell}}$$ 
with $y_j \in R$ (divide the numerator and
denominator in \eqref{Ehubform} by $T_{\hub}^{N_2/m - 1}$).  Here, we are thinking of $\zeta_m$ as an integer
given by taking some arbitrary lift of $\zeta_m \in \FF_p^{\times}$ to
$\ints$.

We look for potential solutions for $G_{\crit}'$ and $G_{\hub}'$ in the forms 
\begin{equation}\label{Egcritdeform}
   G_{\crit}' = \prod_{j=1}^{N_1/m} \prod_{\ell=1}^m (1 - \zeta_m^{-\ell}x_j'T_{\crit}^{-1})^{\zeta_m^{-\ell}a_j},
 \quad x_j':=x_j+\epsilon_j
\end{equation}
and
\begin{equation}\label{Eghubdeform}
G_{\hub}' = \prod_{\ell=0}^{m-1} \left(1 + p^s 
\frac{\sum_{j=1}^{N_2/m-1} y_j' \zeta_m^{-j\ell}T_{\hub}^{-j}}{(1 - \alpha\zeta_m^{-\ell}T_{\hub}^{-1})^{N_2/m - 1}}\right)^{\zeta_m^{-\ell}}, 
\quad y_j' := y_j + \gamma_j,
\end{equation}
where the $\epsilon_j$ and $\gamma_j$ are considered as
indeterminates.  Write
\begin{equation}\label{Ejacobian}
   \frac{G_{\hub}'G_{\crit}'}{G_{\hub}G_{\crit}}I = 1 +\sum_{q=1}^{N_1 + u_{n-1} - 1} c_q T_{\crit}^{-q} + 
   \sum_{q = N_1 + u_{n-1}}^{\infty} p^s c_q T_{\hub}^{-q},
\end{equation}
where $c_q$ is a formal power series in $(\epsilon_j,\gamma_j, b_l)$.
By \eqref{Ejacobian}, the congruence \eqref{Emaincong2}
expresses the $c_q$ relevant to $\equiv'$-congruence for $q < u_n$ in terms
of formal power series in the indeterminates $(\epsilon_j, \gamma_j,
b_l)$.  We take $M$ to be the Jacobian of this system of equations at $0$.
More specifically, let $M$ be the Jacobian matrix (at $\epsilon_j =
\gamma_j = b_l = 0$ for all $j, l$) of the following outputs and inputs:
For the outputs, we take the variables $c_q$ for $q < u_n$, 
where either $q \equiv -1 \pmod{m}$ or $p | q$.  For the input
variables, we take the $\epsilon_j$, the $b_l$ for $pl < N_1 + u_{n-1}$, the
$\gamma_j$, and the $b_l$ for $pl \geq N_1 + u_{n-1}$, in that order.  
The matrix $M$ will be shown to be invertible over $R$ in Proposition
\ref{PMinvert}.  We conclude by Hensel's lemma that
\eqref{Emaincong2} has a (unique) solution with $\epsilon_j, \gamma_j,
b_l\in\m$.  In fact, by the effective Hensel's Lemma, $(\epsilon_j)$, $v(\gamma_j)$, and $v(b_l)$
are all at least as large as $\min_{q < u_n} v(c_q)$.  Given the forms in \eqref{Egcritdeform} and \eqref{Eghubdeform}, this proves the lemma.
\Endproof
 
\begin{rem}
The reason we rescale some of the $c_q$ and the $b_l$ by $p^s$ is to
force $M$ to be invertible.  Our scaling of the $c_q$ motivates the
definition of $v_{\crithub}$ in \S\ref{Stechnical}.
\end{rem}

The rest of this section is dedicated to proving that the matrix $M$
in the proof of Lemma \ref{Lcrithubimprovelem} has entries in $R$ and is invertible over $R$.

%

Let us calculate the entries of $M$, using the notation of the proof of
Lemma \ref{Lcrithubimprovelem}.
To do this, we calculate the partial derivatives of the $c_q$ with respect to the $\epsilon_j$, $\gamma_j$, and $b_l$ at the point
$\epsilon_j = \gamma_j = b_l = 0$ (all partials calculated below are evaluated at this point, and we suppress the point in the notation).  
For $q \equiv -1 \pmod{m}$, as in \eqref{Ederivs}, we have
\begin{equation} \label{Ederivs2} 
  \frac{\partial c_q}{\partial \epsilon_j} = 
  	\begin{cases} \;\;ma_jx_j^{q-1} & q < N_1 + u_{n-1} \\
	\;\;ma_jx_j^{q-1}p^{(r_{\crit} - r_{\hub})q - s} & q \geq N_1 + u_{n-1}
        \end{cases}
\end{equation}
Also, we have
\begin{equation} \label{Ederivs3} 
  \frac{\partial c_q}{\partial b_l} = \;
    \begin{cases} \;\;1, & q=pl \\ \;\;0, & q\neq pl. \end{cases}
\end{equation}
To calculate $\partial c_q/\partial \gamma_j$, first set $G_{\hub, 0}$ and $G_{\hub, 0}'$ equal to the $\ell = 0$ factors of $G_{\hub}$ and $G_{\hub}'$, 
respectively.  Then 
$$\frac{G_{\hub,0}'}{G_{\hub,0}} = 1 + p^s\frac{\sum_{j=1}^{N_2/m-1} \gamma_j T_{\hub}^{-j}}{(1 - \alpha T_{\hub}^{-1})^{N_2/m-1} + p^sC(T_{\hub})}.$$
When this is expanded out as a power series in $T_{\hub}^{-1}$, the coefficient of $T_{\hub}^{-q}$ is 
$$p^s \left(\sum_{j=1}^{N_2/m-1} \gamma_j \alpha^{q-j}\binom{q-j + N_2/m - 2}{N_2/m - 2} + O(p^s)\right),$$
where $O(p^s)$ represents terms with valuation at least $s$.  A computation now yields that
\begin{equation}\label{Ederivs4}
  \frac{\partial c_q}{\partial \gamma_j} = 
  	\begin{cases} \;\; p^{s - (r_{\crit} - r_{\hub})q} \sum_{\ell = 0}^{m-1} \left(\zeta_m^{-(q+1)\ell} \alpha^{q-j}\binom{q-j + N_2/m - 2}{N_2/m - 2} + O(p^s) \right) & q < N_1 + u_{n-1} \\
	\;\; \sum_{\ell = 0}^{m-1} \left(\zeta_m^{-(q+1)\ell} \alpha^{q-j}\binom{q-j + N_2/m - 2}{N_2/m - 2} + O(p^s) \right) & q \geq N_1 + u_{n-1}.
        \end{cases}
\end{equation}
In particular, when $q \equiv -1 \pmod{m}$, we have
\begin{equation}\label{Ederivs5}
  \frac{\partial c_q}{\partial \gamma_j} = 
  	\begin{cases} \;\; mp^{s - (r_{\crit} - r_{\hub})q} \left( \alpha^{q-j}\binom{q-j + N_2/m - 2}{N_2/m - 2} + O(p^s) \right) & q < N_1 + u_{n-1} \\
	\;\;  m \alpha^{q-j}\binom{q-j + N_2/m - 2}{N_2/m - 2} + O(p^s) & q \geq N_1 + u_{n-1}.
        \end{cases}
\end{equation}
It is clear from the above formulas that the entries of $M$ lie in $R$.
Write 
\begin{equation}\label{Emblock}
M = \blockmatrix{M_1}{M_2}{M_3}{M_4}
\end{equation}
as a block matrix, with the columns of 
$M_1$ corresponding to the variables $x_j$ and $b_l$ for $pl < N_1 + u_{n-1}$, 
and the rows of $M_1$ corresponding to the $c_q$ for $q < N_1 +
u_{n-1}$.  Then one checks that $M_1$ 
is a square matrix of size $N_1/m + \lfloor (N_1 + u_{n-1} - 1)/p
\rfloor$ (cf.\ Remark \ref{Rsquare}), and
$M_4$ is square as well as will be seen in the proof of Proposition
\ref{PMinvert} below.  In particular, $M$ is a square matrix.  

\begin{prop}\label{PMinvert}
The matrix $M$ is invertible over $R$.  
\end{prop}

\proof
It suffices to show that the reduction $\bar{M} =
\blockmatrix{\bar{M}_1}{\bar{M}_2}{\bar{M}_3}{\bar{M}_4}$ of $M$ has nonzero determinant.
From \eqref{Ederivs4}, the valuation of $\partial c_q/\partial{\gamma_j}$ for $q < N_1 + u_{n-1}$ is at least
$s - (r_{\crit} - r_{\hub})q$, which is $(r_{\crit} - r_{\hub})(N_1 + u_{n-1} - q) > 0$.  Also, $\partial c_q/\partial b_l$ for $pl \geq N_1 + u_{n-1}$
and $q < N_1 + u_{n-1}$ is $0$ by \eqref{Ederivs3}.  Thus $\bar{M}_2 = 0$.  So $\bar{M}$ is block lower triangular, and
$\det(\bar{M}) = \det(\bar{M}_1)\det(\bar{M}_4)$.  But $M_1$ is just the Jacobian
matrix for the system in \eqref{Emaincong}, where it was shown that $M_1$ is invertible over $R$. 
Thus $\det(\bar{M}_1) \neq 0$.  So we are reduced to showing that $\bar{M}_4$ is invertible.

Each column of $\bar{M}_4$ corresponding to a variable $b_l$ has
a $1$ in the row corresponding to $q = pl$ and a $0$ in each other
position.  Eliminating these columns and the rows where $1$'s appear,
we are left with an 
$(N_2/m -1) \times (N_2/m - 1)$ matrix $\bar{M}_4'$. The entries of
$\bar{M}_4'$ are the reductions of $\partial c_q/ \partial \gamma_j$, where $1 \leq j \leq N_2 -1$, and
$q$ ranges from $N_1 + u_{n-1}$ to $u_n - 1$ over those numbers congruent to $-1 \pmod{m}$ and not
divisible by $p$.  By \eqref{Ederivs5}, 
after multiplying rows and columns by units,
the entry of $\bar{M}_4'$ corresponding to $(q, j)$ is $\binom{q-j + N_2/m - 2}{N_2/m - 2}$, thought of as an element of $\FF_p \subseteq k$.
We will view the binomial coefficients as integers, and show that the determinant is not divisible by $p$.

We will modify $\bar{M}_4'$, without changing its determinant.
For the first modification, moving from left to right, we subtract the
$j=2$ column from the $j=1$ column.  
Then we subtract the $j=3$ column from the $j=2$ column.  We continue until 
we subtract the $j = N_2/m-1$ column from the $j = N_2/m - 2$ column.  This gives a matrix whose entry in the
$(q, j)$ slot is $\binom{q - j + N_2/m - 3}{N_2/m - 3}$, except in the last column, where the entries are 
$\binom{q - j + N_2/m - 2}{N_2/m - 2}$.  For the second modification, we 
repeat this process once more, except that we stop after subtracting
the $j = N_2/m - 2$ column from the $j = N_2/m - 3$ column.
For the third modification, we repeat again, stopping after
subtracting the $j = N_2/m - 3$ column from the $j = N_2/m -4$ column. 
We continue repeating until the $(N_2/m -2)$nd modification, which
consists only of subtracting the $j=2$ column from the $j=1$ column. 
All in all, the $j$th column gets modified $N_2/m - j - 1$ times. 
This leaves us with
a matrix whose entry in the $(q, j)$-slot is $\binom{q-1}{j-1}$.  We apply the formula given on \cite[p.\ 308]{GV:bd}
(the ``alternate expression" when $b = 0$) to get that the determinant of this matrix is
\begin{equation}\label{Egessel}
\frac{\prod_{1 \leq i < j \leq N_2/m - 1}(b_i - b_j)}{1!2!\cdots(N_2/m - 1)!},
\end{equation} 
where the $b_i$ are the values of $q$ corresponding to our $c_q$.

It suffices to check that the numerator in \eqref{Egessel} is not divisible by $p$ (in any case, the denominator is not divisible by $p$ because 
$N_2 \leq mp$ by Assumption \ref{AN1N2}).   The expression $b_i - b_j$ can only take on values $m, 2m, \ldots, N_2 - m$, as 
$(u_n - m) - (N_1 + u_{n-1}) = N_2 - m$, and $u_n - m$ and $N_1 + u_{n-1}$ are the least and greatest values of $q$, respectively.  By
Assumption \ref{AN1N2}, we have $N_2 \leq mp$, so the expression $b_i - b_j$ is never divisible by $p$.  We are done.
\Endproof

\subsection{Completion of the proofs}\label{Scompletion}

The main task in completing the proofs is to turn the $I$ that occurs
in Lemmas \ref{Limprovelem} and \ref{Lcrithubimprovelem}, and that is
very close to a $p$th power,  into an
actual $p$th power.  This will be done through a series of results.
In each case, we will state and prove the result relevant to Lemma
\ref{Lcrithubimprovelem}.  Then we will state the analogous result
relevant to Lemma \ref{Limprovelem}, and mention which modifications
are necessary for the proof to carry through.  As a matter of fact, there are
more straightforward proofs of most of the ``Lemma \ref{Limprovelem}
versions,'' but since we must write the more complicated versions
anyway, we omit the simpler versions to save space.

\begin{lem}\label{Lcrithubimprovelem2}
Suppose $u_n > pu_{n-1}$.  Let $G_{\crit}, G_{\hub} \in
\G_{\crit, g}, \G_{\hub, \alpha}$ respectively.
Let $J \in 1 + T^{-1}\m\{T^{-1}\}$ such that $v_{\crithub}(J-1) > 0$.  
Let $\theta_i$ be as in Lemma \ref{LItoH}.  After a possible finite
extension of $K$ and hub-negligible adjustment to $J$,
there exist $G_{\crit}', G_{\hub}' \in \G_{\crit, g}, \G_{\hub, \alpha}$ respectively, and a polynomial $H \in 1 +
T^{-1}\m[T^{-1}]$ such that if 
$$\frac{J}{(G_{\crit}'G_{\hub}'/G_{\crit}G_{\hub})H^p} = 1 +
\sum_{i=1}^{\infty} c_iT^{-i},$$ then for $0 < i < u_n$, there exists
$\epsilon > 0$ such that
$$v_{\crithub}(c_iT^{-i}) \geq 
\begin{cases} \min\left(v_{\crithub}(J-1) +
  v_{\crithub}'(J-1), \theta_i + \epsilon + \frac{v_{\crithub}'(J-1)}{p}\right) & p | i \text{ or } i
\equiv -1 \pmod{m} \\ 
\min\left(v_{\crithub}(J-1), \theta_i + \epsilon +
  \frac{v_{\crithub}'(J-1)}{p}\right)& \text{otherwise.}
\end{cases}
$$

If $v_{\crithub}'(J - 1) \geq \beta$ for some $0 < \beta < p/(p-1)$, then we
can choose $G_{\crit}'$, $G_{\hub}'$, and $H$ above such that 
$v_{r_{\crit}}(G_{\crit}'/G_{\crit} - 1) \geq \beta$, that
$v_{r_{\hub}}(G_{\hub}'/G_{\hub} - 1) \geq s + \beta$, and that $v_{\crithub}(H^p - 1)
\geq \min(\beta, (p-1)r_{\hub}/p)$.
\end{lem}

\proof
Let $G_{\crit}'$, $G_{\hub}'$, and 
$$I := 1 + \sum_{l=1}^{\lfloor (N_1 + u_{n-1} - 1) /p\rfloor} b_l T_{\crit}^{-pl} +  \sum_{l = \lceil (N_1 + u_{n-1})/p\rceil}^{\lfloor (u_n - 1)/p\rfloor} p^s b_l T_{\hub}^{-pl},$$
be the unique solution guaranteed by Lemma \ref{Lcrithubimprovelem}.  
So $(G_{\crit}'G_{\hub}/G_{\crit}G_{\hub})I \equiv' J \pmod{T_{\hub}^{-u_n}}$.
Set
\[
H := 1 + \sum_{l=1}^{\lfloor(N_1 + u_{n-1} - 1) /p\rfloor} b_l^{1/p} T_{\crit}^{-l} +  \sum_{l = \lceil(N_1 +  u_{n-1})/p\rceil}^{\lfloor(u_n - 1)/p\rfloor} (p^s b_l)^{1/p} T_{\hub}^{-l},
\]
for any choice of $p$th roots.  
Let $L = (G_{\crit}'G_{\hub}'/G_{\crit}G_{\hub})I$.
Then $$\frac{J}{(G_{\crit}'G_{\hub}'/G_{\crit}G_{\hub})H^p} = \left(\frac{J}{L}\right)\left(\frac{I}{H^p}\right).$$  
Now, Lemma \ref{Lcrithubimprovelem} gives us that $v_{\crithub}(I - 1) \geq
v_{\crithub}'(J-1)$. Lemma \ref{LItoH} shows that, if 
$$\frac{I}{H^p} - 1 = \sum_{i=1}^{u_n-1} d_i T^{-i},$$ then for all $0 < i < u_n$
\begin{equation}\label{Ecrithub0}
v_{\crithub}(d_iT^{-i}) > \theta_i + v_{\crithub}'(J-1)/p.
\end{equation}  
Additionally, Lemma \ref{Lcrithubimprovelem}, combined with Corollary
\ref{Ccrithub} and Lemma \ref{Lcrithub2}, gives us that
$v_{\crithub}(L - 1) \geq v_{\crithub}'(J-1) > 0$.  Now, 
$J/L - 1 = (J - L)/L = (J-L)(1 + (L-1))^{-1}$ and $v_{\crithub}(J - L) = v_{\crithub}(J-1 - (L - 1)) \geq v_{\crithub}(J-1)$ by
Lemma \ref{Lcrithub1}(i). 
Using Lemmas \ref{Linvert} and \ref{Lcrithub1}(ii), we have
\begin{equation}\label{Ecrithub1}
v_{\crithub}\left(\frac{J}{L} - 1\right) \geq v_{\crithub}(J-1).
\end{equation}
On the other hand, by construction, $J - L$ has no $T^{-i}$ term if
$i \equiv -1 \pmod{m}$ or $p|i$ and $i < u_n$.  Expanding $L = 1 + (L
-1)$ out as a power series, and again using Lemma \ref{Lcrithub1}, this implies   
\begin{equation}\label{Ecrithub2}
v_{\crithub}'\left(\frac{J}{L} -1\right) \geq v_{\crithub}(J - L) + v_{\crithub}(L-1) \geq v_{\crithub}(J-1) + v_{\crithub}'(J-1).
\end{equation}
Now, we note that $\theta_i$ is a decreasing function of $i$.
Letting $c_i$ be as in the lemma, and using \eqref{Ecrithub0}, this
has the consequence that, for $i <u_n$,
\begin{equation}\label{Ecrithub3}
v_{\crithub}(c_iT^{-i}) \geq 
\begin{cases} \min\left(v_{\crithub}'\left(\frac{J}{L} - 1\right), \theta_i + \epsilon + \frac{v_{\crithub}'(J-1)}{p}\right) & p|i \text{ or
  } i \equiv -1 \pmod{m} \\
\min\left(v_{\crithub}\left(\frac{J}{L} - 1\right), \theta_i + \epsilon + \frac{v_{\crithub}'(J-1)}{p}\right) & \text{otherwise.}
\end{cases}
\end{equation}
for some $\epsilon > 0$.

Combining \eqref{Ecrithub1},
\eqref{Ecrithub2}, and \eqref{Ecrithub3} proves the first part of the lemma.  The last
statement about $G_{\crit}'$ and $G_{\hub}'$ follows easily from Lemma
\ref{Lcrithubimprovelem}.
Lemma \ref{Lcrithubimprovelem} also shows that $v_{\crithub}(I - 1)
\geq \beta$.  Since $\theta_i$ is nonincreasing in $i$, and  
$H^p$ and $I$ have degree less than $u_n$ in
$T^{-1}$, Lemma \ref{LItoH} (along with Lemma \ref{Lcrithub1}(i)) 
shows that $v_{\crithub}(H^p - 1) > \min(\beta, \theta_{u_n - 1})$.
One calculates that $\theta_{u_n - 1} = (p-1)r_{\hub}/p$, and this
completes the proof.
\Endproof

\begin{lem}\label{Lpreimprovelem2}
Let $G_{\crit} \in \G_{\crit, g}$.
Let $J \in 1 + T_{\crit}^{-1}\m\{T_{\crit}^{-1}\}$.  
Let $\theta_i$ be as in Remark \ref{RItoH}.  After a possible finite
extension of $K$ and crit-negligible adjustment to $J$,
there exist $G_{\crit}' \in \G_{\crit, g}$ and a polynomial $H \in 1 +
T^{-1}\m[T^{-1}]$ such that if 
$$\frac{J}{(G_{\crit}'/G_{\crit})H^p} = 1 +
\sum_{i=1}^{\infty} c_iT^{-i},$$ then for $0 < i < N_1 + u_{n+1}$, there exists
$\epsilon > 0$ such that
$$v_{r_{\crit}}(c_iT^{-i}) \geq 
\begin{cases} \min\left(v_{r_{\crit}}(J-1) +
  v_{r_{\crit}}'(J-1), \theta_i + \epsilon + \frac{v_{r_{\crit}}'(J-1)}{p}\right) & p | i \text{ or } i
\equiv -1 \pmod{m} \\ 
\min\left(v_{r_{\crit}}(J-1), \theta_i + \epsilon +
  \frac{v_{r_{\crit}}'(J-1)}{p}\right)& \text{otherwise.}
\end{cases}
$$

If $v_{r_{\crit}}'(J - 1) \geq \beta$ for some $0 < \beta < p/(p-1)$, then we
can choose $G_{\crit}'$ and $H$ above such that 
$v_{r_{\crit}}(G_{\crit}'/G_{\crit} - 1) \geq \beta$, and that $v_{r_{\crit}}(H^p - 1)
\geq \min(\beta, (p-1)r_{\crit}/p)$.
\end{lem}

\proof
The proof is the same as that of Lemma \ref{Lcrithubimprovelem2},
replacing Lemma \ref{Lcrithubimprovelem} by Lemma \ref{Limprovelem},
Lemma \ref{LItoH} by Remark \ref{RItoH}, Lemma \ref{Lcrithub1} and
Corollary \ref{Ccrithub} by Remark \ref{Rcrithub}, $v_{\crithub}$ and
$v_{\crithub}'$ by $v_{r_{\crit}}$ and $v_{r_{\crit}}'$, ``hub-negligible'' by ``crit-negligible,'' omitting the second summations in $I$ and $H$, replacing all
$u_n$'s by $N_1 + u_{n-1}$'s, and omitting all mentions of $G_{\hub}$ and $G_{\hub}'$.
\Endproof

\begin{lem}\label{Lcrithubimprovelem3}
Let $0 < \sigma < p/(p-1)$.  Let $G_{\crit},
G_{\hub} \in G_{\crit, g}, \G_{\hub, \alpha}$, respectively. 
Let $J \in 1 + T^{-1}\m\{T^{-1}\}$ such that $v_{\crithub}(J -1) > 0$.
After a possible finite extension of $K$
and hub-negligible adjustment to $J$,
there exist $G_{\crit}', G_{\hub}' \in \G_{\crit, g}, \G_{\hub, \alpha}$
respectively, and a polynomial $H \in 1 +
T^{-1}\m[T^{-1}]$,
such that 
$$v_{\crithub}'\left(\frac{J}{(G_{\crit}'G_{\hub}'/G_{\crit}G_{\hub})H^p} - 1 \right) \geq \sigma.$$

We can choose $G_{\crit}'$, $G_{\hub}'$, and $H$ above such that 
$$v_{\crithub}((G_{\crit}'G_{\hub}'/G_{\crit}G_{\hub})H^p - 1) >
\min(v_{\crithub}'(J-1), (p-1)r_{\hub}/p).$$
\end{lem}

\proof
We will build
$G_{\crit}'$, $G_{\hub}'$, and $H$ through successive approximation.
Let $\theta_i$ be as in Lemma \ref{LItoH}, and let $\eta_i = (p/(p-1))\theta_i$.
We make the following observations.  First, $\theta_{u_n} = \eta_{u_n} =
0$.  Second, the $\eta_i$ form a decreasing sequence.
Third, if $v_{\crithub}(c_iT^{-i}) > \eta_i$ for some $c_i
\in K$, then $v_0(c_iT^{-i}) = v(c_i) > p/(p-1)$.  Fourth, 
if $x > \eta_i$, then $\theta_i + x/p > \eta_i$.  

By the first observation above, we know that $v_{\crithub}'(J-1) > \eta_j$ for
some $0 <  j \leq u_n$.
Let $G_{\crit, 1}$, $G_{\hub, 1}$, and $H_1$ be the $G_{\crit}'$, $G_{\hub}'$, and $H$ guaranteed by Lemma
\ref{Lcrithubimprovelem2} (after making a hub-negligible adjustment to $J$), and set $J_1 := (G_{\crit, 1}'G_{\hub,
  1}'/G_{\crit}G_{\hub})H_1^p$.  It follows from Lemma
\ref{Lcrithubimprovelem2}, Corollary \ref{Ccrithub}, and Lemma
\ref{Lcrithub2} that $v_{\crithub}(J_1 - 1) > 0$.  Thus
$v_{\crithub}(J_1) = 0$.  Also, $v_{\crithub}(J_1^{-1}) = 0$ as a
consequnce of Lemma \ref{Linvert}, so $v_{\crithub}(J/J_1 - 1) =
v_{\crithub}((J-J_1)(J_1^{-1}))>0$ by Lemma \ref{Lcrithub1}(ii).
Write $J/J_1 = 1 + \sum_{i=1}^{\infty} d_iT^{-i}$.  For $i \geq j$
and either $p | i$ or $i \equiv -1 \pmod{m}$, Lemma
\ref{Lcrithubimprovelem2} and the second, third, and fourth observations above
show that 
$v_0(d_iT^{-i}) > p/(p-1)$.  For $i < j$, there exists $\epsilon > 0$
such that 
$$v_{\crithub}(d_iT^{-i}) \geq 
\begin{cases} \min\left(v_{\crithub}(J-1) +
  v_{\crithub}'(J-1), \theta_i + \epsilon + \frac{v_{\crithub}'(J-1)}{p}\right) & p | i \text{ or } i
\equiv -1 \pmod{m} \\ 
\min\left(v_{\crithub}(J-1), \theta_i + \epsilon + 
  \frac{v_{\crithub}'(J-1)}{p}\right)& \text{otherwise.}
\end{cases}
$$

If $A_1 := \sum_{i \in I} d_iT^{-i}$ where $I \subseteq [j, \infty)$
consists of those integers congruent to $0 \pmod{p}$ or $-1 \pmod{m}$,
or greater than or equal to $u_n$, 
then $A_1$ is hub-negligible by Lemma \ref{Lnegligibletest}(i) and
Remark \ref{Rnegligible}.    Since
$v_{\crithub}(J_1) = 0$, Lemma \ref{Lnegligibletest2}(i) shows that
$J_1A_1$ is hub-negligible.  So we may (and do) replace $J$ with $J -
J_1A_1$, and we assume that $d_i = 0$ for $i \in I$.

Since the $\theta_i$ form a decreasing
sequence, we have
$$v_{\crithub}'\left(\frac{J}{J_1}-1\right) \geq \min\left(v_{\crithub}(J-1) + v_{\crithub}'(J-1), \theta_{j-1} +
    \epsilon + \frac{v_{\crithub}'(J-1)}{p}\right)$$
and
$$v_{\crithub}\left(\frac{J}{J_1}-1\right) \geq \min\left(v_{\crithub}(J-1), \theta_{u_n-1} +
    \epsilon + \frac{v_{\crithub}'(J-1)}{p}\right).$$

For $l > 1$, define $G_{\crit, l}'$, $G_{\hub, l}'$, and $H_l$
inductively as the $G_{\crit}'$, $G_{\hub}'$, and $H$ guaranteed by Lemma
\ref{Lcrithubimprovelem2} with $J/J_{l-1}$ in place of $J$ and
$G_{\crit, l-1}$ and $G_{\hub, l-1}$ in place of $G_{\crit}$ and
$G_{\hub}$ (note that, since $v_{\crithub}(J_{l-1}) = 0$ for the same
reason that $v_{\crithub}(J_1) = 0$, Lemma \ref{Lnegligibletest}(i)
shows that the hub-negligible adjustment to $J/J_{l-1}$
required for Lemma \ref{Lcrithubimprovelem2} can be achieved by making
a hub-negligible adjustment to $J$). Define
$$J_l = \frac{G_{\crit,l}'G_{\hub, l}'}{G_{\crit}G_{\hub}}(H_1 \cdots
H_l)^p$$ 
so that
$$\frac{J}{J_l} = \frac{J/J_{l-1}}{(G_{\crit, l}'G_{\hub, l}'/G_{\crit,l-1}'G_{\hub, l-1}')H_i^p}.$$ 
At each stage, we replace
$J$ with $J-J_lA_l$, where $A_l$ is the part of $J/J_l$ consisting of
terms of degree $i$ in $T^{-1}$, where $j \leq i \leq u_n - 1$ and
either $p|i$ or $i \equiv -1 \pmod{m}$.  As before, this is a
hub-negligible adjustment.  By Lemma \ref{Lcrithubimprovelem2}, there exists $\epsilon > 0$ such that
$$v_{\crithub}'\left(\frac{J}{J_l} -1\right) \geq \min\left(v_{\crithub}\left(\frac{J}{J_{l-1}}-1\right) +
  v_{\crithub}'\left(\frac{J}{J_{l-1}}-1\right), \theta_{j-1} + \epsilon +\frac{v_{\crithub}'(J/J_{l-1}-1)}{p}\right)$$ and
$$v_{\crithub}\left(\frac{J}{J_l} - 1\right) \geq \min\left(v_{\crithub}\left(\frac{J}{J_{l-1}}-1\right),
  \theta_{u_n-1} + \epsilon + \frac{v_{\crithub}'(J/J_{l-1}-1)}{p}\right).$$
Since $\eta_{j-1} = (p/(p-1))\theta_{j-1}$, there exists some $l_{j-1}$ for which
$v_{\crithub}'(J/J_{l_{j-1}} - 1) > \eta_{j-1}$.

Replacing $j$ by $j-1$, we can repeat the entire process again.
Induction now shows that, after further hub-negligible adjustments to
$J$, we get down to the case $j=1$.  That is, there
exists $l_1$ such that $v_{\crithub}'(J/J_{l_1} -1) > \eta_1$.  Replacing $J$ with
$J-J_{l_1}A_{l_1}$ as above, we obtain that $v_{\crithub}'(J/J_{l_1} - 1) = \infty$.  
In particular, setting $G_{\crit}'$ and $G_{\hub}'$ equal to 
$G_{\crit, l_1}'$ and $G_{\hub, l_1}'$ respectively, and setting $H = (H_1 \cdots
H_{l_1})$, gives the desired solution.
 
To prove the last statement, note that Lemma \ref{Lcrithubimprovelem2}
shows that all $H_i$ satisfy $v_{\crithub}(H_i - 1) > \min(v_{\crithub}'(J-1), (p-1)r_{\hub}/p)$.  By
Corollary \ref{Ccrithub}, $v_{\crithub}(H - 1)$ has the same property.  Lemma
\ref{Lcrithubimprovelem2} and Lemma \ref{Lcrithub2} imply that
$v_{\crithub}(G_{\crit}'/G_{\crit} - 1)$ and
$v_{\crithub}(G_{\hub}'/G_{\hub})$ also have this property.  Combining all this with
Corollary \ref{Ccrithub} proves the last statement of the lemma.
\Endproof

\begin{lem}\label{Limprovelem2}
Let $0 < \sigma < p/(p-1)$.  Let $G_{\crit} \in G_{\crit, g}$. 
Let $J \in 1 + T_{\crit}^{-1}\m\{T_{\crit}^{-1}\}$.
After a possible finite extension of $K$
and crit-negligible adjustment to $J$,
there exist $G_{\crit}' \in \G_{\crit, g}$
and a polynomial $H \in 1 + T_{\crit}^{-1}\m[T_{\crit}^{-1}]$,
such that $$v_{r_{\crit}}'\left(\frac{J}{(G_{\crit}'/G_{\crit})H^p} - 1 \right) \geq \sigma.$$

We can choose $G_{\crit}'$ and $H$ above such that 
$v_{r_{\crit}}((G_{\crit}'/G_{\crit})H^p - 1) >
\min(v_{r_{\crit}}'(J-1), (p-1)r_{\crit}/p)$. 
\end{lem}

\proof
The proof is the same as that of Lemma \ref{Lcrithubimprovelem3},
replacing Lemma \ref{Lcrithubimprovelem2} by Lemma \ref{Lpreimprovelem2},
Lemma \ref{LItoH} by Remark \ref{RItoH}, 
Corollary \ref{Ccrithub} by Remark \ref{Rcrithub},
Lemma \ref{Lnegligibletest}(i) by Lemma \ref{Lnegligibletest}(ii), 
Lemma \ref{Lnegligibletest2}(i) by Lemma \ref{Lnegligibletest2}(ii), 
$v_{\crithub}$ and $v_{\crithub}'$ by $v_{r_{\crit}}$ and $v_{r_{\crit}}'$, ``hub-negligible'' by ``crit-negligible,'' 
$u_n$ by $N_1 + u_{n-1}$, and omitting all mentions of $G_{\hub}$ and $G_{\hub}'$.
\Endproof

\begin{lem}\label{Lcrithublessdiscrepancy}
Let $0 < \sigma < p/(p-1)$.  Let $J \in 1 +
T^{-1}\m\{T^{-1}\}$ with $v_{\crithub}(J-1) > 0$.  After a possible
finite extension of $K$ and hub-negligible adjustment to $J$,
there exists $J' \in 1 + T^{-1}\m\{T^{-1}\} \cap \KK$ such 
$J' =' J$, and $J'$ has $r_{\hub}$-discrepancy valuation at least $\sigma$. Furthermore,
$v_{\crithub}(J'-1) > 0$.
\end{lem}

\proof 
Suppose the $r_{\hub}$-discrepancy valuation of $J$ is at
least $\sigma_0 \geq 0$.  By Lemma \ref{Lcrithubimprovelem3}, after
making a hub-negligible adjustment to $J$, there exist
$G_{\crit}$, $G_{\crit}'$, $G_{\hub}$, $G_{\hub}'$, and $H$ (as in
that lemma) such
that 

\begin{equation}\label{E5}
v_{\crithub}'\left(\frac{J}{(G_{\crit}'G_{\hub}'/G_{\crit}G_{\hub})H^p}  -
  1 \right) \geq \sigma.
\end{equation}
For shorthand, write $G' = G_{\crit}'G_{\hub}'$ and $G =
G_{\crit}G_{\hub}$.  Also, by Lemma \ref{Lcrithubimprovelem3}, 
\begin{equation}\label{E6}
v_{\crithub}\left(\frac{G'}{G}H^p - 1\right) \geq \beta :=
\min(v_{\crithub}'(J-1), (p-1)r_{\hub}/p).
\end{equation} 
Equations \eqref{E5} and \eqref{E6} also hold with $v_{r_{\hub}}$ and $v_{r_{\hub}}'$ in place of
$v_{\crithub}$ and $v_{\crithub}'$ by Lemma \ref{Lcrithub2}.

By Corollary \ref{Cdiscrepancy}, since $(G'/G)H^p$ has infinite
$r_{\hub}$-discrepancy valuation, we have
$v_{r_{\hub}}(J/(G'/G)H^p  - 1) \geq \min(\sigma, \sigma_0)$. 
Thus 
\begin{equation}\label{Ecrithubdiscrepancy}
v_{r_{\hub}}'\left(J - \frac{G'}{G}H^p\right) = v_{r_{\hub}}'\left(\left(\frac{J}{(G'/G)H^p} - 1\right) \frac{G'}{G}H^p\right) \geq \min(\sigma, \sigma_0 + \beta).
\end{equation}
Now, replace all terms of $J$ of degree not congruent to $-1 \pmod{m}$ or $0 \pmod{p}$ in $T^{-1}$ with the corresponding terms of $(G'/G)H^p$.   
Since $\K_1((G'/G)H^p)$ is $\psi$-equivariant, \eqref{Ecrithubdiscrepancy} shows that our new $J$ has $r_{\hub}$-discrepancy valuation 
at least $\min(\sigma, \sigma_0 + \beta)$ and lies in $\KK$.   By
\eqref{E6}, we still have $v_{\crithub}(J-1) \geq \beta$.  Repeating this process, we eventually obtain 
$J$ with $r_{\hub}$-discrepancy valuation at least $\sigma$.  This is the $J'$ we seek.
\Endproof

\begin{lem}\label{Llessdiscrepancy}
Let $0 < \sigma < p/(p-1)$.  Let $J \in 1 +
T_{\crit}^{-1}\m\{T_{\crit}^{-1}\}$.  After a possible
finite extension of $K$ and crit-negligible adjustment to $J$,
there exists $J' \in 1 + T^{-1}\m\{T^{-1}\} \cap \KK$ such 
$J' =' J$, and $J'$ has $r_{\crit}$-discrepancy valuation at least $\sigma$. Furthermore,
$v_{r_{\crit}}(J'-1) > 0$.
\end{lem}

\proof
The proof is the same as that of Lemma \ref{Lcrithublessdiscrepancy},
replacing Lemma \ref{Lcrithubimprovelem3} by Lemma \ref{Limprovelem2},
$v_{\crithub}$, $v_{r_{\hub}}$, $v_{\crithub}'$, and $v_{r_{\hub}}'$
by $v_{r_{\crit}}$, $v_{r_{\crit}}$, $v_{r_{\crit}}'$, and
$v_{r_{\crit}}'$, respectively, ``hub-negligible'' by ``crit-negligible,'' 
and omitting all mentions of $G_{\hub}$ and $G_{\hub}'$.
\Endproof

We recall the main proposition to be proved:

\begin{prop}[Proposition \ref{Pcrithubreduce}]\label{Pcrithubreduce2} 
Suppose $N_1 = (p-1)u_{n-1} - m$ (this is consistent with Assumption
\ref{AN1N2}).
Let $G_{\crit}, G_{\hub} \in G_{\crit, g}, \G_{\hub, \alpha}$, respectively. 
Let $r \in [0, r_{\hub}) \cap \QQ$, and let $f \in t^{1-m}k[t^{-m}]$ have degree less than $u_n$ in 
$t^{-1}$, which we regard as the reduction of $T_r$ in $\kappa_r$ (\S\ref{Sgeomsetup}).  Assume $f$ has no terms of degree divisible by $p$.  
Let $\beta = p/(p-1) - u_nr$.  After a possible finite extension of $K$, 
there exist $G_{\crit}', G_{\hub}' \in \G_{\crit, g}, \G_{\hub, \alpha}$
respectively, and $F \in \KK$ with $v_r(F) = 0$ and $[F]_r = f$ 
such that 
$$\frac{G_{\crit}'G_{\hub}'}{G_{\crit}G_{\hub}} \equiv 1 - p^{\beta}F \pmod{(\KK^{\times})^{p}}.$$
\end{prop}

\proof 
We first remark that if $A \in \KK$ such that $v_r(A) = 0$, then
$[p^{\beta}A]_r$ is unaffected by hub-negligible adjustments to $A$.
Essentially, this is the reason for defining hub-negligible as we do.

Let $F'$ be a polynomial in $T^{-1}$ of the same degree as $f$ such that $v_r(F') = 0$, that $[F']_r = f$, and that
$F'$ has no terms of degree divisible by $p$.  
Now, $$v_{\crithub}(p^{\beta}F') \geq \min(\beta - \deg(f)(r_{\hub}
-r) - s,\ \beta - (N_1 + u_{n-1})(r_{\crit} - r)).$$ 
Since $\deg(f) < u_n$ and $r \geq 0$, one calculates that
$v_{\crithub}(p^{\beta}F')$ is positive.
Choose $\sigma$ such that $\beta - (r_{\hub} - r) < \sigma < p/(p-1)$.  
By Lemma \ref{Lcrithublessdiscrepancy}, after making an adjustment to
$F'$ resulting in a hub-negligible adjustment to $p^{\beta}F'$, there
exists $F'' \in p^{-\beta}T^{-1}\m\{T^{-1}\}$ such that 
$p^{\beta}F'' =' p^{\beta}F'$,
that $1 - p^{\beta}F''$ has $r_{\hub}$-discrepancy valuation at
least $\sigma$, and that $v_{\crithub}(1 - p^{\beta}F'') > 0$.

Now, since $1 - p^{\beta}F''$ has $r_{\hub}$-discrepancy valuation
at least $\sigma$, it has $r$-discrepancy valuation at least $\sigma +
(r - r_{\hub}) > \beta$.  By Lemma \ref{Ldiscrepancy} (noting that
$[p^{\beta}F']_r$, and thus $[p^{\beta}F'']_r$, has no terms of degree divisible by $p$), we have that
$[p^{\beta}F'']_r$ contains only terms of degree congruent to $-1
\pmod{m}$ in $t^{-1}$.  Since the same is true by construction for 
$[p^{\beta}F']_r$, we have $[F'']_r = [p^{\beta}F'']_r =
[p^{\beta}F']_r = f$.

Furthermore, Lemma \ref{Lcrithubimprovelem3}
yields $G_{\crit}', G_{\hub}' \in \G_{\crit,g}, \G_{\hub, \alpha}$, respectively, and $H \in 1+ T^{-1}\m[T^{-1}]$ such that  
after making an adjustment to $F''$ resulting in a hub-negligible
adjustment to $p^{\beta}F''$ (which does not change
$[p^{\beta}F'']_r = f$ or $v_r(p^{\beta}F'') = \beta$), we obtain
\begin{equation}\label{E7}
v_{\crithub}'\left(\frac{1 -
    p^{\beta}F''}{(G_{\crit}'G_{\hub}'/G_{\crit}G_{\hub})H^p} -
  1\right) \geq \sigma.
\end{equation}
By Lemma \ref{Lcrithub2}, \eqref{E7} also holds when
$v_{\crithub}'$ is replaced by $v_{r_{\hub}}'$.
Since the fraction in \eqref{E7} has $r_{\hub}$-discrepancy
valuation $\geq \sigma$, Corollary \ref{Cdiscrepancy} shows that 
$v_{r_{\hub}}'$ can even be replaced by $v_{r_{\hub}}$ in \eqref{E7}.  We conclude that
$$v_{r_{\hub}}\left(\frac{G_{\crit}'G_{\hub}'}{G_{\crit}G_{\hub}}H^p -  (1 - p^{\beta}F'')\right) \geq \sigma.$$

In particular, 
\begin{equation}\label{E8}
v_r \left(\frac{G_{\crit}'G_{\hub}'}{G_{\crit}G_{\hub}}H^p -
 (1 - p^{\beta}F'')\right) \geq \sigma + (r_{\hub} - r) > \beta.
\end{equation}

Finally, let $F$ be such that $(G_{\crit}'G_{\hub}'/G_{\crit}G_{\hub})
H^p = 1 - p^{\beta}F$. 
Since $[F'']_r = f$, we need only show that $v_r(F) = 0$ and $[F]_r = [F'']_r$.  This follows from
\eqref{E8} and the fact that $v_r(p^{\beta}F'') = \beta$.
\Endproof 

\begin{rem}
Being able to replace $v_{r_{\hub}}'$ by $v_{r_{\hub}}$ in \eqref{E7}
in order to clear denominators is the
essential reason why we need the concept of discrepancy valuation.
\end{rem}

Proposition \ref{Preduce}, which we recall below, now follows easily.

\begin{prop}[Proposition \ref{Preduce}]
Let $G \in \G_{\crit, g}$, let $r \in [r_{\hub}, r_{\crit}) \cap \rats$, and let $f \in t^{1-m}k[t^{-m}]$ be a polynomial of degree less than $N_1 + u_{n-1}$ in $t^{-1}$, which we regard as the reduction of $T_r$ in $\kappa_r$ (\S\ref{Sgeomsetup}).  Assume $f$ has no terms of degree divisible by $p$.  
Let $\beta = (N_1 + u_{n-1})(r_{\crit} - r)$.  After a possible finite extension of $K$, 
there exist $G' \in \G_{\crit, g}$ and $F \in \KK$ with $v_r(F) = 0$ and $[F]_r = f$ 
such that 
$$\frac{G'}{G} \equiv 1 - p^{\beta}F \pmod{(\KK^{\times})^{p}}.$$
\end{prop}

\proof
The proof is the same as that of Proposition \ref{Pcrithubreduce2},
replacing Lemma \ref{Lcrithublessdiscrepancy} by Lemma
\ref{Llessdiscrepancy},
Lemma \ref{Lcrithubimprovelem3} by Lemma \ref{Limprovelem2},
$v_{\crithub}$, $v_{r_{\hub}}$, $v_{\crithub}'$, and $v_{r_{\hub}}'$
by $v_{r_{\crit}}$, $v_{r_{\crit}}$, $v_{r_{\crit}}'$, and
$v_{r_{\crit}}'$, respectively, ``hub-negligible'' by ``crit-negligible,'' 
choosing $\beta - (r_{\crit} - r) < \sigma < p/(p-1)$, and omitting all mentions of $G_{\hub}$ and $G_{\hub}'$.
\Endproof

\bibliographystyle{alpha}
\bibliography{main}

\end{document}